\documentclass[a4paper,
11pt,  reqno]{amsart}
\usepackage[margin=1.2in,marginparwidth=1.5cm, marginparsep=0.5cm]{geometry}

\usepackage{amsmath,amssymb,amscd,amsthm,amsxtra, esint}
\usepackage{mathrsfs}
\usepackage[most]{tcolorbox}
\usepackage{color}
\usepackage{hyperref}
\usepackage{subcaption}
\usepackage{float}
\usepackage[normalem]{ulem}

\definecolor{gr}{rgb}   {0.,   0.69,   0.23 }
\definecolor{bl}{rgb}   {0.,   0.5,   1. }
\definecolor{mg}{rgb}   {0.85,  0.,    0.85}
\definecolor{yl}{rgb}   {0.8,  0.7,   0.}
\definecolor{or}{rgb}  {0.7,0.2,0.2}

\newcommand{\Gr}{\color{gr}}

\IfFileExists{bbm.sty}{\usepackage{bbm}}{%
  \newcommand{\mathbbm}[1]{\mathchoice%
    {\hbox{\sffamily 1\kern-.40em 1}}%
    {\hbox{\sffamily 1\kern-.40em 1}}%
    {\hbox{\scriptsize\sffamily 1\kern-.30em 1}}%
    {\hbox{\tiny\sffamily 1\kern-.25em 1}}}%
}
\usepackage{forest}

\newcommand{\jap}[1]{\langle #1 \rangle}
\newcommand{\I}{\hspace{0.2mm}\textup{I}\hspace{0.2mm}}
\newcommand{\II}{\textup{I \hspace{-2.9mm} I}}
\newcommand{\III}{\textup{I \hspace{-2.9mm} I \hspace{-2.9mm} I}}
\newcommand{\IV}{\textup{I \hspace{-2.9mm} V}}
\newcommand{\5}{\textup{V}}

\newcommand{\noi}{\noindent}

\newcommand{\ft}{\widehat}
\newcommand{\ind}{\mathbbm 1}
\newcommand{\les}{\lesssim}
\newcommand{\dl}{\delta}
\newcommand{\F}{\mathcal{F}}
\newcommand{\eps}{\varepsilon}

\newcommand{\wt}{\widetilde}

\newcommand{\dx}{\partial_x}
\newcommand{\dt}{\partial_t}

\allowdisplaybreaks[2]

\newtheorem{theorem}{Theorem} [section]

\newtheorem{lemma}[theorem]{Lemma}
\newtheorem{proposition}[theorem]{Proposition}

\newtheorem{definition}[theorem]{Definition}

\theoremstyle{remark}
\newtheorem{remark}[theorem]{Remark}

\DeclareMathOperator{\Id}{\rm Id}

\usepackage{mathtools}
\mathtoolsset{showonlyrefs=true}

\newcommand{\R}{\mathbb{R}}
\newcommand{\CC}{\mathbb{C}}
\newcommand{\T}{\mathbb{T}}
\newcommand{\bul}{\bullet}

\let\Re=\undefined\DeclareMathOperator*{\Re}{Re}
\let\Im=\undefined\DeclareMathOperator*{\Im}{Im}

\let\P= \undefined
\newcommand{\P}{\mathbf{P}}

\newcommand{\TT}{\mathcal{T}}
\newcommand{\Tr}{\mathrm{T}}

\newcommand{\al}{\alpha}
\newcommand{\be}{\beta}

\newcommand{\Dl}{\Delta}

\newcommand{\G}{\Gamma}

\newcommand{\ld}{\lambda}

\newcommand{\s}{\sigma}

\newcommand{\bfrak}{\mathfrak{b}}
\newcommand{\gfrak}{\mathfrak{g}}

\newcommand{\pb}{\mathbf{p}}
\newcommand{\gb}{\mathbf{g}}
\newcommand{\sbf}{\mathbf{s}}

\newcommand{\Ft}{{\mathcal{F}}}

\newcommand{\ta}{\theta}

\renewcommand{\l}{\ell}

\newcommand{\ges}{\gtrsim}

\newcommand{\jb}[1]
{\langle #1 \rangle}

\renewcommand{\S}{\mathcal{S}}

\newcommand{\ff}{\mathfrak{f}}

\newcommand{\M}{\mathcal{M}}

\newcommand{\N}{\mathbb{N}}

\DeclareMathOperator{\mul}{mul}
\DeclareMathOperator{\mulB}{mul_{\BB}}

\tikzset{
	dotA/.style={
		transform shape, fill,circle,inner sep=1.5pt, label distance=-1pt,
		font={\normalsize }
	},
	>=stealth,
}

\tikzset{
	dotAB/.style={
		transform shape, fill,circle,inner sep=1.5pt, label distance=-10pt,
		font={\small }
	},
	>=stealth,
}

\makeatletter
\def\DeclareSymbolNew#1#2#3#4#5{\expandafter\gdef\csname MH@symb@#1\endcsname{\tikz[baseline=#2,scale=0.5, every node/.style={dotA},
		level distance=15pt,
		level 1/.style={sibling distance=#3pt},
		level 2/.style={sibling distance=#4pt},
		level 3/.style={sibling distance=15pt}
		]{#5}}}
\def\<#1>{\csname MH@symb@#1\endcsname}
\makeatother

\def\distC{15}

\DeclareSymbolNew{1}{-7}{\distC}{\distC}
{
	\node   (root) {}
	child
	child
	;
	\node  at (root-1) {};
	\node at (root-2) {};
}

\DeclareSymbolNew{11}{-7}{\distC}{\distC}
{
	\node   (root) {}
	child
	child
	;
	\node  at (root-1) {};
	\node at (root-2) {};
}

\newcommand{\cherryS}{
	\resizebox{!}{1.3ex}{%
		\<1>%
	}
}

\newcommand{\cherry}[3]{
	\begin{tikzpicture}[
		baseline=-8, scale=0.6, every node/.style={dotA},
		level distance=15pt,
		level 1/.style={sibling distance=\distC pt}
		]
		\node  [label=above:{$#1$}]  (root) {}
		child
		child
		;
		\node [label= left:{$#2$}]  at (root-1) {};
		\node [label= right:{$#3$}] at (root-2) {};
	\end{tikzpicture}
}

\DeclareSymbolNew{01}{-7}{\distC}{\distC}
{
	\node  [circle,draw=black, fill=white, inner sep=0pt,minimum size=7pt]  (root) {}
	child  child;
	\node  at (root-1) {};
	\node at (root-2) {};
}

\DeclareSymbolNew{21}{-12}{\distC}{\distC}
{
	\node   (root) {}
	child {child child} child;
	\node  at (root-1) {};
	\node at (root-2) {};
	\node  at (root-1-1) {};
	\node at (root-1-2) {};
}

\DeclareSymbolNew{22}{-12}{\distC}{\distC}
{
	\node   (root) {}
	child child {child child};
	\node  at (root-1) {};
	\node at (root-2) {};
	\node  at (root-2-1) {};
	\node at (root-2-2) {};
}

\DeclareSymbolNew{31}{-12}{\distC}{\distC}
{
	\node   (root) {}
	child {child {child child} child} child;
	\node  at (root-1) {}; \node at (root-2) {};
	\node at (root-1-1) {}; \node at (root-1-2) {};
	\node  at (root-1-1-1) {}; \node at (root-1-1-2) {};
}

\DeclareSymbolNew{32}{-12}{\distC}{\distC}
{
	\node   (root) {}
	child {child child {child child} } child;
	\node  at (root-1) {}; \node at (root-2) {};
	\node at (root-1-1) {}; \node at (root-1-2) {};
	\node  at (root-1-2-1) {}; \node at (root-1-2-2) {};
}

\DeclareSymbolNew{33}{-12}{21}{\distC}
{
	\node   (root) {}
	child {child child} child {child child};
	\node  at (root-1) {}; \node at (root-2) {};
	\node at (root-1-1) {}; \node at (root-1-2) {};
	\node  at (root-2-1) {}; \node at (root-2-2) {};
}

\DeclareSymbolNew{34}{-12}{\distC}{\distC}
{
	\node   (root) {}
	child child {child {child child} child  } ;
	\node  at (root-1) {}; \node at (root-2) {};
	\node at (root-2-1) {}; \node at (root-2-2) {};
	\node  at (root-2-1-1) {}; \node at (root-2-1-2) {};
}

\DeclareSymbolNew{35}{-12}{\distC}{\distC}
{
	\node   (root) {}
	child child {child child {child child} } ;
	\node  at (root-1) {}; \node at (root-2) {};
	\node at (root-2-1) {}; \node at (root-2-2) {};
	\node  at (root-2-2-1) {}; \node at (root-2-2-2) {};
}

\newcommand{\Ical}{\mathcal{I}}

\newcommand{\Gc}{\mathcal{G}}

\numberwithin{equation}{section}
\numberwithin{theorem}{section}

\newcommand{\MM}{\mathrm{M}}
\newcommand{\NN}{\mathrm{N}}
\newcommand{\Ncal}{\mathcal{N}}
\newcommand{\BB}{\mathrm{B}}
\newcommand{\D}{\mathrm{D}}
\newcommand{\DD}{\mathcal{D}}

\newcommand{\simao}[1]{}
\newcommand{\joao}[1]{}

\makeatletter
\@namedef{subjclassname@2020}{%
	\textup{2020} Mathematics Subject Classification}
\makeatother

\title{Local well-posedness for the Novikov-Veselov equation at subcritical regularity}

\title[LWP for the Novikov-Veselov equation at subcritical regularity]{Local well-posedness for the Novikov-Veselov equation at subcritical regularity}
\author[A.~Chapouto \and S.~Correia \and J.~P.~Ramos]
{Andreia Chapouto \and Sim\~ao Correia \and Jo\~ao Pedro Ramos}

\date{\today}

\address{Andreia Chapouto,
	CNRS, Laboratoire de math\'ematiques de Versailles, UVSQ, Universit\'e Paris-Saclay, CNRS, 45 avenue des \'Etats-Unis, 78035 Versailles Cedex, France, and Monash University, Clayton, VIC 3800, Australia}

\email{andreia.chapouto@monash.edu}

\address{Sim\~ao Correia,
	Center for Mathematical Analysis, Geometry and Dynamical Systems,
	Department of Mathematics,
	Instituto Superior T\'ecnico, Universidade de Lisboa
	Av. Rovisco Pais, 1049-001 Lisboa, Portugal
}

\email{simao.f.correia@tecnico.ulisboa.pt}

\address{Jo\~ao Pedro Ramos, 
	IMPA, Rio de Janeiro, Brazil}

\email{joao.ramos@impa.br}

\subjclass[2020]{35A01, 35A22, 35Q53, 37K10}

\keywords{Novikov-Veselov equation, gauge transform,
	local well-posedness.
}
\begin{document}

\begin{abstract}

We establish local well-posedness of the Novikov--Veselov equation (NV), at every fixed value of the energy parameter, for initial data in $H^s(\R^2)$ in the entire subcritical Sobolev range $s>-1$. 
The argument is perturbative and does not invoke the integrable structure, instead adapting the strategy in our previous work on KdV (2026). An infinite normal-form expansion indexed by binary trees is combined with a gauge transform that bypasses the obstructions in the associate bilinear estimate, and the resulting gauged equation -- carrying nonlinear contributions of every order -- is closed by a contraction in a Fourier restriction space. The analytic engine is a family of frequency-restricted multilinear bounds, supported by a detailed description of the two-dimensional resonance geometry attached to (NV).

\end{abstract}

\maketitle

\section{Introduction}

\subsection{Setting of the problem and main results}
We consider the Novikov-Veselov equation on $\R^2$:
\begin{equation}\tag{NV}\label{eq:NV}
	\partial_t u + (\partial^3 + \overline{\partial}^3)u + \tfrac34\big(\partial(u\overline{\partial}^{-1}\partial u) 
    + \overline{\partial}(u\partial^{-1}\overline{\partial}u)\big) 
    + E(\overline{\partial}^{-1}\partial^2 u + \partial^{-1}\overline{\partial}^2 u )=0,
\end{equation}
where $u$ is real-valued, $\partial=\frac12(\partial_x+i\partial_y)$ and $\overline{\partial}=\frac12(\dx-i\partial_y)$ are the Wirtinger derivatives, and $E\in\R$ denotes the 
\emph{energy} parameter.

The Novikov-Veselov equation appeared in the work of Novikov and Veselov \cite{NovikovVeselov} (and implicitly, of Manakov
\cite{Manakov}) 
associated with the
two-dimensional inverse scattering problem for the Schr\"odinger operator, at fixed energy $E$. 
This equation is of great analytical interest as it is closely related to various canonical dispersive models; see \cite{CMMPSS15} for a survey. 
Firstly, $y$-independent solutions to \eqref{eq:NV} solve the Korteweg-de Vries equation (KdV), thus placing \eqref{eq:NV}
among the standard higher-dimensional
dispersive generalizations of KdV, together with the Kadomtsev-Petviashvili equation (KP) and the
Zakharov-Kuznetsov equation (ZK); see, for instance,
\cite{bourg_kpii,GuoPengWang,  HHK, IonescuKenigTataru,
Kinoshita_ZK, takaokaTzvetkov} for representative results for those models.
Moreover, after anisotropic rescalings, \eqref{eq:NV} formally converges to KP-I and KP-II as $E\to +\infty$ and $E \to -\infty$, respectively \cite{Grinevich00, KK17, kazeykinamunoz1, kazeykinamunoz2}.
In comparison {with KP and ZK}, the nonlocal structure of 
\eqref{eq:NV} and the rather involved geometry of its two-dimensional resonances give the equation a distinct analytic flavor.

A natural starting point  to understand the well-posedness for \eqref{eq:NV} is its scaling heuristics.
Given $\ld>0$, define 
\begin{align}\label{scaling}
(S_\ld u) (t,x, y) = \ld^2 u (\ld^3 t , \ld x, \ld y). 
\end{align}
Then, if $u$ is a solution to 
\eqref{eq:NV} with energy $E\in\R$,  $S_\ld u$ solves \eqref{eq:NV} with energy $\ld^2 E$. 
When $E=0$, \eqref{scaling} is a genuine symmetry of the equation, which induces the scaling critical exponent $s_c = -1$, due to the invariance of the homogeneous Sobolev norm $\dot{H}^{-1}(\R^2)$ under \eqref{scaling}. For $E\ne 0$, invariance of \eqref{eq:NV} under scaling holds only within the family of equations with energy parameters $E\in\R$. Nevertheless, the scaling heuristics extend to the whole family and we expect well-posedness of 
\eqref{eq:NV} to hold in $H^s(\R^2)$ for $s > s_c =-1$, regardless of the value of $E\in\R$.

In studying \eqref{eq:NV}, there are two natural viewpoints. One is via complete integrability: inverse scattering,
Miura-type transforms, and the connection to the modified Novikov-Veselov hierarchy can give existence and long-time results by exploiting the algebraic structure of the equation in a decisive way.
The other is dispersive/perturbative, where \eqref{eq:NV} is treated as a semilinear equation 
amenable to robust Fourier-analytic methods. 
In this work, we show that a perturbative approach suffices to recover the full subcritical regime for local well-posedness.

On the integrable side, the construction of
solutions through inverse scattering and related techniques goes back to
\cite{BLMP,MusicPerry, NPT_mNV,  Perry, Tsai}. The more recent results are for the zero-energy case ($E=0$), where the Miura map
\cite{Bogdanov}
\begin{equation}\label{miura}
\mathbf{M}[v]= 2\partial v + |v|^2,
\end{equation}
takes solutions of the \emph{modified} Novikov-Veselov equation
\begin{equation}
\partial_t v + (\partial^3+\overline{\partial}^3) v + \tfrac34\big(\partial(v\overline{\partial}^{-1}\partial |v|^2) + \overline{\partial}(v\partial^{-1}\overline{\partial}|v|^2) + v\partial^{-1}\overline{\partial}(\overline{v}\overline{\partial}v) + v\overline{\partial}^{-1}\partial(\overline{v}\partial v)\big)=0\label{eq:mNV}\tag{mNV}
\end{equation}
satisfying $\partial v=\overline{\partial} v$, to solutions of \eqref{eq:NV}. 
Since \eqref{eq:mNV} has cubic nonlinearity, its low-regularity theory is better understood:
Schottdorf \cite{Schottdorf} first showed small-data global well-posedness in the critical space $L^2(\R^2)$ for mNV with a simpler nonlinearity, which was extended to large-data and the full model by Nachman, Perry and Tataru
\cite{NPT_mNV}. Then, using the Miura map in \eqref{miura}, \cite{NPT_mNV} showed  global well-posedness of \eqref{eq:NV} (with $E=0$) for soliton free data $u_0\in\dot{H}^{-1}(\R^2)+L^1(\R^2)$, namely such that $-\Dl + u_0 \ge 0$.
A parallel mechanism is well known for KdV \cite{BuckmasterKoch, CKSTT_sharp_kdv,   KPST_miura, KT06,   KillipVisan_kdv, Miura}, but the analogy is only partial: for KdV, complete integrability provides coercive conservation laws, and these eventually yield a global theory down to $H^{-1}(\R)$ \cite{KillipVisan_kdv}, which is sharp due to \cite{Molinet_ill_line}. For \eqref{eq:NV}, the conserved quantities are not coercive, in accordance with the existence of blow-up solutions \cite{TaimanovTsarev}. 
Thus, while the integrable structure is extremely informative, it does not seem sufficient to settle the semilinear
low-regularity Cauchy problem for \eqref{eq:NV}.

There have also been recent perturbative developments in the analysis of \eqref{eq:NV} \cite{adamsgrunrock,angelopolous, kazeykinamunoz1, kazeykinamunoz2}. 
Local well-posedness of \eqref{eq:NV} in $H^s(\R^2)$ for $s>\frac12$ was shown for $E\ne0$ by Kazeykina and Mu\~noz \cite{kazeykinamunoz1, kazeykinamunoz2}, via refined stationary-phase estimates adapted to the geometry of the dispersive symbol, and for $E=0$ by Angelopoulos \cite{angelopolous}, using Bourgain-type spaces. 
The latter was improved to $s>-\frac34$ by Adams and Gr\"unrock \cite{adamsgrunrock}.
All of these results (including those from integrability) leave a gap in the subcritical regime, namely
\begin{align}
-1< s\le \tfrac12 \quad \text{for } E\neq 0,
\qquad
-1< s\le -\tfrac34 \quad \text{for } E=0, 
\end{align}
whose threshold is known to be optimal for fixed-point arguments (excluding the endpoint $s=-1$) due to the failure of analyticity of the data-to-solution map for $s<-1$ and $E=0$ in \cite{angelopolous}.
Our main result closes the above gap completely.

\begin{theorem}\label{thm:main_intro}
Fix $E\in\R$ and $s>-1$. For every $u_0\in H^s(\R^2)$ there exist a time $T=T(\|u_0\|_{H^s}, |E|)>0$ and a function $u\in C([0,T];H^s(\R^2))$ with the following property: for any sequence $u_{0,n}\in H^\infty(\R^2)$ with $u_{0,n}\to u_0$ in $H^s(\R^2)$, the corresponding smooth solutions $u_n$ to \eqref{eq:NV} -- provided by \cite{angelopolous, kazeykinamunoz1, kazeykinamunoz2} -- are defined on $[0,T]$ and satisfy $u_n\to u$ in $C([0,T];H^s(\R^2))$. The resulting data-to-solution map $u_0\in H^s(\R^2)\mapsto u\in C([0,T];H^s(\R^2))$ is analytic.
\end{theorem}

Theorem~\ref{thm:main_intro} is purely perturbative, not relying on the integrable nature of \eqref{eq:NV}. This extends the method of \cite{chcora} for KdV to the $2$-$d$ setting, combining infinite normal form reductions, a gauge transformation, and $X^{s,b}$-analysis (see \eqref{xsb}). 
We include some additional remarks before proceeding to a description of our method. 

\begin{remark}
A direct fixed point argument in Fourier restriction norm spaces $X^{s,b}$ suffices to extend the result of \cite{adamsgrunrock} to any energy level $E$ for $s>-\frac34$ (Proposition~\ref{prop:bilin_s-34}), which is already new for the non-zero energy regime, improving on \cite{kazeykinamunoz1, kazeykinamunoz2}. However, this is insufficient to cover the full subcritical regime, as the crucial bilinear estimate fails for $s\le - \frac78$ (Proposition~\ref{PRO:ill}). 
\end{remark}

\begin{remark}\label{REM:thm}
(i) 

Theorem \ref{thm:main_intro} joins a short list of quadratic two-dimensional dispersive models whose subcritical Sobolev theory was established by perturbative methods alone, together with the nonlinear Schr\"odinger equation (NLS) with $u^2$ nonlinearity \cite{BDS08} and KP-II \cite{hadac}. 
For KP-II, a fixed point argument in $X^{s,b}$ suffices, due to its favorable resonance function. For quadratic NLS, \cite{BDS08} requires $X^{s,b}$-spaces adapted to the resonant interactions.
For \eqref{eq:NV}, which has a non-trivial resonant set $\{\Phi =0 \}$, the direct approach fails for $s\le - \frac78$ (Proposition~\ref{PRO:ill}); instead, Theorem~\ref{thm:main_intro} exploits further structure of \eqref{eq:NV} through normal form reductions.

\smallskip 

\noi
(ii) Theorem~\ref{thm:main_intro} leaves open the critical endpoint $s=-1$.

Critical well-posedness has been achieved for several \emph{cubic} $2$-$d$ dispersive models using the $U^2$--$V^2$ machinery of Koch--Tataru \cite{HHK, KochTataru}, such as \eqref{eq:mNV} \cite{NPT_mNV, Schottdorf}, the Dysthe equation \cite{MPS_Dysthe}, and the modified Zakharov-Kuznetsov equation \cite{CorreiaKinoshita}.

In our approach, the obstruction is explicit:
each of the nonlinear terms in \eqref{gnv} presents critical decay in frequency, leading to arbitrarily many \emph{simultaneous} logarithmic divergences. Moreover, the gauge breaks down for $s<-1$ (Lemma~\ref{LEM:D}~(ii)).

In the zero-energy case, \cite{NPT_mNV} obtained global well-posedness of \eqref{eq:NV}, with a continuous solution map, for soliton-free data in $\dot{H}^{-1}(\R^2) + L^1(\R^2)$, namely the range of the Miura map $\mathbf{M}$ \eqref{miura}. However, a global result cannot extend for all critical data, in view of the finite-time blow-up solutions of \cite{TaimanovTsarev}. 

The local well-posedness for general data at critical regularity remains an open problem.

\end{remark}

\begin{remark}\label{remark:E}
(i) Theorem~\ref{thm:main_intro} works for all $E\in\R$, without distinguishing its sign, in contrast with \cite{kazeykinamunoz1, kazeykinamunoz2}, where the case $E<0$ benefits from a more favorable resonance relation. In the proof, we first use the scaling \eqref{scaling} to reduce to small $|E|$, so that the resonance function $\Phi = \Phi_0 + E \Phi_E$ in \eqref{Phi} is dominated by the zero-energy part $\Phi_0$, for large frequencies. 
In \cite{kazeykinamunoz1, kazeykinamunoz2}, the low frequencies are responsible for the difficulties in their analysis, 
which we discuss in Remark~\ref{REM:low-freq} in the context of our method.

\smallskip
\noi (ii) As recalled earlier, after an anisotropic rescaling, \eqref{eq:NV} formally converges to KP-I and KP-II when $E\to +\infty$ and $E\to -\infty$, respectively; see \cite{Grinevich00, KK17, kazeykinamunoz1}.
The limiting equation as $E\to -\infty$ is KP-II, which is known to be analytically well-posed, as is \eqref{eq:NV} for all finite $E$ in Theorem~\ref{thm:main_intro}. 
In the other limit, $E\to+\infty$, \eqref{eq:NV} approaches KP-I which is quasilinear, so one should expect a singular limit. Recently, \cite{CHLO} evidenced a singular transition from a quasilinear to an analytic flow in the shallow-water limit from the intermediate long wave equation to KdV. Here, it would be of interest to study the reverse transition from a semilinear to a quasilinear flow in the \eqref{eq:NV} to KP-I limit.

\end{remark}

\begin{remark}
In \cite{adamsgrunrock}, Adams-Gr\"unrock also establish local well-posedness of \eqref{eq:NV} with $E=0$ for mean-zero initial data in $H^s(\T^2)$ for $s>-\frac15$, where $\T^2$ denotes the $2$-$d$ torus. 
One of their open questions is if the mean-zero condition can be removed, since, unlike KdV where mean-zero result immediately implies any prescribed mean result, for \eqref{eq:NV} prescribing the mean is equivalent to solving \eqref{eq:NV} with non-zero energy. It would be of interest to see if our result could extend to the periodic setting to lead to an improvement in this direction. 
\end{remark}

\subsection{Overview of the method}\label{SUBSEC:method}

The method of proof of Theorem~\ref{thm:main_intro} is shared with our companion paper on KdV \cite{chcora}, but the outcome is qualitatively different: where the KdV analysis produces an improvement in the Fourier-Lebesgue scale, the same mechanism applied to \eqref{eq:NV} covers the full subcritical regime for the $L^2$-based \emph{Sobolev} theory.
Three ingredients are used in concert: infinite normal-form reductions, a gauge transform (motivated by the latter reductions), and multilinear estimates of Fourier-restriction type.
The normal form reductions are in the spirit of  Babin--Ilyin--Titi \cite{bit} and of the later developments
\cite{GKO13, Kishimoto22, KOY20, OW21}, while the analytic estimates rely on the frequency-restricted framework in \cite{COS}. We refer the reader to \cite[Section 1]{chcora} for a detailed description of the method.

We first describe the normal form reduction procedure for \eqref{eq:NV}. Given a solution $u$ to \eqref{eq:NV}, we consider (the spatial Fourier transform of) its interaction representation
\begin{align}\label{inter}
\ft{v}(t):=\F_{x}\big(U(-t)u(t)\big),
\end{align}
where $U(t)$ is the linear propagator associated with \eqref{eq:NV}. 
Then, the Duhamel formulation for $\ft{v}$ in \eqref{inter} can be written as\footnote{\label{foonote}Here and below, we drop powers of $2\pi$ coming from Fourier normalization and coefficients coming from \eqref{eq:NV}.}
\begin{equation}\label{eq:NV2}
		\ft v(t,\xi)=\ft v(0,\xi)+  i \int_0^t \!\!\!\int_{\xi=\xi_1+\xi_2} \!\!\! e^{it'\Phi(\xi,\xi_1,\xi_2)}K(\xi_1,\xi_2) \Phi_0(\xi,\xi_1,\xi_2) \, \ft v(t',\xi_1) \ft v(t',\xi_2)\,d\xi_1\,dt', 
\end{equation}
with frequencies $\xi,\xi_1,\xi_2 \in \R^2\simeq \CC$, $K(\xi_1,\xi_2) = \frac{\xi_1 \cdot \xi_2}{|\xi_1|^2 |\xi_2|^2}$, and resonance function composed of energy-free and energy components, respectively:
\begin{align}\label{Phi}
\Phi(\xi,\xi_1,\xi_2)& =\Phi_0(\xi,\xi_1,\xi_2)+E\,\Phi_E(\xi,\xi_1,\xi_2),\qquad \xi=\xi_1+\xi_2, 
\\\label{Phi0}
\Phi_0(\xi,\xi_1,\xi_2)
&
=\Re(-\xi^3+\xi_1^3+\xi_2^3)=-3\Re(\xi\xi_1\xi_2),
\\
\Phi_E(\xi,\xi_1,\xi_2)
&
=-\Re\Big(-\tfrac{\xi^3}{|\xi|^2}+\tfrac{\xi_1^3}{|\xi_1|^2}+\tfrac{\xi_2^3}{|\xi_2|^2}\Big).
\end{align}

\noi 
The main obstacle in estimating the nonlinear term in \eqref{eq:NV2} comes from the derivative loss encoded in the multiplier $K \Phi_0$, responsible for the regularity obstructions in the literature \cite{adamsgrunrock,angelopolous,kazeykinamunoz1,kazeykinamunoz2}. We hope to control this loss by exploiting the oscillations coming from the exponential factor $e^{i t \Phi}$ by means of the relation
\begin{align}\label{exp}
e^{it\Phi}=\frac{\partial_t e^{it\Phi}}{i\Phi}.
\end{align}

\noi 
Using \eqref{exp} and integrating-by-parts in time in \eqref{eq:NV2} naively would create artificial singularities near the resonant set $\{\Phi=0\}$, due to the singularities of $K$ and also the mismatch between $\Phi$ in the exponential and $\Phi_0$ in the multiplier. Instead, we follow a 2-step process: first, we partition the convolution region into a \textit{good part} and a \textit{bad part}, localized to $A_\dl$ (containing the set where $K$ is singular) and $A^c_\dl$ (containing the high$\times$high interactions, where the distance of frequencies to the origin is at least $1/\dl$ for a small parameter $\dl>0$), respectively; see \eqref{defA} for the definition. 
Secondly, we use the crucial identity $K\Phi_0 = K\Phi - EK \Phi_E$ from \eqref{Phi} in the bad part, obtaining a multiplier of the form
\begin{align}\label{intro-goodbad}
(K \Phi_0 \ind_{A_\dl} - E K \Phi_E \ind_{A^c_\dl}) + K\Phi \ind_{A^c_\dl}. 
\end{align}
The first contributions above are localized to the regions where their respective multipliers are under control and no further action is needed. For the last contribution, we have fixed the earlier mismatch and the localization guarantees boundedness (and indeed smallness) of $K$, allowing us to  perform the integration-by-parts trick.

The procedure is then recursive: splitting the region of integration, 
rewriting the multiplier as in \eqref{intro-goodbad}, and using the identity \eqref{exp} on the bad parts together with integration-by-parts in time produces, at each step, a \emph{boundary} contribution and a residual term (carrying the $\dt \ft v$ terms). Then, feeding the equation for $\dt\ft v$ (see \eqref{eq:NV2}) generates a new collection of terms,  
which can be systematically indexed by binary trees. 
The resulting terms with a time integral can be successfully estimated in Fourier restriction norm spaces $X^{s,b}$ for $s>-1$. However, the boundary terms are evaluated at fixed times and cannot be estimated in these spaces, 
creating the actual obstruction to closing a contraction argument.
As in \cite{chcora}, the entire series of boundary terms can be eliminated by a single nonlinear change of variables, the gauge transform $\Gc_\dl$ (see \eqref{eq:gauge} for an explicit form of its inverse). This leads to the gauged Novikov-Veselov equation
\begin{equation}\label{gnv}
\partial_t z + (\partial^3+\overline{\partial}^3)z + E(\overline{\partial}^{-1}\partial^2 + \partial^{-1}\overline{\partial}^2)z
= \sum_{n=2}^\infty 
N_n[z]
\tag{$\mathcal{G}_\delta$-NV}
\end{equation}
where $N_n[z]$ denotes a nonlinear term of homogeneity $n$ in $z$. We postpone a detailed description of the nonlinearity in \eqref{gnv} to Subsection~\ref{sec:nfe}; see in particular \eqref{z1}.

Equation \eqref{gnv} can be solved via 
a fixed-point argument in $X^{s,b}$-spaces (Proposition~\ref{prop:lwp_rnv}), which relies on 
 the  multilinear estimate
\begin{equation}\label{intro:nonl}
\Bigl\|\int_0^t U(t-t') N_n[z_1, 
\ldots, z_{n}] (t') \,dt'\Bigr\|_{X^{s,b}_T}
\le T^{0+} C^n \delta^{(s+1)(n-2)-3}\prod_{k=1}^{n}\|z_k\|_{X^{s,b}_T},
\end{equation} 
 for $s>-1$ and every $n\ge 2$ (Proposition~\ref{prop:bourgL2}). 
The decisive feature of this bound is the gain of a positive power of the gauge parameter $\dl$ (for $s>-1$ and $n\gg1$), ensuring summability of the nonlinear contributions by taking $\dl$ sufficiently small, depending only on the size of the data. 
This multilinear analysis is done systematically via the frequency-restricted-estimate machinery in \cite{COS},
a detailed study of the \eqref{eq:NV} resonance geometry and a systematic organization of tree
contributions according to their algebraic type.

To carry the well-posedness statement from \eqref{gnv} back to \eqref{eq:NV}, we rely on three ingredients. 
Firstly, the gauge $\Gc_\dl$ and its inverse are bi-Lipschitz, in fact analytic in $H^{s}(\R^2)$ for $s>-1$ (Lemma~\ref{LEM:D}). 
Secondly, the smooth solutions $u$ to \eqref{eq:NV} of \cite{angelopolous,kazeykinamunoz1,kazeykinamunoz2} agree with the smooth solutions $\Gc_\dl^{-1}[z]$, where $z$ solves \eqref{gnv} with matching gauged data. Lastly, combining the two previous points, we are then able ot extend the \eqref{eq:NV} flow by conjugating the gauge transform with the \eqref{gnv} solution map.

\bigskip

We complete this introduction with some technical remarks. 

\begin{remark}\label{REM:low-freq}
In \cite{kazeykinamunoz1, kazeykinamunoz2}, one of the main difficulties came from the energy component of the dispersion $\Phi_E$ dominating at low frequencies and singular at the origin. In our argument, $\Phi_E$ appears  explicitly in the multiplier localized to the high$\times$high regime (see the second term in \eqref{intro-goodbad}), where it can be regarded as a small perturbation, since 
$| K \Phi_0 \ind_{A_\dl}| \les |\xi|$ and $|E K \Phi_E \ind_{A^c_\dl}| \les |E| \dl^2 |\xi|$. Moreover, the arguments used in showing \eqref{intro:nonl} rely only on the geometric properties of sublevel sets of the full resonance \eqref{Phi} away from the $\Phi_E$ singular set. 
Compared to the KdV analysis in \cite{chcora}, these contributions are new and an intrinsic feature of \eqref{eq:NV}. 

\end{remark}

\begin{remark}\rm \label{REM:small}
The summability of the series of nonlinearities in \eqref{gnv} requires smallness. In \cite{chcora}, this was done via a reduction to small data through \textit{scaling}. Here, smallness is encoded directly in the gauged equation \eqref{gnv} through the parameter $\dl$ defining the regions $A_\dl, A_\dl^c$ (see \eqref{defA}), leading to the gain of powers of $\dl$ in \eqref{intro:nonl}. Although we use scaling \eqref{scaling} to reduce to small $|E|$ (see Remark~\ref{remark:E}~(i)), we do not require a reduction to small data in the proof of Theorem~\ref{thm:main_intro}. Therefore, the argument as presented here carries over to equations without a scaling symmetry.

\end{remark}

\subsection{Organization and notation}

The remainder of the paper is structured as follows. Section~\ref{sec:prem} sets up the tree formalism underlying the infinite normal-form expansion, derives the normal-form equation for \eqref{eq:NV}, and constructs the gauge transform $\Gc_\delta$, culminating in the gauged equation \eqref{z1}. Section~\ref{sec:multi} recalls the Fourier restriction space $X^{s,b}$ and the frequency-restricted framework, and proves the  multilinear bound in Proposition~\ref{prop:bourgL2}.
Section~\ref{sec:lwp} develops the local theory for the gauged equation (Proposition~\ref{prop:lwp_rnv}), proves the relevant properties of the gauge map (Lemma~\ref{LEM:D}), the equivalence of \eqref{eq:NV} and \eqref{z1} for smooth solutions (Proposition~\ref{PRO:equiv}), and   Theorem~\ref{thm:main_intro}. Appendix~\ref{sec:s-34} handles
the direct fixed-point argument in the regime $s>-\tfrac34$ and shows failure of the relevant bilinear estimate for $s\le -\frac78$.

\bigskip

We conclude with some notation used throughout the paper. The plane $\R^2$ is identified with $\CC$ via $\xi=(\eta,\zeta)\simeq\eta+i\zeta$.
For $A,B\in\R$, the symbol $A\lesssim B$ means $A\le C B$ for some $C>0$, while $A\sim B$ stands for the two-sided bound $A\lesssim B\lesssim A$. If $C<10^{-10}$, we write $A\ll B$.
We also use $A \lor B = \max(A,B)$ and $A\land B = \min(A,B)$. 
For $\xi_1,\xi_2 \in \CC^2$, the slightly stronger comparability $\max(|\xi_1|,|\xi_2|)\gg|\xi_1-\xi_2|$ is denoted by $\xi_1\simeq \xi_2$. 
We write $\alpha+$ and $\alpha-$ for $\alpha+\varepsilon$ and $\alpha-\varepsilon$, respectively, with $\varepsilon>0$ chosen sufficiently small. We use the Japanese bracket $\jap{\xi}:=(1+|\xi|^2)^{1/2}$, and we denote by $\mathcal{F}$ and $\mathcal{F}^{-1}$ (or $\widehat{\,\cdot\,}$ and ${\,\cdot\,}^{\vee}$) the spatial Fourier transform and its inverse on $\R^2$. 
Lastly, we use $H^s(\R^2)$ for the usual $L^2$-based Sobolev space and, given $T>0$, we use the shorthand notation $L^\infty_T H^s_x$ for $L^\infty([0,T];H^s(\R^2))$.

\subsection*{Acknowledgements}

The authors are deeply indebted to Axel Gr\"unrock and Jorge Silva for his encouragement and for several important discussions. A.~C. is thankful for the kind hospitality of Instituto Superior T\'ecnico, where part of this work was developed.

\smallskip

A.~C. was partially supported by CNRS--INSMI through a grant ``PEPS Jeunes chercheurs et jeunes chercheuses 2025''. S.~C. was partially supported by Funda\c{c}\~ao para a Ci\^encia e Tecnologia (FCT) through CAMGSD (project UID/04459/2025). A.~C. and S.~C. were also supported by FCT through the grant FCT/Mobility/1330917020/2024-25.  S.~C. and J.~P.~R. were also supported by FCT through project SHADE (project 2023.17881.ICDT, DOI 10.54499/2023.17881.ICDT).

\subsection*{AI statement} 
The authors declare that generative AI was not used in the research or preparation of this manuscript.

\section{Preliminaries}\label{sec:prem}

\subsection{Tree formalism}	
\label{SUB:tree}
We recall some definitions pertaining to rooted binary trees from our earlier paper \cite{chcora}.

\begin{definition}
\label{DEF:tree}
\rm 
A \emph{binary tree} $\TT$ is a partially ordered set (with partial order $\le$) satisfying the following properties:

\begin{enumerate}
\item[(i)] $\TT$ has a unique minimal element called the \textit{root node} labelled by $\emptyset$ (i.e., $\emptyset \le a$ for all $a \in \TT$). 

\smallskip 
\item[(ii)] Given $a\in \TT\setminus \{ \emptyset \}$, 
there exists a unique $b \in \TT$ such that: $b\ne a$, $b \le a$, and $b \le c \le a$ implies $c\in\{a,b\}$. In that case, we say that $b$ is the \textit{parent} of $a$ and $a$ is a \textit{child} of $b$. 

\smallskip 
\item[(iii)] A node $a\in \TT$ is called \textit{terminal} or a \textit{leaf} if it has no children. A \textit{non-terminal node} or a \textit{parent} $a\in\TT$ is a node with exactly two children, one called \textit{left child} and the other \textit{right child}.\footnote{Introducing the notation of left and right children of a node  reflects the lack of ambiguity in the usual planar representation of binary trees.}  
We write $\TT^0, \TT^\infty$ for the sets of parent and terminal nodes of $\TT$, respectively.  Note that $\TT$ is the disjoint union of $\TT^0$ and $\TT^\infty$. 

\smallskip 
\item[(iv)]
Each node in $\TT$ is labelled with a word in the following way: the root node has word~$\emptyset$; and for a given parent node $a\in \TT^0$, its left and right children have labels $a1$ and $a2$,  obtained by concatenating $1$ or $2$ to the right of $a$, respectively. Note that children of the root node have labels 1 and 2.

\end{enumerate}

\noi
We interchangeably use $\TT$ to denote the representation of the tree as a finite connected graph and to denote the set of words used to label its nodes. Note that the latter is also a partially ordered set, with words $a,b$ satisfying $a\le b$ if there exists a word $c$ such that $b = ac$.

\end{definition}

Given $j\in\N$, we use $\Tr_j$ to denote the family of all binary trees
with $j$ non-terminal nodes, which we call the $j$-th generation of trees. From Definition~\ref{DEF:tree}, it follows that, given $\TT\in\Tr_j$,
\begin{align}
| \TT^0 |  = j \quad \text{and} \quad |\TT^\infty| = j+1, 
\end{align}

\noi
while the number of distinct trees in $\Tr_j$ is given by the $j$-th Catalan number
\begin{align}\label{eq:catalan}
| \Tr_j | = \frac{1}{j+1} \binom{2j}{j}  = \frac{(2j)!}{(j+1)! j !}.
\end{align}

For readibility, we will reserve the symbols $\bul$, $\star$, $\circ$, and $\ast$ to denote nodes in a tree and their respective associated word. 
Next, we introduce some additional notation to identify different nodes in a tree.

\begin{definition}\rm
\label{DEF:tree2}
Let $\TT$ be a binary tree in the sense of Definition~\ref{DEF:tree} and 
$\bul, \star$ be two nodes in $\TT$. 

\begin{enumerate}

\smallskip
\item[(i)]
If $\bul \in \TT\setminus \{ \emptyset\}$, we write $\pb(\bul)$ to denote the \emph{parent node} of $\bul$. If $\pb(\bul)\neq \emptyset$, we write $\gb(\bul)=\pb(\pb(\bul))$ and call it the \textit{grandparent node} of $\bul$. 

\smallskip
\item[(ii)] 
If $\bul, \star\in \TT \setminus\{ \emptyset \}$ are distinct and share the same parent, they are called \textit{siblings}. We write $\sbf(\bul)$ for the sibling of $\bul$. 

\smallskip
\item[(iii)] 
Given $\bul\in \TT^0$, we say $\bul$ is a \textit{final parent} if  $\bul1, \bul2 \in \TT^\infty$. The set of final parents is denoted by $\TT^{0,\ff}$.

\end{enumerate}

\end{definition}

\begin{remark}\rm 
Note that, from Definition~\ref{DEF:tree}, the maps $\pb, \gb, \sbf$, for the parent, grandparent, and sibling of a given node, are well-defined in their respective domains, since if either the parent, grandparent, or sibling node exists, it is unique. 

\end{remark}

We next introduce concatenation and subtraction of trees.

\begin{definition} 
\label{DEF:tree3}
\rm

Let $\TT$ be a binary tree in the sense of Definition~\ref{DEF:tree}.

\begin{enumerate}

\smallskip
\item[(i)] (Concatenation)  
Given $ \star \in \TT^\infty$, we define
$$
\TT \overset{\star}{\oplus} \<1>
: = \TT \cup \{\star1, \star2\}
$$
as the result of adding the children $\star1,\star2$ of $\star$ to the tree $\TT$. 
Note that this operation now makes $\star$ a final parent in the tree $\TT \overset{\star}{\oplus} \<1>$.

\smallskip
\item[(ii)] (Subtraction) 
Given a final parent $\star  \in \TT^{0,\ff}$, we define the tree
$$
\TT \overset{\star }{\ominus } \<1>
: = \TT \setminus \{ \star1, \star2  \}
$$
as the result of eliminating the children $\star1, \star2$ of $\star$
from $\TT$. Note that $\star$ becomes a terminal node in the new tree $\TT \overset{\star}{\ominus } \<1>$. 
More generally, given $C\subseteq \TT^{0,f}$, we define
\begin{align}
\label{treeSubtr}
\TT \overset{C}{\ominus} \<1> : =
 \bigcap_{\star \in C} ( \TT \overset{\star}{\ominus}  \<1>)
= \bigcap_{\star\in C} (\TT \setminus \{\star1,\star2\})
= \TT \setminus \Big( \bigcup_{\star\in C} \{\star1,\star2\} \Big), 
\end{align}
which corresponds to the set where all the children of the final parents $\star\in C$ have been removed from the tree.

\end{enumerate}

\end{definition}

In order to move the complexity of the normal form reduction, performed in Subsection~\ref{sec:nfe},  to combinatorial aspects on binary trees, we need some additional definitions that will help us encode the precise structure of each nonlinear term.

\begin{definition}\rm 

Let $j\in\N$ and $\TT \in \Tr_j$. 

\begin{enumerate}
\item[(i)] (Index function) 
We define the index function $\vec\xi : \TT\to \CC^\TT$ such that $\vec\xi = (\xi_\bul)_{\bul\in\TT}$, where
\begin{align}
\label{parent}
\xi_\bul = \xi_{\bul1} + \xi_{\bul2}, \quad \bul\in \TT^0,
\end{align}
and $\bul1$ and $\bul2$ denote the children of $\bul$. Moreover, we define the set $\G(\TT)$ as the set of convolution relations on the tree $\TT$:
\begin{align}
\label{Gamma}
\G(\TT) : = \big\{ \vec\xi \in \CC^\TT: \ \xi_\bul = \xi_{\bul1} + \xi_{\bul2}, \ \bul \in \TT^0 \big\}. 
\end{align}

\smallskip

\item[(ii)] (Zero-energy phase function) 
We use $\Psi_0(\TT)$ to denote the \textit{zero-energy phase function} associated with the tree $\TT$, given by
\begin{align}
\label{psi0-tree}
\Psi_0(\TT)
:=\Re\bigg(- \xi^3 + \sum_{\bul\in \TT^\infty} \xi_\bul^3\bigg)
\end{align}
where $\TT^\infty$ denotes the set of leaves of $\TT$. 
Given a parent node $\star\in\TT^0$, we also introduce the partial zero-energy phase function for the tree $\cherry{\star}{}{}=\{\star, \star1, \star2\}$:
\begin{align}
\label{psi-circ}
\Psi_0(\star) = \Psi_0(\cherry{\star}{}{}) : = \Re \big( - \xi^3_\star + \xi^3_{\star1} + \xi^3_{\star2} \big). 
\end{align}
Note that the following additivity property holds for $\Psi$:
\begin{equation}\label{eq:additivity}
\Psi_0(\TT) = \sum_{\star \in \TT^0} \Psi_0(\star). 
\end{equation}
\item[(iii)] (Energy phase function) Given a binary tree $\TT$, we also define the \textit{energy phase function} as
\begin{align}
\label{psiE-tree}
\Psi_E(\TT)
:=
-\Re\bigg(- \frac{\xi^3}{|\xi|^2} + \sum_{\bul\in \TT^\infty} \frac{\xi_\bul^3}{|\xi_\bul|^2}\bigg).
\end{align}
The partial energy phase function is given as
\begin{align}
\label{psi-circ-E}
\Psi_E(\star) = \Psi_E(\cherry{\star}{}{}) = 
-\Re \big(- \tfrac{\xi_{\star}^3}{|\xi_{\star}|^2} + \tfrac{\xi_{\star1}^3}{|\xi_{\star1}|^2}  + \tfrac{\xi_{\star2}^3}{|\xi_{\star2}|^2} \big) .
\end{align}
\item[(iv)] (Phase function) For a fixed binary tree $\TT$, the (total) phase function is then defined as
\begin{equation}\label{PsiT}
\Psi(\TT)=\Psi_0(\TT) + E \Psi_E(\TT).
\end{equation}
Given $\star\in \TT^0$, we write
\begin{align}\label{Psi-star}
\Psi(\star)=\Psi_0(\star) + E \Psi_E(\star). 
\end{align}
\item[(v)] 
Lastly, for $\star\in \TT^0$, we define the multiplier $K_\star$ given by
\begin{equation}
\label{Kmul}
K_{\star} = K(\xi_{\star1},\xi_{\star2}) := \frac{\xi_{\star1}\cdot\xi_{\star2}}{|\xi_{\star1}|^2|\xi_{\star2}|^2},
\end{equation}
where $\xi_{\star1}\cdot\xi_{\star2}$ is the $\R^2$-inner product.

\end{enumerate}

\end{definition}

Under the above definitions, we introduce a correspondence between binary trees and oscillatory integrals. 
Given $j\in\N$, a binary tree $\TT \in \Tr_j$, a multiplier $m:\CC^{\TT}\to \CC$ on $\TT$, and a
function $v: \R\times \CC \to\CC$, we define the associated integral term
\begin{equation}\label{eq:defiI}
\Ical(\TT, m; v)(t) :=   i \int_{\G(\TT)}
e^{it\Psi(\TT)}m(\vec{\xi})
\prod_{\bul\in \TT^\infty}  \ft v_\bul (t) d\vec{\xi},
\end{equation}
where 
we use the short-hand notation  $\ft v_\bul(t) := \ft v(t,\xi_\bul)$ for a label $\bul$. 

\subsection{Deriving the normal form equation}\label{sec:nfe}

Using the tree formalism in Subsection~\ref{SUB:tree}, we now derive the normal form equation. 
We apply a normal form reduction in the form of ``differentiation-by-parts" as in \cite{bit}, only in a problematic region where we cannot show the nonlinear estimate at low regularity.
In particular, recalling from \eqref{Gamma} that $\G(\<1>)=\big\{(\xi,\xi_1,\xi_2)\in\CC^3:\ \xi=\xi_1+\xi_2\big\}$ is the convolution set of the elementary tree $\<1>$, we split it into the following two subsets $A$ and $A^c=\G(\<1>)\setminus A$:
\begin{align}
\label{A}\begin{split}
A & := \big\{ (\xi,\xi_{1}, \xi_{2}) \in \G(\<1>): \  |\xi_{1}| \land  |\xi_{2}| \le 1/\delta  \big\}, \\
A^c & := \big\{ (\xi,\xi_{1}, \xi_2) \in \G(\<1>): \  |\xi_{1}| \land |\xi_{2}| > 1/\delta  \big\},
\end{split}
\end{align}
where $A^c$ denotes the bad region where we perform a reduction.

Using the tree-indexing in Subsection~\ref{SUB:tree} and the notations in \eqref{psi0-tree}, \eqref{psi-circ}, and \eqref{Kmul}, the equation satisfied by $\ft v= \F_{x,y}(U(-t)u)$, in the interaction representation, can be written as follows:
\begin{align}
\label{duhamel2}
\dt\ft v(t) =  i \int_{\xi = \xi_1 + \xi_2} e^{it \Psi(\cherryS)}  K \Psi_0(\<1>)   (\widehat{v}_{1} \widehat{v}_{2}) (t)
d\xi_{1}.
\end{align}
Absolute multiplicative constants have been normalized to $1$; see the footnote on p.~4.
Then, splitting the region of integration into $A\cup A^c$ as in \eqref{A}, we may write \eqref{duhamel2} as 
\begin{align}
\dt\ft v
&
= 
\Ical\big( \<1>, \ind_A K \Psi_0(\<1>); v 
\big)
+ \Ical\big( \<1>, \ind_{A^c} K \Psi_0(\<1>) ;v 
\big)
\label{NFE000}
\end{align}
where $\Ical$ is as in \eqref{eq:defiI}. In order to perform a differentiation-by-parts in time on the term which is restricted to the bad region $A^c$, we need to make the total phase $\Psi(\<1>)$ appear. To that end, we simply recall that 
\[
\Psi(\<1>) = \Psi_0(\<1>) + E \Psi_E (\<1>).
\]
This allows us to write 
\begin{align}
\dt\ft v
&
= 
\Ical\big( \<1>, \ind_A K \Psi_0(\<1>); v 
\big)
+ \Ical\big( \<1>, \ind_{A^c} K \Psi (\<1>) ;v 
\big)
- \Ical \big(\<1>, \ind_{A^c} K \cdot E\Psi_E(\<1>);v\big) 
\\
& 
= \Ical \big(\<1>, \ind_A K \Psi_0(\<1>) - \ind_{A^c} K \cdot E  \Psi_E(\<1>);v\big) + \Ical(\<1>, \ind_{A^c} K \Psi(\<1>);v) 
\\
& 
=: \F_{x,y}[ \NN ( \<1>; v )] +  \F_{x,y}[ \MM(\<1>; v) ].
\label{NFE001}
\end{align}
In order to move on with our development, we use the relation
$$
\dt\frac{e^{it\Psi(\cherryS)}}{i\Psi(\<1>)}=e^{it\Psi(\cherryS)},
$$ 
which allows us to perform the first differentiation-by-parts step in our chain. We obtain
\begin{align}
\label{NFE01}
\begin{split}
\dt\ft v
&=
\F_{x,y}[ \NN ( \<1>; v )] 
+ \dt\bigg[ \int_{\xi = \xi_1 + \xi_2} e^{it \Psi(\cherryS)} 
[\ind_{A^c}
 K ] 
(\ft v_{1} \ft v_{2} )(t)
d\xi_{1} 
\bigg]
\\
&
\quad 
-
   \int_{\xi = \xi_1 + \xi_2} e^{it \Psi(\<1>)} 
[ \ind_{A^c} K ]
\partial_{t} (\ft v_{1} \ft v_{2} )(t)
d\xi_{1} 
\\
& =
\F_{x,y}[ \NN ( \<1>; v )] 
+  \dt\bigg[ \int_{\xi = \xi_1 + \xi_2} e^{it \Psi(\cherryS)} 
[\ind_{A^c}
 K ] 
(\ft v_{1} \ft v_{2} )(t)
d\xi_{1} 
\bigg]
\\
&
\quad 
- i
 \int_{\substack{\xi = \xi_1 + \xi_2 \\ \xi_1 = \xi_{11}+ \xi_{12}  }} \exp\Big(it \Psi(\<21>) \Big)
[ \ind_{A^c} K K_1 \Psi_0(1) ]
 (\ft v_{11} \ft v_{12} \ft v_2 )(t)
d\xi_{1} d\xi_{11} 
\\
&
\quad 
-i  \int_{\substack{\xi = \xi_1 + \xi_2 \\ \xi_2 = \xi_{21}+ \xi_{22}  }} \exp\Big(it \Psi(\<22>) \Big)
[ \ind_{A^c} K K_2 \Psi_0(2) ]
 (\ft v_{1} \ft v_{21} \ft v_{22} )(t)
d\xi_{1} d\xi_{21}
\end{split}
\end{align}
where we replaced $\dt \ft v_j$ by \eqref{duhamel2} to obtain the last equality. This can be rewritten as
\begin{align}
\label{NFE00}
\begin{split}
\dt\ft v 
&
=
 \F_{x,y}[ \NN ( \<1>; v )] 
+ 
\dt\Ical \big( \<1>,  - i \ind_{A^c} K; v \big) 
\\
&
\quad 
+
\Ical \Big( \<21>, - \ind_{A^c} KK_1\Psi_0(1) ; v  \Big) 
+
\Ical \Big( \<22>, - \ind_{A^c} KK_2\Psi_0(2) ; v  \Big) 
. 
\end{split}
\end{align}

In order to iterate the procedure taking \eqref{duhamel2} into \eqref{NFE00}, we use the following result. It provides a systematic way to perform a normal form reduction on integral terms of the form $\Ical$ in \eqref{eq:defiI}, which translates to an operation on trees (based on tree concatenation as defined in Definition~\ref{DEF:tree3}(ii)) and suitable modifications of the multipliers. We refer the reader to the companion KdV paper~\cite[Lemma~2.6]{chcora} for the proof, which works in the exact same fashion in our setting.

\begin{lemma}[Differentiation-by-parts on trees]
\label{LEM:IBP}
Let $j\in\N$, $\TT\in\Tr_j$, $m: \CC^{\TT}\to \mathbb{C}$ a multiplier on $\TT$,
and $v$ be a smooth function such that $\ft v$ satisfies \eqref{duhamel2}. Then, the following identity holds
\begin{align}\label{IBP}
\begin{split}
\Ical (\TT,m; v)
= \dt \Ical \Big(\TT,\frac{m}{ i \Psi(\TT)} ; v \Big)
+ \sum_{\star \in \TT^\infty} \Ical 
\Big(\TT \overset{\star}{\oplus} \<1>
, - K_\star \Psi_0(\star) \frac{m}{\Psi(\TT)} ; v \Big)
,
\end{split}
\end{align}
where $\Ical$ is as in \eqref{eq:defiI}, $\TT \overset{\star}{\oplus} \<1>$ denotes tree concatenation at the final node $\star\in \TT^\infty$ as in Definition~\ref{DEF:tree3}(ii), and $\Psi(\TT)$, $\Psi_0(\star)$, $K_\star$ are as in \eqref{PsiT},  \eqref{psi-circ}, and~\eqref{Kmul}.

\end{lemma}

As in \eqref{NFE01}, we split each region of frequency integration between a good and a bad set, only performing a differentiation-by-parts when restricted to the latter. Generalizing the definition of $A, A^c$ in \eqref{A}, given a tree $\TT$ and a parent node $\bul \in \TT^0$, we define the analogous sets $A_\bul$ and $A_\bul^c$ as
\begin{align}
A_\bul& 
:=
 \big\{ (\xi_\bul,\xi_{\bul1}, \xi_{\bul2}) \in \CC^3: \ \xi_\bul =
\xi_{\bul1}
+ \xi_{\bul2}, \ |\xi_{\bul1}| \land |\xi_{\bul2}|  \le 1/\delta 
 \} , \\
A_\bul^c & := \big\{ (\xi_\bul,\xi_{\bul1}, \xi_{\bul2}) \in \CC^3: \ \xi_\bul =
\xi_{\bul1} +
\xi_{\bul2},\
|\xi_{\bul1}| \land |\xi_{\bul2}| \ge 1/\delta \big\} 
.
\label{defA}
\end{align}
Note that, when $\bul = \emptyset$ is the root node, $A_\bul =A$ as in \eqref{A}. 
We will often omit the dependence on the frequencies $\xi_\bul, \xi_{\bul 1}, \xi_{\bul2}$, and use the short-hand notation $\ind_{A_\bul}=\ind_{A_\bul} (\xi_\bul,\xi_{\bul1},\xi_{\bul2})$.

Returning to \eqref{NFE00}, we continue to derive the second iteration of the normal form equation.
Unlike the first normal form reduction from equation \eqref{duhamel2} to \eqref{NFE00}, where we split the integral term into $A$ and $A^c$ (see \eqref{A}) and differentiate-by-parts in the latter region, we now have an additional step, which we describe in the following. Considering the two last terms in \eqref{NFE00}, we first split their regions of integration $(\xi_\bul,\xi_{\bul1}, \xi_{\bul2}) \in A_\bul \cup A_\bul^c$, where $\bul = 1,2$, respectively:
\begin{align}
\label{NFE001b}
\dt\ft v
&
= 
 \F_{x,y}[ \NN ( \<1>; v )]  
+ 
\dt \Ical \big( \<1>, -i  \ind_{A^c} K; v \big)
\\
&
\quad 
+
\Ical \Big( \<21>, - \ind_{A^c \cap A_1} KK_1\Psi_0(1) ; v  \Big)  
+
\Ical \Big( \<21>, - \ind_{A^c \cap A_1^c} KK_1\Psi_0(1) ; v  \Big) 
\\
&
\quad 
+
\Ical \Big( \<22>, - \ind_{A^c \cap A_2} KK_2\Psi_0(2) ; v  \Big) 
+
\Ical \Big( \<22>, - \ind_{A^c \cap A_2^c} KK_2\Psi_0(2) ; v  \Big) .
\end{align}
As a first preparatory step, we add and subtract the `energy' terms, in order to complete the phases in the 
terms restricted to $A^c \cap A^c_j$, $j=1,2$. 
Indeed, from \eqref{Psi-star}, we first write 
\[
\Psi_0(\star) = \Psi(\star) - E \Psi_E(\star)
\]
which allows us to write 
\begin{align}
\label{NFE0011}
\dt\ft v
&
= 
 \F_{x,y}[ \NN ( \<1>; v )]  
+ 
\dt \Ical \big( \<1>, -i  \ind_{A^c} K; v \big)
\\
&
\quad 
+
\Ical \Big( \<21>, K K_1 \ind_{A^c} \cdot \left( E \ind_{A_1^c}  \Psi_E(1) - \ind_{A_1} \Psi_0(1) \right) ; v  \Big)  
+
\Ical \Big( \<21>, - \ind_{A^c \cap A_1^c} KK_1\Psi(1) ; v  \Big) 
\\
&
\quad 
+
\Ical \Big( \<22>, KK_2 \ind_{A^c} \left( E \ind_{A_2^c} \Psi_E(2) - \ind_{A_2} \Psi_0(2)\right) ; v  \Big) 
+
\Ical \Big( \<22>, - \ind_{A^c \cap A_2^c} KK_2\Psi(2) ; v  \Big). 
\end{align}

As before, we keep the integral contributions which either come from a ``lower-order'' or energy term, or which are localized to a ``good set" $A_\bul$. Regarding the contributions restricted to $A_\bul^c$ for all $\bul\in \TT^0$, despite the restriction to the ``bad regions'', we will keep a suitable part of these terms before performing another differentiation-by-parts. In particular, by using the additivity of the phase $\Psi$ in \eqref{eq:additivity}, we write
\begin{align}
\Psi(\star) 
= \big[ \Psi(\star) + \Psi(\<1>) \big]- \Psi(\<1>) 
=
\Psi(\TT) - \Psi(\<1>),
\quad (\star, \TT)\in \left\{ (1, \<21>) , (2, \<22>) \right\},
\end{align}
and keep the contributions with partial phase $-\Psi(\<1>)$, only differentiating-by-parts the terms with total phase $\Psi(\TT)$. Then, from \eqref{NFE001}, we obtain the second iteration of the normal form equation
\begin{align}
\dt\ft v
&
= 
 {\F_{x,y}}[ \NN ( \<1>; v )] 
+ 
\dt \Ical \big( \<1>, - i  \ind_{A^c} K; v \big)
\\
&
+
\Ical \Big( \<21>,  K K_1 \ind_{A^c} \cdot \left( \ind_{A_1^c} \Psi(\<1>) + E \ind_{A_1^c}  \Psi_E(1) - \ind_{A_1} \Psi_0(1) \right); v  \Big) 
\\
&
+
\Ical \Big( \<21>, - \ind_{A^c \cap A_1^c} KK_1\Psi(\<21>) ; v  \Big) 
\\
&
+
\Ical \Big( \<22>, KK_2 \ind_{A^c} \left( \ind_{A_2^c} \Psi(\<1>) + E \ind_{A_2^c} \Psi_E(2) - \ind_{A_2} \Psi_0(2)\right) ; v  \Big) 
\\
&
+
\Ical \Big( \<22>, - \ind_{A^c \cap A_2^c} KK_2\Psi(\<22>) ; v  \Big) 
\\
&
= : 
 {\F_{x,y}}[ \NN ( \<1>; v )]  
+ 
 {\F_{x,y}}[ 
\BB(\<1>;v)
] 
\\
&
\quad 
+ 
 {\F_{x,y}}[ 
\NN(\<21>; v )
] 
+
 {\F_{x,y}}[ 
\MM(\<21>; v )
] 
+
 {\F_{x,y}}[ 
\NN(\<22>; v)
]  
+
 {\F_{x,y}}[ 
\MM(\<22>;v)
]  
.
\label{NFE002}
\end{align}
In the Duhamel expansion of $v$, based on \eqref{NFE002}, $\NN$ corresponds to time-integral terms, $\BB$ to boundary terms, and $\MM$ are remainders, to which we apply \eqref{IBP} at the next step. Thus, in the following, we refer to them as integral terms, boundary terms, and remainders, respectively.

In the subsequent iterations of the normal form equation, as in \eqref{NFE00} and \eqref{NFE002}, for each tree $\TT$, we want to have a unique integral term $\NN(\TT)$ and boundary term $\BB(\TT)$, 
where this must encode all the possible ways of arriving at such a tree structure. 
From Lemma~\ref{LEM:IBP}, we note that given a tree $\TT \in \Tr_j$ for $j\ge2$, after differentiation-by-parts, the resulting contributions associated with the tree $\TT$ can come from any term of the form $\Ical(\TT \overset{\star}{\ominus} \<1>)$ where $\star \in \TT^{0, \ff}$, i.e., that the children of $\star$ are leaves in $\TT$. Consequently,  we need to perform differentiation-by-parts on the remainder terms $\MM(\TT \overset{\star}{\ominus} \<1>)$ simultaneously and combine any contributions associated with the same tree before preparing the integral terms to be kept and those which will be differentiated in the following step. 

The following result provides a systematic way of repeating the process described above, with an explicit description of the contributions and multipliers involved at each finite normal form equation. 
We note that we see exactly one boundary term and one integral term for each tree of a given generation $\Tr_j$, with the complexity of the contributions being encoded in the corresponding multipliers. The proof is a direct adaptation of the corresponding iteration in \cite[Proposition 2.7]{chcora}; the only structural difference is the energy contribution, which through the splitting $K\Psi_0=K\Psi-E\,K\Psi_E$ enters the remainder terms and is propagated perturbatively, leaving the tree combinatorics unchanged.

\begin{proposition}
\label{PRO:NFE}
Let $J\in\N$, $u$ be a smooth solution to \eqref{eq:NV}, and $\ft v={\F_{x,y}}(U(-t)u)$ be the corresponding solution to \eqref{duhamel2}.
Then, 
the $J$-th normal form equation for $v$ is 
\begin{align}
\label{NFE1}
\dt v=
\sum_{j=1}^{J-1} \sum_{\TT\in \Tr_j} \BB(\TT; v) 
+ \sum_{j=1}^{J} \sum_{\TT\in \Tr_j} \NN(\TT; v) 
+ \sum_{\TT\in \Tr_{J}}
\MM(\TT; v),
\end{align}
where the operators are defined via their spatial Fourier transforms
\begin{align}
\label{NFE1a}
\begin{split}
{\F_{x,y}} [\BB(\TT;v)]  & = \dt\Ical(\TT, - i \mulB(\TT); v),\\
{\F_{x,y}} [\NN(\TT;v)]  & = \Ical(\TT, \mul(\TT) ; v  ), \\
{\F_{x,y}} [\MM(\TT;v)] & = \Ical(\TT, \Psi(\TT)\mulB(\TT); v),
\end{split}
\end{align}

\noi
with $\Ical$ as in \eqref{eq:defiI}.
Moreover,
the multipliers $\mul(\TT),\mulB(\TT)$ for $\TT\in \Tr_j$ are given~by
\begin{align}
\mulB(\TT)
& =  (-1)^{j+1}  \prod_{\bul \in \TT^0}  \ind_{A^c_{\bul}} K_\bul
,
\\
\mul(\TT) & =  (-1)^{j+1} \bigg(\prod_{\bul \in \TT^0} K_\bul \bigg) \bigg\{  \sum_{\star \in \TT^{0,f}}
\left(\ind_{A_{\star}} \Psi_0(\star) - \ind_{A_{\star}^c} E\Psi_E(\star)\right)\bigg(\prod_{\substack{\bul \in \TT^0 \setminus\{\star\}
}} \ind_{A^c_{\bul}} \bigg)\\
&
\hspace{4cm}
-
\bigg( \prod_{\bul \in \TT^0} \ind_{A_\bul^c} \bigg)
\sum_{\star\in \TT^0 \setminus \TT^{0,\ff}} \Psi(\star)
 \bigg\},
 \label{NFE0b}
\end{align}
where $K_\bul$, $\Psi_0(\star)$, $\Psi_E(\star)$, and $\Psi(\star)$ are as in \eqref{Kmul}, \eqref{psi-circ}, \eqref{psi-circ-E}, and \eqref{Psi-star}. 

\end{proposition}

Taking $J\to\infty$ in Proposition~\ref{PRO:NFE} formally yields the \emph{infinite normal form equation}
\begin{align}
\label{NFE}
\dt v= \sum_{j=1}^{\infty} \sum_{\TT\in \Tr_j} \left(\NN(\TT;v) + \BB(\TT;v) \right),
\end{align}
with $\NN(\TT;v)$ and $\BB(\TT;v)$ as in \eqref{NFE1a}.

Our route from~\eqref{NFE} to the gauged equation has three stages. First,  we define a nonlinear change of variables $v\mapsto w$, based on the expression for the first boundary term in \eqref{NFE}, which eliminates  \emph{every} boundary term via an algebraic identity (Lemma~\ref{LEM:wbd}). Second, an analogous cancellation on the integral terms greatly reduces their algebraic complexity. In the end, only five different types of nonlinear terms will remain - see the precise description in Lemma~\ref{LEM:wint}.
Third, the interaction representation $z(t):=U(t)w(t)$ recasts the system as the gauged Novikov--Veselov equation~\eqref{z1}. The algebra runs in close parallel with~\cite{chcora}, with one structural novelty: alongside the dispersive phase $\Psi_0$, the energy phase $E\Psi_E$ contributes its own family of multipliers (Type~$\III$), so the inclusion-exclusion bookkeeping is carried out on the pair $(\Psi_0,E\Psi_E)$ rather than on a single phase.

\subsection{The gauged NV equation}

It will be convenient to introduce a multilinear version of $\Ical$ from~\eqref{eq:defiI}. Fix $j\in\N$ and $\TT\in\Tr_j$, and order the $j+1$ leaves of $\TT$ lexicographically by their words, left-to-right in any planar drawing:
\begin{align}
\label{lexic}
\bul_1\le \bul_2\le \dots \le \bul_{j+1}.
\end{align}
For smooth functions $v_1,\dots,v_{j+1}:\R\times \CC\to\R$, set
\begin{equation}
\label{eq:defiImulti}
\Ical(\TT, m; v_1,\dots,v_{j+1}) (t)
= i \int_{\Gamma(\TT)} e^{it\Psi(\TT)} m(\vec{\xi}) \prod_{\ell=1}^{j+1} \ft v_\l (t,\xi_{\bul_\ell})\, d\vec{\xi}
,
\end{equation}
with $\Psi(\TT)$ as in \eqref{PsiT}.
Specialising $v_1=\cdots=v_{j+1}=v$ recovers $\Ical(\TT, m; v)$ of~\eqref{eq:defiI}. The multilinear counterparts of $\BB(\TT)$ and $\NN(\TT)$ from~\eqref{NFE1a} are then
\begin{align*}
{\F_{x,y}} [\BB(\TT; v_1, \ldots, v_{j+1}) ]  &= \dt\Ical (\TT,- i \mulB(\TT) ; v_1,\dots,v_{j+1}),\\
{\F_{x,y}} [\NN(\TT; v_1, \ldots, v_{j+1} ) ]  &= \Ical(\TT,\mul(\TT) ; v_1, \ldots, v_{j+1} ).
\end{align*}

Define $v\mapsto w$ implicitly by
\begin{align}
\label{w}
v = w + \Ical(\<1>, - i \mulB(\<1>); w) =: w + \D[w, w],
\end{align}
with $\Ical$ as in~\eqref{eq:defiI}. By design, $\partial_t\D[w,w]$ reproduces precisely the quadratic boundary contribution $\BB(\<1>;w)$, with $\BB$ as in~\eqref{NFE1a} (cf.~\eqref{NFE002}), and thus
\begin{align}
\label{wa}
\dt v = \dt w + \BB(\<1>; w).
\end{align}
Plugging~\eqref{wa} and~\eqref{w} into~\eqref{NFE} gives the closed~equation
\begin{align}
\dt w
& = -\BB(\<1>;w) + \BB(\<1>;v) + \NN(\<1>;v) + \sum_{j=2}^\infty \sum_{\TT\in\Tr_j} \Big( \NN(\TT; v)  + \BB(\TT; v) \Big) \nonumber \\
& =  \BB\big(\<1>; w, \D[w, w] \big) + \BB \big(\<1>; \D[w, w],w\big) + \BB\big(\<1>; \D[w,w], \D[w, w]\big) \nonumber\\
& \qquad + \NN(\<1>; w ,w )  + \NN\big(\<1>; w, \D[w, w]\big) + \NN\big(\<1>; \D[w, w], w\big) + \NN\big(\<1>; \D[w, w],\D[w, w]\big) \nonumber\\
& \qquad + \sum_{j=2}^\infty \sum_{\TT\in \Tr_j} \Big( \NN\big(\TT; w+ \D[w, w] \big) +  \BB\big(\TT; w + \D[w, w]\big) \Big).
\label{w1}
\end{align}

Because $\D[w,w]$ is itself quadratic in $w$, expanding~\eqref{w1} produces, for each underlying tree, several terms of the same homogeneity but supported on different  \emph{colourings} (by a colouring of a tree we mean an assignment to each of its leaves of either $w$ or $\D[w,w]$).  After a careful grouping of the resulting boundary terms, depending on the underlying tree structure, a striking algebraic phenomenon is revealed: \emph{every} boundary term is annihilated.

\begin{lemma}[Boundary cancellation]
\label{LEM:wbd}
Let $j\in\N$, $\TT\in \Tr_j$, $\star \in \TT^{0,\ff}$, and let $w_1, \ldots, w_{j+1}: \R\times \CC \to \R$ be smooth. We order the leaves of $\TT$ as in~\eqref{lexic}, and let $1 \le k \le j$ be such that $(\bul_k, \bul_{k+1})$ are the children of $\star$ in $\TT$. Then
\begin{equation}
\label{wbd0}
\BB\big(\TT ; w_1, \ldots, w_{j+1} \big)
+ \BB\big(\TT\overset{\star}{\ominus} \<1>; w_1,  \ldots, w_{k-1} , \D[w_k,w_{k+1}] , w_{k+2}, \ldots,  w_{j+1} \big) \equiv 0.
\end{equation}
Consequently, all the boundary terms in~\eqref{w1} vanish:
\begin{multline}
\label{wbd1}
  \BB \big(\<1>; w, \D[w, w] \big)
+ \BB\big(\<1>; \D[w, w],w\big)
+ \BB\big(\<1>;\D[w, w], \D[w, w]\big)
\\
+ \sum_{j=2}^\infty \sum_{\TT\in \Tr_j}  \BB\big(\TT; w + \D[w, w]\big)  \equiv 0  .
\end{multline}
\end{lemma}

\begin{proof}
Write $\TT' := \TT\overset{\star}{\ominus}\<1>$, so that $(\TT')^0 = \TT^0\setminus\{\star\}$ has cardinality $j-1$. Since $\D = \Ical(\<1>, -i\mulB(\<1>); \cdot)$, the multiplier of the second summand in~\eqref{wbd0} splits as $-i\mulB(\TT')\cdot \mulB(\cherry{\star}{}{})$, with the auxiliary factor at $\star$ given by $\mulB(\cherry{\star}{}{}) = \ind_{A^c_\star}K_\star$. Phase additivity~\eqref{eq:additivity} forces $\Psi(\TT')+\Psi(\cherry{\star}{}{}) = \Psi(\TT)$, so both summands integrate the same exponential over the same domain (after re-indexing the leaves by~\eqref{lexic}). Adding the multipliers gives the one-line identity
\begin{align}
\label{wbd-cancel}
\mulB(\TT) + \mulB(\TT')\,\mulB(\cherry{\star}{}{})
= 
\big[(-1)^{j+1} + (-1)^{j}\big] \prod_{\bul \in \TT^0} \ind_{A^c_\bul} K_\bul = 0,
\end{align}
which proves~\eqref{wbd0}.

For~\eqref{wbd1}, we group every boundary contribution in~\eqref{w1} whose uncoloured underlying tree is $\TT$ (compare with~\eqref{wbd-cancel}); the cumulative multiplier reads
\begin{align*}
\wt{\mul}_\BB(\TT) 
&= 
\sum_{P\subseteq \TT^{0,\ff}} \mulB\big(\TT\overset{P}{\ominus}\<1>\big) \prod_{\star\in P}\mulB(\cherry{\star}{}{})
\\
&= 
(-1)^{j+1}\bigg(\prod_{\bul\in \TT^0}\ind_{A_\bul^c}K_\bul\bigg)\sum_{P\subseteq \TT^{0,\ff}}(-1)^{|P|},
\end{align*}
and vanishes by $\sum_{k=0}^n\binom{n}{k}(-1)^k=0$ with $n=|\TT^{0,\ff}|\ge 1$. Identical bookkeeping is carried out term-by-term in~\cite[Lemma~2.8]{chcora}.
\end{proof}

\begin{remark}
\label{remark:factorial}
The cancellation in Lemma~\ref{LEM:wbd} relies in a crucial way on the fact that the constraint $A^c_\bul$ on $(\xi_\bul,\xi_{\bul1},\xi_{\bul2})$ depends only on the local node $\bul$. By contrast, in the standard normal form approach (e.g.~\cite{GKO13, kishimoto}) the resonant set involves all ancestors of $\bul$, which entangles the trees with their growth history, preventing a cancellation of the form in Lemma~\ref{LEM:wbd}.
Moreover, it obliges one to introduce a $j$-dependent sequence of small parameters to tame the factorial growth of the number of ordered trees. We refer the reader to~\cite{chcora} for a more complete discussion.
\end{remark}

The same mechanics leading to the cancellation for the boundary terms $\BB$ in Lemma~\ref{LEM:wbd} also propagates to the integral terms $\NN$ in \eqref{w1}, with only five families of contributions remaining, namely $\NN_\l$ for $\l=1, \ldots, 5$ which we refer to as Type $\I$-$\5$ terms, respectively.

\begin{lemma}
\label{LEM:wint}
Let $w:\R \times \R^2 \to \R$ be smooth and $\NN$ as in~\eqref{NFE1a}. Then
\begin{align}
\label{wint0a}
\sum_{j=1}^\infty \sum_{\TT\in \Tr_j} \NN\big(\TT; w+ \D[w, w] \big)
=
\sum_{j=1}^\infty \sum_{\TT\in\Tr_j} \sum_{\ell=1}^5  \NN_\ell(\TT; w),
\end{align}
where, for $\l = 1,\dots,5$,
\begin{align}
\label{wint0b}
{\F_{x,y}} \big[ \NN_\ell (\TT; w) \big] = \Ical \big(\TT, {\mul}_\l(\TT); w \big),
\end{align}
with $\Ical(\TT, m; w)$ as in~\eqref{eq:defiI}, and the multipliers given by
\begin{align}
\label{eq:multiplicadores}
\mul_1(\<1>) &= K\big(\ind_A  \Psi_0(\<1>) - \ind_{A^c} E\Psi_E(\<1>)\big), \\
\mul_2(\TT) &=(-1)^{j+1} \ind_{A_\star} K_\star\Psi_0(\star) \bigg( \prod_{\bul \in \TT^0 \setminus\{\star\}}  K_\bul  \ind_{A_\bul^c} \bigg),
   && \hspace{-1.0cm}\text{if} \  \TT^{0,\ff} = \{\star\}, \ \star\neq \emptyset, \nonumber\\
\mul_3(\TT) &=(-1)^{j} E  \ind_{A_\star^c} K_\star\Psi_E(\star) \bigg( \prod_{\bul \in \TT^0 \setminus\{\star\}}  K_\bul  \ind_{A_\bul^c} \bigg),
   && \hspace{-1.0cm}\text{if} \  \TT^{0,\ff} = \{\star\}, \ \star\neq \emptyset, 
\nonumber\\
\mul_4(\TT) &=(-1)^{j}  K_{\pb(\star)} 
\Psi_0(\pb(\star))
\bigg( \prod_{\bul \in \TT^0 \setminus\{\pb(\star)\}}  K_\bul  \ind_{A_\bul^c} \bigg),
   && \hspace{-1.0cm}\text{if} \  \TT^{0,\ff} = \{\star\},  \ \star\neq \emptyset, 
\label{wint0c}\\
\mul_5(\TT) &=(-1)^{j+1}  K_\star \Psi_0(\star) \bigg( \prod_{\bul \in \TT^0 \setminus\{\star\}}  K_\bul  \ind_{A_\bul^c} \bigg),
   && \hspace{-1.0cm}\text{if} \  \TT^{0,\ff} = \{\star1, \star2\}, 
\nonumber\\
\mul_\l(\TT) &= 0, \quad \l=1,2,3,4,5, && \hspace{-1.0cm}\text{otherwise}, \nonumber
\end{align}
with $\TT\in \Tr_j$, $j \ge 1$, and $K$, $\Psi_0$, and $\Psi_E$ and as in~\eqref{Kmul}, \eqref{psi-circ}, and \eqref{psi-circ-E}, respectively. 
We refer to the contributions $\NN_\l(\TT)$ as terms of Type $\I, \II, \III, \IV, \5$, respectively.
\end{lemma}

\begin{proof}
The argument follows the inclusion-exclusion strategy of~\cite[Lemma~2.10]{chcora} 
verbatim, except for one new ingredient: the bookkeeping induced by the energy phase $E\Psi_E$, which we record below.

\medskip

For $j=1$,
on the l.h.s.\ of~\eqref{wint0a}, only term is $\NN(\<1>;w,w)$. Hence $\NN_1(\<1>;w)=\NN(\<1>;w,w)$ and $\mul_1(\<1>) = \mul(\<1>) = K(\ind_A\Psi_0(\<1>) - \ind_{A^c} E\Psi_E(\<1>))$, matching~\eqref{eq:multiplicadores}.

Next, fix $j\ge 2$ and $\TT\in\Tr_j$. Following~\cite[Lemma~2.10]{chcora}, collect all l.h.s.\ contributions whose underlying uncoloured tree (after expanding each $\D[w,w]$) is $\TT$. Indexing such expansions by $C\subseteq \TT^{0,\ff}$, the total multiplier is
\[
\wt{\mul}(\TT)
= \sum_{C\subseteq \TT^{0,\ff}} \mul\big(\TT\overset{C}{\ominus}\<1>\big) \prod_{\star\in C}\mulB(\cherry{\star}{}{}).
\]
The splitting $\Psi(\star) = \Psi_0(\star) + E\Psi_E(\star)$ in the definition of $\mul$ (see~\eqref{NFE0b}) together with $\mulB(\cherry{\star}{}{}) = \ind_{A_\star^c}K_\star$ and $(\TT\overset{C}{\ominus}\<1>)^0 = \TT^0\setminus C$ decomposes
\begin{align*}
\wt{\mul}(\TT) &= \wt{\mul}_{\gfrak,0}(\TT) + \wt{\mul}_{\gfrak,E}(\TT) + \wt{\mul}_{\bfrak}(\TT), \quad \text{with}
\\
\wt{\mul}_{\gfrak,0}(\TT)
&
= (-1)^{j+1} \bigg(\prod_{\bul\in \TT^0}K_\bul\bigg) \sum_{C\subseteq \TT^{0,\ff}}
(-1)^{|C|} \sum_{\star\in (\TT\overset{C}{\ominus}\cherryS)^{0,\ff}}
\ind_{A_\star}\Psi_0(\star) \prod_{\bul\in \TT^0\setminus\{\star\}}\ind_{A_\bul^c}
, 
\\
\wt{\mul}_{\gfrak,E}(\TT)
&
=
(-1)^j E \bigg(\prod_{\bul\in \TT^0}K_\bul\bigg)
\sum_{C\subseteq \TT^{0,\ff}}(-1)^{|C|}
\sum_{\star\in (\TT\overset{C}{\ominus}\cherryS)^{0,\ff}}
\ind_{A_\star^c}\Psi_E(\star)
\prod_{\bul\in \TT^0\setminus\{\star\}}\ind_{A_\bul^c}
, 
\\
\wt{\mul}_{\bfrak}(\TT)
&
=(-1)^j 
\bigg(\prod_{\bul\in \TT^0}\ind_{A_\bul^c}  K_\bul\bigg)
\sum_{C\subseteq \TT^{0,\ff}}(-1)^{|C|} 
\sum_{\star \in \TT^0} \Psi(\star)
\ind_{\star \in (\TT\overset{C}{\ominus}\cherryS)^{0} \setminus (\TT\overset{C}{\ominus}\cherryS)^{0,\ff} }
.
\end{align*}

The simplifications of $\wt{\mul}_\bfrak(\TT)$ and $\wt{\mul}_{\gfrak,0}(\TT)$ are imported, line for line, from~\cite[Lemma~2.10]{chcora}, yielding
\begin{align}
\label{eq:mulb-nv}
\wt{\mul}_{\bfrak}(\TT) &=(-1)^j \bigg(\prod_{\bul\in \TT^0}K_\bul\ind_{A_\bul^c}\bigg) \sum_{\star\in \TT^0\setminus \TT^{0,\ff}} \Psi(\star)
\Big( \ind_{\TT^{0,\ff}=\{\star1\}} + \ind_{\TT^{0,\ff}=\{\star2\}} - \ind_{\TT^{0,\ff}=\{\star1,\star2\}}\Big),\\
\label{eq:mulg0-nv}
\wt{\mul}_{\gfrak,0}(\TT)
&= (-1)^{j+1} \sum_{\star\in \TT^0} \ind_{A_\star}K_\star\Psi_0(\star) \bigg(\prod_{\bul\in \TT^0\setminus\{\star\}}K_\bul\ind_{A_\bul^c}\bigg) \nonumber\\
&\quad \times\Big( \ind_{\TT^{0,\ff}=\{\star\}} - \ind_{\TT^{0,\ff}=\{\star1\}} - \ind_{\TT^{0,\ff}=\{\star2\}} + \ind_{\TT^{0,\ff}=\{\star1,\star2\}}\Big).
\end{align}

\medskip

Moving now to $\wt{\mul}_{\gfrak,E}$, 
fix $\star\in \TT^0$. 
Compared to $\wt{\mul}_{\gfrak,0}$, the indicator $\ind_{A_\star}$ becomes $\ind_{A_\star^c}$ and the factor $\Psi_0$ becomes $E\Psi_E$, while the index set of the inner alternating sum stays put, namely

\[
\big\{ C \subseteq \TT^{0,\ff} \, : \, \star \in (\TT\overset{C}{\ominus}\<1>)^{0,\ff} \big\}.
\]
Three cases arise.

\noindent\underline{Case 1: $\star\in\TT^{0,\ff}$.} As $\star\in (\TT\overset{C}{\ominus}\<1>)^{0,\ff}$ holds iff $\star\notin C$, we have
\[
\sum_{\substack{C\subseteq \TT^{0,\ff} \\ \star\in (\TT\overset{C}{\ominus}\cherryS)^{0,\ff}}} (-1)^{|C|}
= \sum_{C\subseteq \TT^{0,\ff}\setminus\{\star\}} (-1)^{|C|}
= \ind_{\TT^{0,\ff}=\{\star\}}.
\]

\noindent\underline{Case 2: $(\star1,\star2)\in\TT^{0,\ff}\times\TT^\infty$ (or the symmetric case).} Here $\star\in(\TT\overset{C}{\ominus}\<1>)^{0,\ff}$ requires $\star1\in C$, so the alternating sum collapses to $-\ind_{\TT^{0,\ff}=\{\star1\}}$ (resp.\ $-\ind_{\TT^{0,\ff}=\{\star2\}}$).

\noindent\underline{Case 3: $\star1,\star2\in\TT^{0,\ff}$.} Now $\star\in(\TT\overset{C}{\ominus}\<1>)^{0,\ff}$ requires $\{\star1,\star2\}\subseteq C$, yielding $\ind_{\TT^{0,\ff}=\{\star1,\star2\}}$.

Assembling the three cases gives
\begin{align}
\label{eq:mulgE-nv}
\wt{\mul}_{\gfrak,E}(\TT)
&= (-1)^j E \sum_{\star\in \TT^0} \ind_{A_\star^c}K_\star\Psi_E(\star) \bigg(\prod_{\bul\in \TT^0\setminus\{\star\}}K_\bul\ind_{A_\bul^c}\bigg) \nonumber\\
&\quad \times\Big( \ind_{\TT^{0,\ff}=\{\star\}} - \ind_{\TT^{0,\ff}=\{\star1\}} - \ind_{\TT^{0,\ff}=\{\star2\}} + \ind_{\TT^{0,\ff}=\{\star1,\star2\}}\Big).
\end{align}

\medskip

Having exploited the cancellations induced by the combinatorial structure in ~\eqref{eq:mulb-nv},~\eqref{eq:mulg0-nv}, and~\eqref{eq:mulgE-nv}, we now split the nonlinear terms according to the shape of $\TT^{0,\ff}$:
\begin{itemize}
\item if $\TT^{0,\ff}=\{\star\}$ with $\star\neq \emptyset$, only the indices $\star$ and $\pb(\star)$ produce nonzero contributions. At $\star$, $\wt{\mul}_{\gfrak,0}$ and $\wt{\mul}_{\gfrak,E}$ combine to $\mul_2+\mul_3$ (Types~$\II$ and~$\III$). At $\pb(\star)$, only $\wt{\mul}_\bfrak$ contributes, giving $\mul_4$ (Type~$\IV$).
\item if $\TT^{0,\ff}=\{\star1,\star2\}$, only the index $\star$ survives; the three pieces combine to $\mul_5$ (Type~$\5$).
\item otherwise, every alternating sum in~\eqref{eq:mulb-nv}--\eqref{eq:mulgE-nv} vanishes, so $\wt{\mul}(\TT)=0$.
\end{itemize}
Consequently, for every $\TT\in \Tr_j$ with $j\ge 2$ the total contribution carried by $\TT$ equals $\sum_{\ell=2}^5 \Ical(\TT,\mul_\ell(\TT);w)$. Summing over $j\ge 1$ and $\TT\in\Tr_j$ yields~\eqref{wint0a}.
\end{proof}

\begin{remark}\label{REM:tree-size}
At fixed $j\ge 2$, Types~$\II$--$\IV$ contribute $2^{j-1}$ trees of size $j$, and Type~$\5$ contributes $2^{j-1}$ trees of size $j+2$; this is dramatically smaller than $|\Tr_j|\sim 4^j/j^{3/2}$. The decisive feature of Lemma~\ref{LEM:wint}, however, is not the cardinality but the rigid shape of the survivors, which forces a very restricted family of frequency interactions and allows for a systematic derivation of frequency-restricted estimates in Section~\ref{sec:multi}.
\end{remark}

Substituting Lemmas~\ref{LEM:wbd} and~\ref{LEM:wint} into~\eqref{w1} produces the 
normal form equation for~$w$:
\begin{align}
\label{w2}
\dt w & =  \sum_{j=1}^\infty \sum_{\TT \in \Tr_j} \sum_{\l=1}^5  \NN_\l(\TT;w).
\end{align}
Setting
\begin{align}
\label{z0}
z(t) := U(t) w(t),
\end{align}
with $U(t)$ the linear Novikov--Veselov propagator, a direct calculation shows that $z$ obeys the \emph{gauged Novikov--Veselov equation}
\begin{align}
\dt z + (\partial^3+\overline{\partial}^3)z + E(\overline{\partial}^{-1}\partial^2 + \partial^{-1}\overline{\partial}^2)z
&= \sum_{j=1}^\infty \sum_{\TT \in \Tr_j} \sum_{\l=1}^5  U(t) \NN_\l(\TT;U(-t)z)
\nonumber\\
&=: \sum_{j=1}^\infty \sum_{\TT \in \Tr_j} \sum_{\l=1}^5  \Ncal_\l(\TT;z).\label{z1}
\end{align}
Equation~\eqref{z1} is what we call the \textit{gauged} NV equation, clarifying the meaning of the nonlinearity in \eqref{gnv}. This equation is the focus of the next two sections: Section~\ref{sec:multi} establishes multilinear $X^{s,b}$-bounds for each $\Ncal_\l(\TT;z)$ that are summable in $j$ thanks to the rigid shape of the surviving trees, and Section~\ref{sec:lwp} deduces local well-posedness of~\eqref{z1} for $H^s(\R^2)$ data throughout the full subcritical range, before transferring the result back to (NV) via the gauge map.

\section{Functional setting and multilinear estimates}\label{sec:multi}

Following the conventions of \cite{bourg1, bourg2, tao_book}, the Fourier restriction spaces tailored to the linearised flow are defined as follows. For exponents $s,b\in\R$, set
\begin{equation}\label{xsb}
X^{s,b}=\left\{ u\in \S'(\R\times \R^2) : \jap{\tau-\phi(\xi)}^b\jap{\xi}^s  \widehat{u}\in L^2_{\tau,\xi}(\R\times \R^2)\right\},
\end{equation}
equipped with the induced norm, where the phase function $\phi$ is given by
\begin{align}\label{phi}
\phi(\xi): = \Re(\xi^3) - E \Re(\tfrac{\xi^3}{|\xi|^2}),\quad \xi\in\R^2\simeq \CC.
\end{align}
 Its restriction to a finite-time interval is, as usual,
\begin{equation}
X^{s,b}_T=\left\{ u\in \mathcal{D}'((0,T)\times \R^2) : u\equiv v\big|_{(0,T)}\mbox{ for some }v\in X^{s,b}  \right\},\quad T>0.
\end{equation}

The local well-posedness of the gauged NV equation \eqref{z1} follows from a standard contraction mapping argument in $X^{s,b}_T$ from the following multilinear estimates. 

\begin{proposition}\label{prop:bourgL2}
Let $-1 < s \le 0 $, $E\in\R$ with $|E|<1$, $0<\delta\ll 1$ and $T>0$.
Then, there exists $b>\frac12$, such that for all $j \in\N$, $\TT \in \Tr_j$, and $\Ncal_\l(\TT)$, $\l=1,\dots,5$, as in \eqref{z1} (see also Lemma~\ref{LEM:wint}), the estimate holds
	\begin{align}\label{eq:bourg_bi_l2}
		\left\|\int_0^t U(t-t') \Ncal_\l(\TT; z_1,\dots, z_{j+1})(t') dt'  \right\|_{X^{s,b}_T}  \le T^{0+} C^j\delta^{(s+1)(j-1)-3}\prod_{k=1}^{j+1}\|z_k\|_{X^{s,b}_T}
	\end{align}
    for some constant $C>0$ independent of $j$ and $\TT$.
\end{proposition}

\begin{remark}\label{rmk:exponent_optimal}
For $-1 < s \le 0 $, the seven steps in the proof of Proposition \ref{prop:bourgL2} actually produce the  bounds $\delta^{-2+s}$ in Step 1, for $j=1$; $\delta^{-3 + (s+1)(j-1)}$ in Steps~2-3, for $j\ge 2$; 
$1$ in Step~4, with $j=2$; 
$\dl^{(s+1)(j-3)}$ in Step~5, with $j\ge3$;
$\dl^{s+1-}$ in Step~6, with $j=3$;
and $\delta^{\frac12(s+1)(2j+1)}$ in Steps~7, for $j\ge 4$.
 The uniform exponent in \eqref{eq:bourg_bi_l2} is convenient for the geometric-series argument in Proposition~\ref{prop:lwp_rnv}, whose summability is guaranteed by~\eqref{eq:escolhadelta}.
\end{remark}

We show Proposition~\ref{prop:bourgL2} by a reduction to frequency-restricted estimates. Below, we collect the relevant auxiliary results, 
both the bilinear (Lemma~\ref{lem:fre_bi}) and higher order multilinear versions (Lemma~\ref{lem:fre_tri}).

\begin{lemma}[{\cite[Lemma 3.1]{CLS}}]\label{lem:fre_bi}
Let $s\in\R$, $m_1$ be as in \eqref{eq:multiplicadores}, and $\Psi(\<1>)$ as in \eqref{PsiT}.
Suppose that, for all $\al\in \R$ and $M>1$, there exist $0<\beta<1$ and $C>0$ such that one of the following holds:
\begin{enumerate}
    \item[(i)] for all distinct $ k, \l\in\{\emptyset,1,2\}$\footnote{The empty-set symbol $\emptyset$ here corresponds to the frequency component $\xi$, without index.}
    and $1<M\lesssim |\alpha|$,
    \begin{equation}\label{eq:frequad}
        \sup_{\xi_k}\int_{\xi=\xi_1+\xi_2}
        \frac{|m_1(\<1>)|^2\jap{\xi}^{2s}}{\jap{\xi_{1}}^{2s}\jap{\xi_{2}}^{2s}}\mathbbm{1}_{|\Psi(\cherryS)-\alpha|<M} d\xi_{\l} \le C\jap{\alpha}^\beta M^\beta;
    \end{equation}

    \item[(ii)] we have that
    \begin{equation}\label{eq:frequad3}
        \sup_{\xi}\int_{\xi=\xi_1+\xi_2}\frac{|m_1(\<1>)|^2\jap{\xi}^{2s}}{\jap{\xi_{1}}^{2s}\jap{\xi_{2}}^{2s}}\mathbbm{1}_{|\Psi(\cherryS)-\alpha|<M} d\xi_{1}  \le C\jap{\alpha}^\beta M;
    \end{equation}

	\item[(iii)] there exist distinct $k,\l \in \{1,2\}$ such that
		\begin{equation}\label{eq:frequad2}
        \sup_{\xi_k}\int_{\xi=\xi_1+\xi_2}\frac{|m_1(\<1>)|^2\jap{\xi}^{2s}}{\jap{\xi_{1}}^{2s}\jap{\xi_{2}}^{2s}}\mathbbm{1}_{|\Psi(\cherryS)-\alpha|<M} d\xi_{\l}  \le C\jap{\alpha}  M.
    \end{equation}
\end{enumerate}
Then, there exists $b>\frac12$ such that \eqref{eq:bourg_bi_l2} holds for $\ell=1$ and $\TT=\<1>$ (with implicit constant $C^{\frac12}$).
\end{lemma}

\begin{remark}\label{remark:frebi} A close inspection of the proof of \cite[Lemma 3.1]{CLS} shows that, in item (i), for a fixed $k\in\{\emptyset, 1, 2\}$, one may replace \eqref{eq:frequad} with
\begin{equation}\label{eq:frequad_2d}
\sup_{\zeta_\ell}\int_{\zeta=\zeta_1+\zeta_2} \sup_{\substack{\eta_k, \alpha, M \\1\le M\lesssim |\alpha|}}
        \left(\int \frac{|m_1(\<1>)|^2\jap{\xi}^{2s}}{\jap{\xi_{1}}^{2s}\jap{\xi_{2}}^{2s}}\mathbbm{1}_{|\Psi(\cherryS)-\alpha|<M} d\eta_{\l}\cdot \jap{\alpha}^{-\beta} M^{-\beta}\right) d\zeta_k \le C, 
\end{equation}
where $\xi_j = (\eta_j, \zeta_j)$, $j\in\{\emptyset, 1, 2\}$. 
Heuristically, in the duality form of \eqref{eq:bourg_bi_l2}, one fixes the integration in the $\zeta$ variables, applies the frequency-restricted reduction to the $\eta$ variables, and then performs a Cauchy-Schwarz inequality to bound the $\zeta$-integrals. This refinement, available in two dimensions, will be useful in Appendix~\ref{sec:s-34}.

\end{remark}

\begin{lemma}[{\cite[Lemma 2]{COS}}]\label{lem:fre_tri}
Let $s\in\R$, $j\in \N$, $\ell=2,\dots,5$, and $\TT\in \Tr_j$. Let $B\subset\TT^\infty \cup\{\emptyset\}$ be a proper non-empty subset, and suppose that $\M_1, \M_2:\CC^{j+1} \to \R^+$ are multipliers satisfying
\begin{align}
\frac{|m_\ell(\TT)|^2\jap{\xi}^{2s}}{\prod_{\bul\in\TT^\infty} \jap{\xi_{\bul}}^{2s}}
\les
\M_1(\vec\xi)\,\M_2(\vec\xi),
\end{align}
where $\vec\xi=(\xi_\bul)_{\bul\in\TT^\infty}$.
If there exist $0<\be<1$ and $C_1,C_2>0$ such that
\begin{align}
\label{L2hypo}
\sup_{\xi_{\bul\in B}} \int_{\Gamma(\TT)} \mathcal{M}_1 (\vec\xi)\, \mathbbm{1}_{|\Psi(\TT)-\alpha|<M} \,d\xi_{\bul\notin B}
&
\le C_1M^\beta
,
\\
\label{L2hypo2}
\sup_{\xi_{\bul\notin B}} \int_{\Gamma(\TT)} \mathcal{M}_2 (\vec\xi)\, \mathbbm{1}_{|\Psi(\TT)-\alpha|<M} \,d\xi_{\bul\in B}
&
\le C_2 M^\beta,
\end{align}
for all $M\ge 1$ and $\al\in\R$, where $\Psi(\TT)$ is as in \eqref{PsiT}, 
then there exists $b>\frac12$, depending only on $\be$, such that \eqref{eq:bourg_bi_l2} holds (with implicit constant $(C_1C_2)^{\frac12}$). 
\end{lemma}

The endpoint $\beta=1$ in the hypotheses of Lemmas~\ref{lem:fre_bi}-\ref{lem:fre_tri} is not allowed, but it can be reached at the price of a small loss in the largest frequency. The following technical statement formalises this trade-off and will be invoked tacitly in every step of the proof.

\begin{lemma}[{\cite[Lemma 3.3]{CLS}}]\label{lem:eta1}
    Let $K:\R^n\times \R^m\times \R\to \R^+$ be a positive measurable function such that, for some $N\in \mathbb{N}$,
    $$
    |K(x,y,M)|\lesssim \max\{1,|x|,|y|\}^N,\quad \forall x\in\R^n,\ y\in \R^m, \ M\in \R.
    $$
    Suppose that
    $$
    \sup_y \int K(x,y,M)\max\{1,|x|,|y|\}^{0+} dx \lesssim M,\quad \mbox{for all }M>1.
    $$
    Then, there exists $0<\beta<1$ such that
    $$
    \sup_y \int K(x,y,M) dx \lesssim M^{\beta},\quad \mbox{for all } M>1.
    $$
\end{lemma}

We turn next to the geometric heart of the proof: the behaviour of the resonance function $\Psi(\TT)$ in \eqref{PsiT} on the support of the tree integral. The dichotomy below (non-stationary in $\xi_1$ off a thin set, Morse with parameters on it) is precisely the analytic content driving the optimal Strichartz estimates obtained in~\cite{kazeykinamunoz1, kazeykinamunoz2} for the NV propagator.

\begin{lemma}\label{lem:prop_resonance} 
Fix $|E|\ll 1$, $0<\delta'\ll 1$ and $k\ge 2$. Consider the mapping $\Theta:\CC^{k}\to \R$
\begin{equation}\label{resonance0}
\Theta(\xi_1,\dots, \xi_{k})= \sum_{j=1}^{k+1} \phi(\xi_j)\quad \text{and} \quad \sum_{j=1}^{k+1} \xi_j=0, 
\end{equation}
with $\phi$ as in \eqref{phi}. 

\noi{\rm (i)} Over the region
\begin{equation}
V=\big\{(\xi_1,\ldots, \xi_{k})\in \CC^{k} :\  \tfrac{1}{\delta'}<\max\{|\xi_1|, |\xi_{k+1}|\},\  |\xi_1\pm\xi_{k+1}|>\delta'\max\{|\xi_1|, |\xi_{k+1}|\}\big\},
\end{equation}
we have
$$
|\nabla_{\xi_1}\Theta|  
\ge 2(\delta')^2\max\{|\xi_1|^2, |\xi_{k+1}|^2\}.
$$

\noi{\rm (ii)} For a fixed\footnote{In the applications below, $\mu$ will be a small but fixed universal constant, e.g., $\mu=10^{-10}$.}  $0<\mu<1$,consider the region
$$
U=\big\{(\xi_1,\ldots, \xi_{k})\in \CC^{k} :\  \tfrac{1}{\delta'}<|\xi_1|,\ \mu|\xi_2|< |\xi_1|< \tfrac{1}{\mu}|\xi_2|,\ 
|\xi_1-\xi_{k+1}|<\delta' |\xi_1|\big\}.
$$
Then there exist $0 < \dl' < 1 $, open sets $U_1, U_2$ with $U_1\cup U_2=U$ and changes of coordinates $\varphi_1, \varphi_2$ satisfying, for each $\ell=1,2$,
\begin{align} \label{eq:nondeg}
\varphi_\ell:U_\ell&\mapsto B_{\delta'}(0)\times \CC^{k-1}, \qquad \text{with} \qquad 
|\det D^2_{\xi_1}\varphi_\ell| \sim |\xi_2|^{-2},\quad \ell=1,2,
\\ \vec{\xi}&\mapsto (q_1, \xi_2, \dots, \xi_k),
\end{align}
and
$$
\Theta(\xi_1, \ldots, \xi_k)=|\xi_2|^3\Re q_1^2 - \frac14\phi_0\bigg(\sum_{j=2}^k \xi_j\bigg)+ \sum_{j=2}^{k} \phi(\xi_j).
$$

\end{lemma}

\begin{proof}
We first split the dispersion $\phi$ in \eqref{phi} into energy and energy-free parts
\begin{align}
\phi(\xi) =   - E \Re(\tfrac{\xi^3}{|\xi|^2}) +  \Re\xi^3 = : E\phi_E(\xi) + \phi_0(\xi). 
\end{align}
Fixing $\xi_2, \ldots, \xi_k$, from \eqref{resonance0}, we write $\xi_{k+1}$ as a function of $\xi_1$, $\xi_{k+1}  = -\xi_1 -\sum_{j=2}^k \xi_j$. 
On the region $V$, a direct computation yields
\begin{align}
&|\nabla_{\xi_1} (\phi_0(\xi_1)+\phi_0(\xi_{k+1}))| = 3|\xi_1^2-\xi_{k+1}^2|\ge 3(\delta')^2\max\{|\xi_1|^2,|\xi_{k+1}|^2\},\\ &|\nabla_{\xi_1} (\phi_E(\xi_1)+\phi_E(\xi_{k+1}))|\le C.
\end{align}
Combined with the lower bound $\max\{|\xi_1|, |\xi_{k+1}|\}>1/\delta'$, this gives
$$
|\nabla_{\xi_1} \Theta|\ge 3 (\delta')^2\max\{|\xi_1|^2,|\xi_{k+1}|^2\} - C|E| \ge \tfrac12 (\delta')^2\max\{|\xi_1|^2,|\xi_{k+1}|^2\},
$$
for $|E| \le \tfrac{1}{ C }$.

We now turn to (ii).  Define the rescaled variables
$$
p_j=\frac{\xi_j}{|\xi_2|},\quad j=1,\dots, k+1,\quad \text{and }c:=p_1+p_{k+1}=-\sum_{j=2}^{k}p_j
$$
which put $\Theta$ into the form
\begin{align}\label{eq:Theta_rescaled}
\Theta(\xi_1, \dots, \xi_k)&=|\xi_2|^3\bigg(  \sum_{j=1}^{k+1} \phi_0(p_j) - \tfrac{E}{|\xi_2|^2} \sum_{j=1}^{k+1} \phi_E(p_j)\bigg)\\
&=|\xi_2|^3\Big(\phi_0(p_1) + \phi_0(p_{k+1}) - \tfrac{1}{4}\phi_0(c) - \tfrac{E}{|\xi_2|^2}\big(\phi_E(p_1)+\phi_E(p_{k+1})\big) \Big)
\\&\ \quad + \tfrac{1}{4}\phi_0\bigg(\sum_{j=2}^k \xi_j\bigg)+ \sum_{j=2}^{k} \phi(\xi_j).
\end{align}
Over the set $U$,
$$
\mu\le |p_1|\le \tfrac{1}{\mu},\quad |p_{k+1}|\le \tfrac{1}{\mu}(1+\delta'),\quad |p_1+p_{k+1}|>(2-\delta')\mu.
$$
In particular, $c$ takes values on a compact annulus $C  = \{ (2-\dl') \mu \le |c| \le (2+\dl')\frac1\mu \} \subset \CC\setminus\{0\}$. We split $U=U_1\cup U_2$ with
$$
U_1= \{\vec{\xi}\in U: \Re c < |\text{Im} c| \},\quad U_2= \{\vec{\xi}\in U: \Re c >-|\text{Im} c| \}.
$$
Over each $U_\ell$, $\ell=1,2$, choosing appropriately the branch cut, we may define $p_1'=c^{\frac{1}{2}}\left(p_1-\frac{c}{2}\right)$. A short manipulation yields
\begin{align}\label{resonance1a}
\phi_0(p_1)+\phi_0(p_{k+1})-\tfrac{1}{4}\phi_0(c)=3\Re (p_1')^2,\quad |p_1'|<\tfrac{\delta'}{2\mu}\sqrt{\tfrac{2+\delta'}{\mu}}.
\end{align}

To handle the lower-order energy term, it is convenient to introduce the parameter
$$
\varepsilon=\tfrac{E}{|\xi_2|^2},\qquad |\varepsilon|<\tfrac{|E| (\dl')^2}{\mu^2}.
$$
 Combining \eqref{eq:Theta_rescaled}, \eqref{resonance1a} and the choice of $\eps$, the function 
$$
F(p_1', c, \varepsilon) = 3\Re (p_1')^2 -  \varepsilon \bigl(\phi_E(p_1) + \phi_E(c-p_1)\bigr), 
$$
where $ p_1 = p_1(p_1', c) = \frac{c}{2} + c^{-\frac12 } p'_1$,
becomes the object on which the Morse normal form is to be applied. Pick a reference point $c_0\in C$. At $\varepsilon=0$, the map $p_1'\mapsto F(p_1', c_0, 0)$ has an isolated critical point at $p_1'=0$ whose real Hessian has signature $(1,-1)$.
Invoking the Morse lemma with parameters, in the form available in  (\cite[Lemma~C.6.1]{hormander1}; see also \cite[Section~3.4]{duistermaat}), we obtain radii $\eta(c_0),\eta'(c_0)>0$ and a smooth family $\varphi_{c,\varepsilon}:B_{\eta(c_0)}(0)\to B_{\eta'(c_0)}(0)$ of diffeomorphisms $\varphi_{c,\varepsilon}(p_1')=q_1$ such that, for $|c-c_0|+|\varepsilon|<\eta(c_0)$,
$$
F(p_1', c, \varepsilon) = \Re q_1^2,\qquad |\det D^2_{p_1'}\varphi_{c,\varepsilon}|\sim 1.
$$
The next step is to make the radii independent of $c$ via compactness. The family $\{|c-c_0|<\eta(c_0)\}_{c_0\in C}$ covers the compact annulus $C$, so we can extract a finite subcover indexed by $c_0^{(1)},\dots,c_0^{(N)}$ and define
$$
\delta'=\min \big( \min_{1\le j\le N} \eta(c_0^{(j)}),\min_{1\le j\le N} \eta'(c_0^{(j)}), 1 \big) .
$$
Then, the constants depend on $\mu$ and $E$ (through the upper bound $|E|\ll1$) but not on the residual variables $\xi_2,\dots,\xi_k$ nor on the parameter~$k$.

Going back to the original coordinates, we rewrite \eqref{eq:Theta_rescaled} as 
\begin{align}
\Theta (\vec \xi )
=
 |\xi_2|^3 \Re q_1^2 +
\tfrac{1}{4} \phi_0 \bigg( \sum_{j=2}^k \xi_j \bigg) + \sum_{j=2}^k \phi(\xi_j).
\end{align}
The Jacobian \eqref{eq:nondeg} of the resulting map $\varphi_\l:\xi_1\mapsto q_1$ follows from chaining the rescaling $\xi_1=|\xi_2|p_1$ with the unit-Jacobian estimate $|\det D^2_{p_1'}\varphi_{c,\varepsilon}|\sim 1$.
\end{proof}

For convenience, we record the pointwise estimates for the multipliers $\mul_\ell(\TT)$, $\ell=1,\dots,5$, of Lemma \ref{LEM:wint} in the following lemma, whose proof is immediate.

\begin{lemma}\label{lem:pointwise_multipliers}
Given $\Psi$, $\Psi_0$, $\Psi_E$, and $K$ as in \eqref{Psi-star}, \eqref{psi-circ}, \eqref{psi-circ-E}, and \eqref{Kmul}, respectively, and $\xi,\xi_1,\xi_2\in\CC$ with $\xi=\xi_1+\xi_2$, we have
$$
|\Psi_0(\<1>)|\sim |\xi\xi_1\xi_2|,\quad  |\Psi_E(\<1>)|\lesssim \min_{j=\emptyset, 1,2}|\xi_j|,\quad |K|\lesssim \frac{1}{|\xi_1||\xi_2|}.
$$
Moreover, if $|\xi|\gtrsim 1$, then $|\Psi(\<1>)|\lesssim |\xi\xi_1\xi_2|$.
\end{lemma}

\begin{proof}[Proof of Proposition \ref{prop:bourgL2}]

We only present the proof for $-1 < s < -\frac12$, as the proof for $s \ge -\frac12$ follows by a simple modification.
The seven steps that follow are arranged primarily by tree type $\I$-$\5$ and, within each type, by the length $j$ of the tree.

\medskip

\noi\textit{Step 1. Type~I, $j=1$, $\TT=\<1>$: workhorse Lemma~\ref{lem:fre_bi}.} The bilinear estimate at $j=1$ is the only one not subsumed by Lemma~\ref{lem:fre_tri}; we split the multiplier $\mul_1(\<1>)$ in \eqref{eq:multiplicadores} into two pieces and feed  each into Lemma~\ref{lem:fre_bi}. 
Using Lemma~\ref{lem:pointwise_multipliers}, we have
$$
|\mul_1(\<1>)|\lesssim |\xi|\ind_A + | K | |\Psi_E(\<1>)|\ind_{A^c} \lesssim |\xi|\ind_A + \frac{\ind_{A^c}}{\max_{\bul\in\{\emptyset, 1, 2\}}|\xi_\bul|}.
$$
Without loss of generality, assume $|\xi_1|\ge |\xi_2|$.  
For the first term, if $|\xi_1| \les \frac1\dl$, then $|\xi|, |\xi_2| \les \frac1\dl$ and we get
\begin{align*}
\sup_\xi \int_A \frac{|\xi|^2\jap{\xi}^{2s}}{\jap{\xi_1}^{2s}\jap{\xi_2}^{2s}} \ind_{|\Psi(\cherryS)-\alpha|<M} d\xi_2
\les \dl^{-2 + 2s} \int_{|\xi_2| \le 1/\dl} d\xi_2 \les \dl^{-4 + 2s}, 
\end{align*}
since $s>-1$.
Else, $|\xi_2| \le \frac1\dl \les |\xi_1| \sim |\xi|$, and thus, from Lemma~\ref{lem:prop_resonance}(i), $|\nabla_{\xi_2} \Psi(\<1>) | \ges |\xi|^2$. Performing the change of variables\footnote{We implicitly assumed that $|\partial_{\Re \xi_2} \Psi(\<1>)| \ges |\nabla_{\xi_2} \Psi(\<1>) |$, which allows for the change of variables in the first component of $\xi_2$. Otherwise, we have $|\partial_{\Im \xi_2} \Psi(\<1>) | \ges |\nabla_{\xi_2} \Psi(\<1>) | $ and perform the change of variables $\xi_2 \mapsto (\zeta_2, \phi)$ instead. In the following, we assume that the former holds and omit details.} $\xi_2 \mapsto (\phi, \zeta_2)$ with $\phi = \Psi(\<1>)$ and $\zeta_2 = \Im \xi_2$,  
\begin{align*}
\sup_\xi \int_A \frac{|\xi|^{2}\jap{\xi}^{2s}}{\jap{\xi_1}^{2s}\jap{\xi_2}^{2s}} \ind_{|\Psi(\cherryS)-\alpha|<M} d\xi_2 
\lesssim
\dl^{2s} \sup_\xi \int_{|\zeta_2|<1/\delta}  \ind_{|\phi-\alpha|<M} \,d\phi\, d\zeta_2 \lesssim \delta^{-1 + 2s}M.
\end{align*}

\noi
and Lemma~\ref{lem:fre_bi}(ii) applies.

For the second term, we have $|\xi_1| \ge |\xi_2| \ge \frac1\dl$. If $\xi_1 \not\simeq \pm \xi_2 $, then $|\xi| \sim |\xi_1|$ and use Lemma~\ref{lem:prop_resonance}(i) to obtain
$$
\sup_\xi \int_{A^c} \frac{\jap{\xi}^{2s}}{\jap{\xi_1}^{2s+2}\jap{\xi_2}^{2s}} \ind_{|\Psi(\cherryS)-\alpha|<M} d\xi_2 \lesssim \sup_\xi \int_{|\zeta_2| > 1/\dl} \frac{1}{|\zeta_2|^{4+2s}} \ind_{|\phi-\alpha|<M}\,d\phi\, d\zeta_2 \lesssim \dl^{3 + 2s } M, 
$$
and Lemma~\ref{lem:fre_bi}(ii)  applies once more. 
Otherwise, $\xi_1 \simeq \xi_2$ or $\xi_1 \simeq -\xi_2$, in which case, for fixed $\xi_1$, $|\nabla_{\xi} \Psi(\<1>) | \sim |\xi_1|^2$, and  the change of variables $\xi \mapsto (\phi, \zeta)$ gives
\begin{align}
\sup_{\xi_1} \int_{A^c} \frac{\jap{\xi}^{2s}}{\jap{\xi_1}^{2s+2}\jap{\xi_2}^{2s}} \ind_{|\Psi(\cherryS) - \al| < M } d\xi 
&
\lesssim \sup_{\xi_1} \int_{A^c} \frac{\jap{\zeta}^{2s}}{\jap{\xi_1}^{4s+4}} \ind_{|\phi-\al| < M } d\phi d\zeta \lesssim \delta^{4s+4}M.
\end{align}
The result now follows from Lemma~\ref{lem:fre_bi}(iii).

\smallskip 

For Types II--V (i.e., $\ell=2,\dots,5$) the workhorse switches from Lemma~\ref{lem:fre_bi} to Lemma~\ref{lem:fre_tri}, whose flexibility in choosing the splitting set $B$ is essential.  We consider the multiplier with an extra $(0+)$-power of the largest frequency and aim to prove a bound by $M$, so that Lemma~\ref{lem:eta1} then gives the desired \eqref{L2hypo}-\eqref{L2hypo2}.
For all the trees $\TT\in \Tr_j$ of each type, by symmetry, we only consider the representative tree which always grows to the left, namely such that all $\bul \in \TT^\infty$ with $\pb(\bul)\not \in \TT^{0,\ff}$ are right-children, satisfying $\bul = \pb(\bul)2$. 

\smallskip

\noi\textit{Step 2. Type~II, $j\ge 2$: workhorse Lemma~\ref{lem:fre_tri} with $B=\{\star1, \star 2\}$.}
Assume $|\xi_{\star1}| \ge |\xi_{\star2}|$. Due to the restriction to $A_\star$, this implies that $|\xi_{\star2}|\le 1/\delta$.  
We have, for $\TT\in \Tr_j$, $j\ge 2$,
\begin{align}
|\mul_2(\TT)| \frac{\jap{\xi}^s }{\prod_{\bul\in\TT^\infty}\jap{\xi_\bul}^s} &\lesssim \frac{|K_{\pb(\star)}K_\star \Psi_{ 0}(\star)|\ind_{A_\star }\jap{\xi}^s}{\jap{\xi_{\star1}}^s\jap{\xi_{\star2}}^s\jap{\xi_{\sbf(\star)}}^s}\prod_{\substack{\bul\in\TT^0\\\bul\neq \star, \pb(\star)}} \frac{|K_\bul|\ind_{A_\bul^c}}{\jap{\xi_{\bul2}}^{s}}\\&
\lesssim \frac{\delta^{s}
\jap{\xi}^s }{\jap{\xi_{\star1}}^{s}\jap{\xi_{\sbf(\star)}}^{s+1}}\prod_{\substack{\bul\in\TT^0\\\bul\neq \star, \pb(\star)}} \frac{\ind_{A_\bul^c}}{\jap{\xi_{\bul1}}\jap{\xi_{\bul2}}^{s+1}}, \label{eq:estimativa_mul2}
\end{align}
where we used \eqref{wint0c} and Lemma~\ref{lem:pointwise_multipliers}.
The worst-case scenario occurs when $\xi_{\star1}$ corresponds to the largest frequency (the specific tree structure plays no further role).

For $\xi$ and $\xi_{\bul}$, $\bul\in \TT^\infty\setminus\{\star1, \star2\}$, fixed, if $|\xi_1| \les \frac1\dl$, we integrate directly, while for $|\xi_1| \gg \frac1\dl$, we use Lemma~\ref{lem:prop_resonance}(i) and perform the change of variables $\xi_{\star2}\mapsto (\phi, \zeta_{\star2})$ with $\phi = \Psi(\TT)$ and $\zeta_{\star2} = \Im \xi_{\star2}$:
\begin{equation}\label{type2j}
\int |\xi_{\star1}|^2\ind_{|\Psi(\TT)-\alpha|<M} \,d\xi_{\star2} \lesssim   \dl^{-4} + \int_{|\zeta_{\star2}|\le 1/\dl} \ind_{|\phi-\alpha|<M} \,d\phi\, d\zeta_{\star2} \lesssim  \dl^{-4} + \delta^{-1}M
\les \dl^{-4}M. 
\end{equation}
On the other hand, for $\xi_{\star1},\xi_{\star2}$ fixed, we do not use the resonance restriction: 
\begin{align*}
& \dl^{2s}\int \frac{\jap{\xi}^{2s}}{\jap{\xi_{\star1}}^{2s+2-}   \jap{\xi_{\sbf(\star)}}^{2s+2}}
\bigg( \prod_{\substack{\bul\in\TT^0\\\bul\neq \star, \pb(\star)}} \frac{\ind_{A_\bul^c} d\xi_{\bul2}}{\jap{\xi_{\bul1}}^2\jap{\xi_{\bul2}}^{2s+2}}\bigg)
\ind_{|\Psi(\TT)-\alpha|<M}\, d\xi\\
\lesssim\  &  \dl^{2s} \int \left(\int \frac{\jap{\xi}^{2s}}{\jap{\xi_{\star1}}^{2s+2-}   \jap{\xi_{\sbf(\star)}}^{2s+2}}d\xi\right)\prod_{\substack{\bul\in\TT^0\\\bul\neq \star, \pb(\star)}} \frac{\ind_{A_\bul^c} d\xi_{\bul2}}{\jap{\xi_{\bul1}}^2\jap{\xi_{\bul2}}^{2s+2}}\\
\lesssim {}
&
\dl^{(2s+2)(j-2) + 2s } 
.
\end{align*}

\noi\textit{Step 3. Type~III, $j\ge 2$: reduction to Type~II (Step~2).}
\\
We observe that, since $|\Psi_E(\star)|\lesssim |\xi_\star|$,
\begin{align}
|\mul_3(\TT)|\frac{\jap{\xi}^s}{\prod_{\bul\in\TT^\infty}\jap{\xi_\bul}^s}  &\lesssim \frac{|EK_{\pb(\star)}K_\star \Psi_E(\star)|\ind_{A_\star^c}\jap{\xi}^s}{\jap{\xi_{\star1}}^s\jap{\xi_{\star2}}^s\jap{\xi_{\sbf(\star)}}^s}\prod_{\substack{\bul\in\TT^0\\\bul\neq \star, \pb(\star)}} \frac{|K_\bul|\ind_{A_\bul^c}}{\jap{\xi_{\bul2}}^{s}}
\\
&\lesssim \frac{ \jap{\xi}^s\ind_{A_\star^c}}{\jap{\xi_{\star1}}^{s+1}\jap{\xi_{\star2}}^{s+1}\jap{\xi_{\sbf(\star)}}^{s+1}}\prod_{\substack{\bul\in\TT^0\\\bul\neq \star, \pb(\star)}} \frac{|K_\bul|\ind_{A_\bul^c}}{\jap{\xi_{\bul2}}^{s}}
\\
&\lesssim 
\frac{ \jap{\xi}^s}{\jap{\xi_{\star1}}^{s}\jap{\xi_{\star2}} \jap{\xi_{\sbf(\star)}}^{s+1}}
\prod_{\substack{\bul\in\TT^0\\\bul\neq \star, \pb(\star)}} \frac{\ind_{A_\bul^c}}{\jap{\xi_{\bul1}}\jap{\xi_{\bul2}}^{s+1}},
\end{align}
which is precisely the multiplier appearing in \eqref{eq:estimativa_mul2} with the extra factor $\jb{\xi_{\star2}}^{-1}$. As such, the frequency-restricted estimates for Type III reduce to those shown in Step 2, where we use the extra factor in $\xi_{\star2}$ as a replacement for its smallness in \eqref{type2j}.

\medskip

\noi \textit{Step 4. Type~IV, $j=2$: workhorse Lemma~\ref{lem:fre_tri} with $B=\{2,11\}$, supported by Lemma~\ref{lem:prop_resonance}.}
With $j=2$ the only representative tree is $\TT=\<21>$ and
\begin{equation}
\label{typeIVmult2}
|\mul_4(\TT)|\frac{\jap{\xi}^s}{\jap{\xi_{11}}^s\jap{\xi_{12}}^s\jap{\xi_2}^s} \lesssim \frac{\jap{\xi}^{s+1}}{\jap{\xi_{11}}^{s+1}\jap{\xi_{12}}^{s+1}\jap{\xi_2}^s}.
\end{equation}
The worst-case scenario corresponds to $|\xi_2|\ge |\xi|\ge |\xi_{11}|\ge |\xi_{12}|$, as the multiplier in the remaining cases is controlled by the above. Moreover, we can always assume that the two largest frequencies satisfy $|\xi_2| \ge |\xi| \ge \frac1\dl$.

The proof proceeds by bounding
\begin{align}\label{eq:type4interp1}
& 
\sup_{\xi,\xi_{12},\alpha} \int \frac{|\xi_2|^{2}}{\jb{\xi_{11}}^{1+}} \ind_{|\Psi(\TT)-\alpha|<M}\,d\xi_{11}
\\
\label{eq:type4interp2}
& \sup_{\xi_2,\xi_{11},\alpha} \int \frac{\jb{\xi}^{2s+2+}}{\jap{\xi_{11}}^{2s+1-}\jap{\xi_{12}}^{2s+2} \jap{\xi_2}^{2s+2}}\ind_{|\Psi(\TT)-\alpha|<M}\,d\xi_{12}.
\end{align}

\smallskip
\noi\textbf{Case A.}
{$(\xi\not\simeq \pm \xi_{12}, \xi_2\not\simeq \pm \xi_{11}) \lor (\xi\not\simeq \pm \xi_{11}, \xi_2\not\simeq \pm \xi_{12})$.}

\noi 
Under the assumption $\xi\not\simeq \pm \xi_{12}, \xi_2\not\simeq \pm \xi_{11}$, Lemma~\ref{lem:prop_resonance} applies both for $(\xi,\xi_{12})$ or $(\xi_2,\xi_{11})$ fixed, to perform the change of variables in $\xi_{11}$ or $\xi_{12}$, respectively. 
For $\xi, \xi_{12}$ fixed,
\begin{equation}\label{eq:typeIVest1}
\eqref{eq:type4interp1}\lesssim \int \frac{1}{\jb{\zeta_{11}}^{1+}}\ind_{|\phi-\alpha|<M}\,d\phi\, d\zeta_{11} \lesssim M.
\end{equation}
For $\xi_2, \xi_{11}$ fixed,
\begin{align}\label{eq:typeIVest2}
\eqref{eq:type4interp2} 
\lesssim 
\int 
\frac{1}{\jb{\zeta_{12}}^{2s+2} \jb{\xi_{11}}^{2s+3-} }
\ind_{|\Psi(\TT)-\alpha|<M}\,d\phi\, d\zeta_{12} \lesssim M.
\end{align}

\noi 
The same strategy applies if $\xi \not\simeq \pm \xi_{11}$ and $\xi_2 \not\simeq \pm \xi_{12}$, by swapping the roles of $\xi_{11}$ and $\xi_{12}$ above.

\smallskip
\noi
\textbf{Case B.} $\lnot (\xi\not\simeq \pm \xi_{12}, \xi_2 \not\simeq \pm \xi_{11}) \land \lnot (\xi\not \simeq \pm \xi_{11} , \xi_2 \not\simeq \pm \xi_{12})$. 

In this case, all frequencies are comparable, so that
$$
|\mul_4(\TT)|\lesssim |\xi|^{-1-2s}.
$$

\noi
Moreover, 
three of the frequencies are equal up to a sign, say $\xi_2\simeq  \xi_{11} \simeq \ld \xi_{12}$,  for $\ld \in \{\pm 1\}$. If the sign is positive, then
$$
\xi_2\simeq \xi_{11}\simeq \xi_{12} \simeq \tfrac{\xi}{3}.
$$
For $\xi_2, \xi_{11}$ fixed, we estimate \eqref{eq:type4interp2} as in Case A (the non-stationarity condition still holds). For $\xi, \xi_{12}$ fixed, we apply Lemma \ref{lem:prop_resonance}(ii)
with $k=3$ and $(\xi_1, \xi_2, \xi_3, \xi_4) := (\xi_{11},  -\xi, \xi_{12}, \xi_{2})$,
to obtain
\begin{align}
\eqref{eq:type4interp1}\lesssim \sup_{\xi, \wt\al} \int_{|q_1|\ll 1} |\xi|^{3-}\ind_{||\xi|^3\Re q_1^2 - \widetilde\alpha|<M} \,dq_1 \lesssim M^{1-}, 
\end{align}
where we used that an integral over $\{q \in \CC^2: |q| \le 1, \, |\Re q^2 -\al|< M  \} $ is bounded by $M^{1-}$; see \cite[Lemma 5]{COS} for a proof.

\noi
Finally, if $\xi_2\simeq \xi_{11}\simeq-\xi_{12}$, then $\xi\simeq -\xi_{12}$. In this case, both \eqref{eq:type4interp1} and \eqref{eq:type4interp2} are handled by Lemma \ref{lem:prop_resonance}(ii) as above for \eqref{eq:type4interp1} and with the choice $(\xi_1,\xi_2,\xi_3,\xi_4) := (\xi_{12}, \xi_2, \xi_{11}, -\xi) $ for \eqref{eq:type4interp2}.

\medskip

\noi \textit{Step 5. Type~IV, $j\ge 3$: workhorse Lemma~\ref{lem:fre_tri}.}

\noi 
Set $\circ=\gb(\star)$. Then
\begin{equation}\label{eq:typeIVmult}
|\mul_4(\TT)|\frac{\jap{\xi}^s}{\prod_{\bul\in\TT^\infty}\jap{\xi_\bul}^s} \lesssim \frac{\jap{\xi}^{s}}{\jap{\xi_{\circ111}}^{s+1}\jap{\xi_{\circ112}}^{s+1}\jap{\xi_{\circ12}}^{s}\jap{\xi_{\circ2}}^{s+1}} \prod_{\bul\in\TT^\infty, \bul\ngeq\circ} \frac{1}{\jap{\xi_\bul}^{s+1}\jap{\xi_{\sbf(\bul)}}}.
\end{equation}
We consider the worst case in terms of frequency relations, namely that  $|\xi_{\circ12}|\ge |\xi_{\circ 2}|\ge |\xi_{\circ 111}| \ge |\xi_{\circ 112}| \ge |\xi|$. 
Moreover, we can assume that $|\xi_{\circ12}|\sim |\xi_{\circ 2}|$. Note that this is a valid assumption since the roles of the frequencies $\xi_{\circ111}, \xi_{\circ112}, \xi_{\circ2}, \xi_\bul$ for $\bul\in\TT^\infty$ with $\bul\not\ge \circ$ are interchangeable (as they are all terminal nodes with the same power in \eqref{eq:typeIVmult}). In fact, we can assume that $|\xi_{\circ2}|  \ge |\xi_{\circ111}| \ge |\xi_{\circ112}| \ge \frac1\dl$ are the 3 largest frequencies in this group, which in turn implies that $|\xi_{\circ12}|\sim |\xi_{\circ2}|$, since otherwise $|\xi_{\circ12}|\gg |\xi_{\circ2}|$ and $|\xi_{\circ12}| \sim |\xi| \le |\xi_{\circ2}| \ll |\xi_{\circ12}|$ which cannot happen.

Consider
\begin{align}
\label{eq:IVinterp1}
&\sup_{\xi_{\circ111}, \xi_{\circ12}}
\int \frac{\jap{\xi}^{2s}}{|\xi_{\circ12}|^{2+3s-}}\left(\prod_{\bul\in\TT^\infty, \bul\ngeq\circ} \frac{d\xi_\bul}{\jap{\xi_\bul}^{2s+2}\jap{\xi_{\sbf(\bul)}}^2}\right) \ind_{|\Psi(\TT)-\alpha|<M}\,d\xi\, d\xi_{\circ112}, 
\\
\label{eq:IVinterp2}
&\sup_{\xi, \xi_{\circ2}, \xi_{\circ112}, \xi_{\bul\in\TT^\infty, \bul\ngeq\circ}}
\int \frac{|\xi_{\circ12}|^{2+s+}}{\jap{\xi_{\circ111}}^{2s+2}\jap{\xi_{\circ112}}^{2s+2}\jap{\xi_{\circ2}}^{2s+2}}\ind_{|\Psi(\TT)-\alpha|<M} \,d\xi_{\circ111}.
\end{align}

\smallskip
\noi\textbf{Case A.}  ($\xi_{\circ12}\not\simeq \pm \xi_{\circ111}$, $\xi_{\circ2}\not\simeq \pm \xi_{\circ112}$)$\lor$($\xi_{\circ12}\not\simeq \pm \xi_{\circ112}$, $\xi_{\circ2}\not\simeq \pm \xi_{\circ111}$)

\noi 
We consider only the first condition, as the same approach also applies if the second condition holds. For $\xi_{\circ12}, \xi_{\circ111}$ fixed, from Lemma~\ref{lem:prop_resonance}~(i) with the change of variables $\xi_{\circ112} \mapsto (\phi, \zeta_{\circ112}) = (\Psi(\TT), \Im \xi_{\circ112})$, we get
\begin{align*}
\eqref{eq:IVinterp1}\lesssim & \int \frac{\jap{\xi}^{2s}}{|\xi_{\circ12}|^{4+3s-}}\left(\prod_{\bul\in\TT^\infty, \bul\ngeq\circ} \frac{d\xi_\bul}{\jap{\xi_\bul}^{2s+2}\jap{\xi_{\sbf(\bul)}}^2}\right) \ind_{|\phi-\alpha|<M}\,d\xi\, d\zeta_{\circ112} d\phi\lesssim \delta^{(2s+2)(j-3)}M,
\end{align*}
where the $j-3$ factors of $\delta^{2s+2}$ from the $|\xi_{\bul 2}|\gtrsim 1/\delta$ constraints associated with each of the $j-3$ terminal vertices above $\circ$ in the tree. 
Similarly, for $\xi, \xi_{\circ2}, \xi_{\circ112}$ and $\xi_{\bul\in \TT^\infty, \bul\ngeq \circ}$ fixed,
\begin{align*}
\eqref{eq:IVinterp2} \lesssim & \int \frac{1}{\jap{\xi_{\circ111}}^{2s+2} \jb{\xi_{\circ112}}^{2s+2} \jb{\xi_{\circ2}}^{s+2-}}\ind_{|\phi-\alpha|<M}\, d\phi\, d\zeta_{\circ111} \lesssim M.
\end{align*}

\smallskip
\noi
\textbf{Case B.} 
$\lnot(\xi_{\circ12}\not\simeq \pm \xi_{\circ111}, \xi_{\circ2}\not\simeq \pm \xi_{\circ112})\land\lnot(\xi_{\circ12}\not\simeq \pm \xi_{\circ112},\xi_{\circ2}\not\simeq \pm \xi_{\circ111}).$

\noi
In this context, we may suppose, without loss of generality, that $|\xi_{\circ111}|\sim |\xi_{\circ12}| \sim |\xi_{\circ2}|$ and that one of the following holds:
$$
\begin{cases}
\xi_{\circ111}\simeq \xi_{\circ12} \simeq  \lambda\xi_{\circ2} \not\simeq \pm \xi_{\circ112}\ \text{ for some }\lambda=\pm1 & \text{(Case B1)}\\
\xi_{\circ111}\simeq \xi_{\circ12},\  \xi_{\circ2} \simeq \xi_{\circ112} & \text{(Case B2)}\\
\xi_{\circ111}\simeq \xi_{\circ12},\  \xi_{\circ2} \simeq -\xi_{\circ112} & \text{(Case B3)}
\end{cases}
$$

In all cases, for $\xi, \xi_{\circ2}, \xi_{\circ112}$ and $\xi_{\bul\in \TT^\infty, \bul\ngeq \circ}$ fixed, we apply Lemma \ref{lem:prop_resonance}(ii) to find
\begin{align*}
\eqref{eq:IVinterp2} \lesssim & \int_{|q_{\circ111}|<1}\jap{\xi_{\circ2}}^{-3s+}\ind_{||\xi_{\circ2}|^{3}\Re q_{\circ111}^2-\tilde\alpha|<M}\, dq_{\circ111} \lesssim M.
\end{align*}
We now focus on \eqref{eq:IVinterp1}. 
In Case B1, we may proceed as in Case A, as $|\nabla_{\xi_{\circ112}} \Psi(\TT)|\gtrsim |\xi_{\circ12}|^2$. In Case B2, we apply instead Lemma \ref{lem:prop_resonance}(ii):
\begin{align*}
\eqref{eq:IVinterp1}\lesssim & \sup_{\wt\al, \xi, \xi_{\bul\in \TT^\infty, \bul\not\geq\circ}} \left(\int_{|q_{\circ112}|<1} |\xi_{\circ12}|^{-3s+} \ind_{||\xi_{\circ12}|^{3}\Re q_{\circ112}^2-\tilde\alpha|<M}dq_{\circ112}\right)\\ &\times  \int \jap{\xi}^{2s}\left(\prod_{\bul\in\TT^\infty, \bul\ngeq\circ} \frac{d\xi_\bul}{\jap{\xi_\bul}^{2s+2}\jap{\xi_{\sbf(\bul)}}^2}\right) \,d\xi\, \lesssim \delta^{(2s+2)(j-3)}M.
\end{align*}
Finally, in Case B3, $|\xi_{\circ}|\sim |\xi_{\circ111} + \xi_{\circ112} + \xi_{\circ12} +\xi_{\circ2}| \sim  |\xi_{\circ111}|$. If $\circ \neq \emptyset$, assume that $\xi_{\sbf(\circ)}\not\simeq -\xi_{\circ2}$ (otherwise exchange the roles of $\xi_{\circ2}$ and $\xi_{\circ112}$,  since $\xi_{\sbf(\circ)} \simeq -\xi_{\circ2} \simeq \xi_{112}$, and thus $\xi_{\sbf(\circ)} \not\simeq - \xi_{\circ112}$). Then 
\begin{align*}
\eqref{eq:IVinterp1}\lesssim & \sup_{\xi, \xi_{\bul},\bul\in \TT^\infty, \bul\not\geq\pb(\circ)} \left(\int |\xi_{\circ 111}|^{-4-5s+} \jb{\xi_{\sbf(\circ)}}^{-2} \ind_{|\Psi(\TT)-\alpha|<M}d\xi_{\circ2} d\xi_{\circ{112}}\right)\\ &\times  \int \jap{\xi}^{2s}\left(\prod_{\bul\in\TT^\infty, \bul\ngeq\pb(\circ)} \frac{d\xi_\bul}{\jap{\xi_\bul}^{2s+2}\jap{\xi_{\sbf(\bul)}}^2}\right) \,d\xi\lesssim \delta^{(2s+2)(j-3)}M,
\end{align*}

\noi 
either via Lemma \ref{lem:prop_resonance}(i) (if $\xi_{\sbf(\circ)}\not\simeq \pm \xi_{\circ2}$) or  Lemma \ref{lem:prop_resonance}(ii) (if $\xi_{\sbf(\circ)}\simeq \xi_{\circ2}$). Finally, if $\circ=\emptyset$, then $\xi\not\simeq \pm \xi_{\circ2}$ and, performing a change of variables $\xi_{\circ2}\mapsto (\phi,\zeta_{\circ2})$, $\phi=\Psi(\TT)$, gives
\begin{align*}
\eqref{eq:IVinterp1}\lesssim & \int |\xi_{\circ111}|^{-2-s+} \ind_{|\Psi(\TT)-\alpha|<M}d\xi_{\circ2} d\xi_{\circ{112}}
\\
\lesssim & \int |\xi_{\circ111}|^{-4-s+} \ind_{|\phi-\alpha|<M}d\phi d\zeta_{\circ2} d\xi_{\circ{112}}
\lesssim  \delta^{s+1 -}M.
\end{align*}

\medskip
\noi\textit{Step 6. Type~V, $\ell=5$, $j=3$, $\TT=\<33>$: workhorse Lemma~\ref{lem:fre_tri}.}
At $j=3$, Type~V is realized by the single tree $\TT=\<33>$, whose multiplier involves four leaves rather than three. 
With $\TT=\<33>$, from \eqref{wint0c}
$$
|\mul_5(\TT)| \frac{\jb{\xi}^s}{\prod_{\bul\in\TT^\infty} \jb{\xi_\bul}^s} \lesssim \frac{\jap{\xi}^{s+1}}{\jap{\xi_{11}}^{s+1}\jap{\xi_{12}}^{s+1}\jap{\xi_{21}}^{s+1}\jap{\xi_{22}}^{s+1}}.
$$
The worst-case scenario occurs when $|\xi|\ge |\xi_{11}|\ge |\xi_{12}|\ge |\xi_{21}|\ge |\xi_{22}|  \ge \frac1\dl$, when we estimate the two quantities:
\begin{align}
\label{eq:Vinterp1}
\sup_{\xi, \xi_{12}, \xi_{22}}
&
\int |\xi|^{1-} \ind_{|\Psi(\TT)-\alpha|<M} \,d\xi_{21}
,
\\
\label{eq:Vinterp2}
\sup_{\xi_{11}, \xi_{21}}
&
\int \frac{ 1 }{|\xi|^{1-}\jap{\xi_{12}}^{2s+2}\jap{\xi_{21}}^{2s+2}\jap{\xi_{22}}^{2s+2}}\ind_{|\Psi(\TT)-\alpha|<M} \,d\xi_{12}\,d\xi_{22}.
\end{align}

\smallskip

\noi\textbf{Case A.} $(\xi\not\simeq \pm \xi_{\l_1}, \xi_{11}\not\simeq \pm\xi_{\l_2})$ for distinct $\l_1, \l_2 \in \{ 12, 21, 22\}$.

\noi  
We assume $\l_1 = 12, \l_2 = 21$, as the remaining cases can be handled analogously, with suitable modifications to the choice of integration variables in \eqref{eq:Vinterp1}-\eqref{eq:Vinterp2}. 
For $\xi_{11}, \xi_{21}, \xi_{22}$ fixed, Lemma~\ref{lem:prop_resonance}~(i) gives
\begin{equation}\label{eq:estVinterp1}
\eqref{eq:Vinterp1} \lesssim \int \frac{1}{|\xi|^{1+}} \ind_{|\phi-\alpha|<M} \,d\phi\,  d\zeta_{12} \lesssim M.
\end{equation}
For $\xi_{11}, \xi_{21}$ fixed,
\begin{equation}\label{eq:estVinterp2}
\eqref{eq:Vinterp2} \lesssim \int \frac{1}{|\xi_{11}|^{3-}\jap{\zeta_{21}}^{4s+4}\jap{\xi_{22}}^{2s+2}} \ind_{|\phi-\alpha|<M}\, d\phi\, d\zeta_{21}\,d\xi_{22} \lesssim \dl^{6(s+1)-} M .
\end{equation}

\smallskip
\noi\textbf{Case B.} $\lnot (\xi \not\simeq \pm \xi_{\l_1}, \xi_{11}\not\simeq \pm \xi_{\l_2})$ for all distinct $\l_1, \l_2 \in\{12, 21, 22\}$. 

\noi In this case, we have that at least 4 of the frequencies in $\{\xi, \xi_{11}, \xi_{12}, \xi_{21}, \xi_{22}\}$ agree up to signs, including either $\xi$ or $\xi_{11}$. In particular, $|\xi| \sim |\xi_{11}| \sim |\xi_{12}|\sim |\xi_{21}| \Gr > \frac1\dl$ and{\Gr, for $s>-1$,} the multiplier is bounded by
\begin{align}
|\mul_5(\TT)| 
\frac{\jb{\xi^s}}{\prod_{\bul\in\TT^\infty} \jb{\xi_\bul}^s} \les \frac{1}{\jb{\xi}^{2s+2} \jb{\xi_{22}}^{s+1}} \les \frac{1}{\jb{\xi}^{2s+2}}.
\end{align}
Notice that, since $\xi = \xi_{11}+ \xi_{12} + \xi_{21} +\xi_{22}$, we cannot have all five frequencies agreeing up to a sign.

\noi
Hence we may suppose, by symmetry, that $\xi \simeq \ld_{11} \xi_{11} \simeq \ld_{12} \xi_{12} \simeq \ld_{21} \xi_{21}$, for $\ld_{11}, \ld_{12}, \ld_{21} \in \{\pm1\}$, and $\xi_{11}\not\simeq \pm \xi_{22}$.
If $\xi \simeq - \xi_{11}$, since $\xi_{21}\not\simeq \pm \xi_{22}$, we can proceed by Lemma~\ref{lem:prop_resonance}~(i)-(ii) for the analogues of \eqref{eq:Vinterp1}-\eqref{eq:Vinterp2} as follows
\begin{align}
\sup_{\xi, \xi_{11},\xi_{21}} \int \frac{1}{|\xi|^{1+}} \ind_{|\Psi(\TT) - \al | < M } d\xi_{22}  
& \les M
, 
\\
\sup_{\xi_{12}, \xi_{22}} \int \jb{\xi}^{-(2s+2)} \ind_{|\Psi(\TT) - \al| < M } d \xi d\xi_{21}  
&
\les \int \int_{|q| < 1 } \jb{\xi_{21}}^{-2s+} \ind_{| |\xi_{21}|^3 \Re q^2 - 
\wt\al | < M } dq d\xi_{21}  
\\
&
\les M \int  \jb{\xi_{21}}^{-2s - 3 + } d\xi_{21}
\les M \dl^{2(s+1)-}. 
\end{align}
Similar estimates hold if $\xi \simeq - \xi_{12}$ or $\xi \simeq - \xi_{21}$. Otherwise, $\xi \simeq \xi_{11} \simeq \xi_{12} \simeq \xi_{21}$, where we proceed as above, where for the second integral we fix $\xi_{21}, \xi_{22}$ and integrate in $\xi_{11}, \xi$, using $\xi_{11} \simeq \xi_{12}$ to apply Lemma~\ref{lem:prop_resonance}~(ii) to the inner integral.

\medskip
\noi\textit{Step 7. Type~V, $\ell=5$, $j\ge 4$: workhorse Lemma~\ref{lem:fre_tri}.}

\noi
From \eqref{wint0c},  the multiplier is 
\begin{align}\label{eq:mul5_high_j}
& 
|\mul_5(\TT)|  \frac{\jb{\xi}^{s}}{\prod_{\bul\in\TT^\infty} \jb{\xi_{\bul}}^s }   
\\
&\lesssim \frac{\jap{\xi}^s}{\jap{\xi_{\star11}}^{s+1}\jap{\xi_{\star12}}^{s+1}\jap{\xi_{\star21}}^{s+1}\jap{\xi_{\star22}}^{s+1}\jap{\xi_{\sbf(\star)}}^{s+1}} \prod_{\bul \in \TT^\infty, \bul\ngeq \pb(\star)} \frac{1}{\jap{\xi_{\sbf(\bul)}}\jap{\xi_\bul}^{s+1}}.
\end{align}
In the worst-case scenario, $|\xi_{\star11}|\gtrsim|\xi_{\star12}|\gtrsim|\xi_{\star21}|\gtrsim|\xi_{\star22}|\gtrsim|\xi_{\sbf(\star)}|\gtrsim |\xi|$. 
We consider
\begin{align}
\label{eq:Vinterp3}
&\sup_{\xi_{\star11}, \xi_{\star21}, \xi_{\sbf(\star)}} 
 \int \frac{\jap{\xi}^{2s}}{|\xi_{\star11}|^{2+3s-}}\prod_{\bul \in \TT^\infty, \bul\ngeq \pb(\star)} \frac{d\xi_\bul}{\jap{\xi_{\sbf(\bul)}}^2\jap{\xi_\bul}^{2s+2}}\,\ind_{|\Psi(\TT)-\alpha|<M}\,d\xi\, d\xi_{\star22}, 
\\
\label{eq:Vinterp4}
&\sup_{\xi, \xi_{\star12}, \xi_{\star22}, \xi_{\bul \in \TT^\infty, \bul\ngeq \pb(\star)}} 
\int \frac{1}{|\xi_{\star11}|^{-s \Gr +}\jap{\xi_{\star21}}^{4s+4}\jap{\xi_{\sbf(\star)}}^{4s+4}} \ind_{|\Psi(\TT)-\alpha|<M}\, d\xi_{\star21}\,d\xi_{\sbf(\star)}.
\end{align}

\noi
The quantity \eqref{eq:Vinterp3} is essentially identical to \eqref{eq:IVinterp1} of Step~5 (with an extra harmless fixed variable) and the estimate for \eqref{eq:Vinterp3} follows the same lines.

\noi
\textbf{Case A. $\xi_{\star11} \not\simeq \pm \xi_{\sbf(\star)}$}
In this case,
we perform the change of variables $\xi_{\sbf(\star)}\mapsto (\phi, \zeta_{\sbf(\star)})$, $\phi=\Psi(\TT)$ (Lemma \ref{lem:prop_resonance}~(i)), to find
$$
\eqref{eq:Vinterp4}\lesssim \int\frac{1}{|\xi_{\star11}|^{-s+2+}\jap{\xi_{\star21}}^{4s+4}\jap{\zeta_{\sbf(\star)}}^{4s+4}} \ind_{|\phi-\alpha|<M}\, d\phi d\zeta_{\sbf(\star)} d\xi_{\star21} \lesssim \delta^{7(1+s)}M.
$$

\noi 
\textbf{Case B.  $\xi_{\star11} \simeq \ld  \xi_{\sbf(\star)}$, $\ld \in \{\pm 1 \}$}
\\
Here, we have $|\xi_{\star11}| \sim |\xi_{\star12}|\sim|\xi_{\star21}| \sim |\xi_{\star22}| \sim |\xi_{\sbf(\star)}|$.
If $\xi_{\star11} \not\simeq \pm \xi_{\star21}$ or $\xi_{\star21} \not\simeq \pm \xi_{\sbf(\star)}$, we can proceed as in Case A, with the correct permutation of the variables $\xi_{\star11}, \xi_{\star 21}, \xi_{\sbf(\star)}$. Thus, it remains to consider $\ld_{\star11} \xi_{\star11} \simeq \ld_{\star21} \xi_{\star21} \simeq \ld_{\sbf(\star)} \xi_{\sbf(\star)}$, for some $\ld_\bul \in \{\pm1\}$, $\bul\in \{\star11, \star21, \sbf(\star)\}$. In particular, at least 2 of the $\ld_\bul$ must agree, so we can assume that $\xi_{\star11} \simeq \xi_{\sbf(\star)}$, and by Lemma~\ref{lem:prop_resonance}~(ii), we find
$$
\eqref{eq:Vinterp4}\lesssim \int \int_{|q_{\sbf(\star)}|<1 } \frac{1}{|\xi_{\star12}|^{7s+6+}} \ind_{||\xi_{\star12}|^3\Re q_{\sbf(\star)}^2-\tilde\alpha|<M}\,  dq_{\sbf(\star)}  d\xi_{\star21}\,\lesssim \delta^{7(1+s)}M.
$$

\medskip

Assembling Steps~1--7 yields, for every $\TT\in\Tr_j$, $j \ge 2$, and every $\ell\in\{1,\dots,5\}$, the estimate \eqref{eq:bourg_bi_l2} with the uniform exponent $\delta^{(s+1)(j-1)-3}$ (cf.\ Remark~\ref{rmk:exponent_optimal}). The proof is finished.
\end{proof}

\section{Proof of the main theorem}\label{sec:lwp}

\subsection{Proof of Theorem \ref{thm:main_intro}} In this subsection, we prove Theorem \ref{thm:main_intro} assuming assuming Proposition \ref{PRO:equiv}, whose proof will be postponed to Subsection \ref{SUBSEC:equiv}. By a suitable rescaling (see \eqref{scaling}), it suffices to consider the case $|E|\ll1$. Throughout this section, we fix $0< \eps\ll1$ and $R\gg 1$, and take $\delta>0$ such that
\begin{equation}\label{eq:escolhadelta}
\delta^{s+1}R\ll \eps.
\end{equation}

We first show well-posedness of \eqref{gnv} (see also \eqref{z1}). 

\begin{proposition}[Local well-posedness for the gauged NV equation]\label{prop:lwp_rnv}
Fix $s>-1$. Given $z_0\in H^s(\R^2)$ with $\|z_0\|_{H^s}<R$, there exist $T=T_{\delta,R}>0$ and a unique
$
z\in X^{s,\frac12 +}_T
$
solution of \eqref{z1} over $[0,T)$ with initial data $z_0$.
\end{proposition}

\begin{proof}
The argument is by Banach fixed point in the standard way; we record the main inequality in order to highlight the dependence on $\delta$.

Fix $b=\frac12+$ as in Proposition~\ref{prop:bourgL2}
and define
$$
B_{R,T}=\left\{ z\in X^{s,b}_T : \|z\|_{X^{s,b}_T}\le 2R \right\},
$$
endowed with the induced distance, together with the integral map
$$
F[z](t): =U(t)z_0 + \int_0^t U(t-t')\sum_{j=1}^\infty \sum_{\TT \in \Tr_j} \sum_{\l=1}^5  \Ncal_\l(\TT;z)(t')\,dt',
$$
see \eqref{z1} and \eqref{wint0a} for the definition of the nonlinear terms.
Applying the multilinear estimates of Proposition \ref{prop:bourgL2} and recalling that $|\Tr_j|$ grows only exponentially in $j$ (cf.~\eqref{eq:catalan}), for any $z\in B_{R,T}$ we obtain
\begin{align}
\|F[z]\|_{X^{s,b}_T}
&\le \|z_0\|_{H^s} + T^{0+}\sum_{j=1}^\infty C^j\delta^{(s+1)(j-1)-3}\|z\|_{X^{s,b}_T }^{j+1} \\
&\le R + T^{0+} C\delta^{-3} R^2 \sum_{j=1}^\infty \big(C \,\delta^{s+1}\, R\big)^{j-1}  \\
&\le R + T^{0+} C \delta^{-3} R^2 \le 2R,
\end{align}
where the geometric series converges thanks to the smallness condition \eqref{eq:escolhadelta}, and $T=T_{\delta,R}$ is chosen sufficiently small in the last inequality. Using the multilinear structure of the nonlinear terms, an analogous computation yields
$$
\|F[z]-F[z']\|_{X^{s,b}_T } \le T^{0+}C\delta^{-3}R\|z-z'\|_{X^{s,b}_T } \le \tfrac12\|z-z'\|_{X^{s,b}_T }.
$$
The result follows from Banach's fixed point theorem.
\end{proof}

\begin{remark}\label{remark:persistence}
If $z_0\in H^{s'}(\R^2)$, for some $s'>s$, a standard persistence-of-regularity argument shows that $z\in X^{s',b}(0,T) \hookrightarrow C([0,T]; H^{s'}(\R^2))$, where $T$ depends only on $R=\|z_0\|_{H^s}$ instead of $R' = \|z_0\|_{H^{s'}}$.
\end{remark}

Next, we define rigorously the gauge transform, whose (formal) definition is
\begin{equation}\label{eq:gauge}
\mathcal{G}_\delta^{-1}[f]=f+ \F_{x,y}^{-1}\left(\int_{\xi=\xi_1+\xi_2} K(\xi_1,\xi_2) \mathbbm{1}_{A^c}(\xi,\xi_1,\xi_2)\widehat{f}(\xi_1)\widehat{f}(\xi_2)d\xi_1\right)
=: 
f + \DD[f, f]
, 
\end{equation}
with $K$ and $A^c$ as in \eqref{Kmul} and \eqref{A}.

\begin{lemma}\label{LEM:D}

Consider the map  $\Gc^{-1}_\dl$ defined in \eqref{eq:gauge} and let $0<\eps<  1$. 

\smallskip

\noi {\rm(i)}
Let $s' \ge s>-1$, $R>0$, and $B^{s',s}_R$ be given by
\begin{align}
\label{Bss}
B^{s',s}_R=\left\{ f\in H^{s'}(\R^2): \|f\|_{H^s} \le R \right\}.
\end{align}
Then, there exists $\dl = \dl(\eps,R)$ with $0< \dl \ll R^{-1/(1+s)}$ such that the map $\Gc_\dl^{-1}: B^{s',s}_R \to B^{s',s}_{(1+\eps)R}$ is continuous with a continuous inverse $\Gc_\dl: B^{s',s}_{(1-\eps)R} \to B^{s',s}_R$.

\smallskip
\noi {\rm(ii)} The map $\Gc_\dl^{-1}$ is unbounded on $H^s(\R^2)$ for $s< - 1$ and $\dl>0$.

\end{lemma}

\begin{proof}
We proceed as in the proof of \cite[Lemma 6.1]{chcora}. 
We first show boundedness of the bilinear operator $\DD$ in \eqref{eq:gauge}: by duality, Lemma~\ref{lem:pointwise_multipliers}, and Cauchy-Schwarz, we have
\begin{align}
\| \DD[f_1,f_2] \|_{H^s} 
& 
\le  \sup_{\|g\|_{H^{-s}}\le 1} \int_{\xi, (\xi_1,\xi_2)\in A^c} \frac{|\ft{g}(\xi) \ft{f}_1(\xi_1) \ft{f}_2(\xi_2)| }{ |\xi_1 \xi_2|} d\xi_1 d\xi
\\
&
\le \sup_{\xi_1} \bigg( \int  \frac{\jb{\xi}^{2s} \ind_{A^c}}{\jb{\xi_1}^{2s+2}   \jb{\xi_2}^{2s+2}  } d\xi \bigg)^\frac12 \| f_1\|_{H^s} \| f_2\|_{H^s}
\\&\le C \dl^{s+1} \| f_1\|_{H^s} \| f_2\|_{H^s}, 
\label{Dbd}
\end{align}
for $s>-1$ and some $C>0$ independent of $s$.

\smallskip 

Given $0<\eps<1$, from \eqref{Dbd}, for $s=s'$, $f_1,f_2 \in B^{s,s}_R$, and $2 C \dl^{s+1} R \le \eps$ , we have
\begin{align}
\| \DD[f_1,f_1] \|_{H^s} & \le C \dl^{s+1} R^2 \le \eps R 
, 
\\
\| \DD[f_1,f_1] - \DD[f_2,f_2] \|_{H^s} & \le 2 C \dl^{s+1} R \| f_1-f_2\|_{H^s} \le \eps \|f_1-f_2\|_{H^s}, 
\end{align}
showing that $\DD$ is Lipshitz from $B^{s,s}_R$ to $B^{s,s}_{\eps R}$ with constant $\eps$. Then, mapping properties of $\Gc^{-1}_\dl$ follow from being a small Lipschitz perturbation of the identity and the existence, continuity, and boundedness of $\Gc_\dl$ follow from a fixed point argument. 
This shows (i) for $s'=s$. 

To update the above result to $s'>s$, it suffices to show boundedness and continuity of the maps $\Gc_\dl^{-1}$ and $\Gc_\dl$ in $H^{s'}(\R^2)$. For $\Gc^{-1}_\dl$, these properties follow from the estimate
\begin{align}
\| \DD[f,f] \|_{H^{s'}} \le C \dl^{s+1} \|f\|_{H^{s'}} \|f\|_{H^{s'}}, 
\end{align}
for $s'\ge s > -1$, which can be shown by modifying \eqref{Dbd}.
It remains to show the analogous result for $\Gc_\dl$. Via a bilinear von Neumann expansion, we have
\begin{align}\label{vonneumann1}
\Gc_\dl [f] = f + \sum_{j = 1}^\infty \sum_{\TT \in \Tr_j} \Ical(\TT, -i \mul_\BB (\TT); f ), 
\end{align}
where $\mul_\BB$ is as in \eqref{NFE0b} and $\Ical$ in \eqref{eq:defiI}. 
Since $\Ical(\TT, -i \mul_\BB (\TT); f )$ is a $j$-fold composition of $\DD$, from \eqref{Dbd} we get
\begin{align}
\| \Ical(\TT, -i \mul_\BB (\TT); f_1, \ldots, f_{j+1} ) \|_{H^s} \le C^j \dl^{(s+1)j} \prod_{k=1}^{j+1} \| f_k\|_{H^s} , 
\end{align}
which, from $|\xi| \le (j+1) \max_{\bul\in\TT^\infty}|\xi_\bul|$, can be upgraded to
\begin{align}
\| \Ical(\TT, -i \mul_\BB (\TT); f_1, \ldots, f_{j+1} ) \|_{H^{s'}} \le C^j \dl^{(s+1)j} \max_{\l=1, \ldots, j+1} \bigg( \|f_\l\|_{H^{s'}} \prod_{\substack{1\le k \le j+1 \\ k \ne \l }} \|f_k\|_{H^s} \bigg). 
\end{align}
Then, for $f_1,f_2\in B^{s',s}_R$, we obtain
\begin{align}
\| \Gc_\dl[f_1] \|_{H^{s'}} & \le \bigg( 1 + \sum_{j=1}^\infty |\Tr_j| (C \dl^{s+1} R)^j \bigg) \|f_1\|_{H^{s'}} \le (1+\eps) \|f_1\|_{H^{s'}}
, 
\label{vonneumann2a}
\\
\| \Gc_\dl[f_1] - \Gc_\dl[f_2] \|_{H^{s'}} & \le \bigg( 1 + (\|f_1\|_{H^{s'}} + \|f_2\|_{H^{s'}}) \sum_{j=1}^\infty |\Tr_j| (C \dl^{s+1} R)^{j-1}  \bigg) \|f_1 - f_2\|_{H^{s'}} 
\\
& 
\le (1 + \|f_1\|_{H^{s'}} + \|f_2\|_{H^{s'}}) \|f_1 - f_2\|_{H^{s'}} , 
\label{vonneumann2b}
\end{align}
for $\dl,R$ as before. Thus, $\Gc_\dl$ maps $B^{s',s}_R$ to $B^{s',s}_{(1+\eps)R}$ continuously.

To prove (ii), take $\widehat{f}(\xi)=\jap{\xi}^{-1-s-} $, so that $f\in H^s(\R^2)$.
 Then, for $|\xi|<1$, 
$$
\Ft_{x,y} (\DD[f])(\xi) = \int_{\xi=\xi_1+\xi_2} K\ind_{A^c_\dl}(\xi_1,\xi_2)\widehat{f}(\xi_1)\widehat{f}(\xi_2)d\xi_1 \sim \int \frac{1}{\jap{\xi_1}^{4+2s+}}d\xi_1 = \infty,
$$
since $s<-1$, and thus $\DD[f] \not\in H^s(\R^2)$. 
\end{proof}

The last ingredient to complete the proof of Theorem \ref{thm:main_intro} is the following equivalence between the original and gauged flows at high regularity.

\begin{proposition}\label{PRO:equiv}
Let $s' \ge s > -1$ with $s'\gg1$ sufficiently large, $R>0$, $0<\dl \ll R^{-\frac{1}{s+1}}$, and $u_0$ satisfying
\begin{align}\label{equiv-0}
u_0\in H^{s'}(\R^2) \quad \text{and} \quad \|u_0\|_{H^s}<R.
\end{align}
Let $u \in X^{s',b}_T$ be the solution to \eqref{eq:NV} with initial data $u_0$ constructed in \cite{angelopolous, kazeykinamunoz1, kazeykinamunoz2}, where $T>0$ is its maximal time of existence, and let $z\in X^{s', b}_{T_{\delta,R}}$ be the solution to \eqref{z1} with initial data $z_0 := \Gc_\delta[u_0]$ provided by Proposition~\ref{prop:lwp_rnv}. Then, $T \ge T_{\delta,R}$ and
\begin{align}\label{equiv-1}
z(t)=\Gc_\delta[u(t)] \quad \text{in} \quad C([0,T_{\delta,R}]; H^{s'}(\R^2)).
\end{align}
\end{proposition}

We postpone the proof of Proposition~\ref{PRO:equiv} to Subsection~\ref{SUBSEC:equiv} and present the proof of Theorem~\ref{thm:main_intro}.

\begin{proof}[Proof of Theorem~\ref{thm:main_intro}]

Let $u_0 \in H^s(\R^2)$ and $(u_{0,n})_{n\in\N} \subset H^{s'}(\R^2)$ be a smooth approximating sequence with $u_{0,n} \to u_0$ in $H^s(\R^2)$, for $s'\gg1$ sufficiently large, and $\|u_{0,n}\|_{H^s}, \|u_0\|_{H^s} < R$. 
From Proposition~\ref{PRO:equiv}, there exists $T =T_{\dl,R}>0$, independent of $n$, solutions $u_n \in X^{s',b}_T$ to \eqref{eq:NV} (agreeing with those in \cite{angelopolous, kazeykinamunoz1, kazeykinamunoz2} and Appendix~\ref{sec:s-34}) and solutions $z_n \in X^{s',b}_T$ to \eqref{gnv} from Proposition~\ref{prop:lwp_rnv}, with data $u_{0,n}$ and $\Gc_\dl[u_{0,n}]$, respectively, satisfying
\begin{align}\label{WP-aux}
\Theta^{\eqref{gnv}} \circ \Gc_\dl[u_{0,n}]
= z_n = \Gc_\dl [u_n] \quad \text{in} \quad C([0,T]; H^{s'}(\R^2)), 
\end{align}
where $\Theta^{\eqref{gnv}}$ denotes the data-to-solution map for \eqref{gnv}.
From the invertibility of $\Gc_\dl$ in Lemma~\ref{LEM:D}, we define the map
\begin{align}
\Theta^{\eqref{eq:NV}} [f] 
= \Gc_\dl^{-1} \circ \Theta^{\eqref{gnv}} \circ \Gc_\dl [f] , 
\end{align}
and note that \eqref{WP-aux} guarantees that $\Theta^{\eqref{eq:NV}}[u_{0,n}] = u_n$. Thus $\Theta^{\eqref{eq:NV}}$ agrees with the \eqref{eq:NV}-solution map in $C([0,T]; H^{s'}(\R^2))$, it is analytic (from Proposition~\ref{prop:lwp_rnv} and Lemma~\ref{LEM:D}), and defines a unique continuous extension of the solution map to $C([0,T]; H^{s}(\R^2)) $ since
\begin{align}
u_n = \Theta^{\eqref{eq:NV}}[u_{0,n}] \to \Theta^{\eqref{eq:NV}} [u_0] =: u 
\quad \text{in} \quad C([0,T]; H^{s'}(\R^2)). 
\end{align}
\qedhere

\end{proof}

\subsection{Proof of Proposition~\ref{PRO:equiv}}\label{SUBSEC:equiv}

We proceed as in \cite{chcora}, presenting the relevant details needed in this setting where the smallness encoded in the $\dl$-parameter appears in the sets $A^c$ in \eqref{defA} instead of via a small data assumption. See Remark~\ref{REM:small} for further details.

Given $n\in\N$, write $\mathbf{P}_n$ for the Fourier projector with
symbol $\ind_{[0,n]}(|\xi|)$, $\xi\in\R^2$, and consider the
\emph{frequency-truncated Novikov-Veselov equation}
\begin{equation}\tag{NV$_n$}\label{eq:NV_n_dens}
\dt u_n + (\partial^3+\overline\partial^3) u_n
+ \tfrac34 \mathbf{P}_n\!\left( \partial(u_n\overline\partial^{-1}\partial u_n)
+ \overline\partial(u_n\partial^{-1}\overline\partial u_n) \right)
+ E(\overline\partial^{-1}\partial^2 u_n+\partial^{-1}\overline\partial^2 u_n)=0,
\end{equation}
with $u_n(0)=u_0$.
Since $\mathbf{P}_n$ commutes
with the linear group, the resonance set $A$ in \eqref{A} and the phase
$\Psi$ in \eqref{Psi-star} are insensitive to
the cut-off; consequently the normal-form reduction in 
Subsection~\ref{sec:nfe}
carries over to \eqref{eq:NV_n_dens} unchanged. Define
\begin{equation}\label{eq:Kmul_n}
K_{\bul,n} := K_\bul\, \ind_{[0,n]}(|\xi_\bul|),
\qquad \bul\in\TT^0,
\end{equation}
the truncated analogue of $K_\bul$ in \eqref{Kmul}, and let
$\Gc_{\dl,n}$ be the \emph{truncated gauge transform} determined by
\begin{equation}\label{eq:gauge_n_dens}
\Gc_{\dl,n}^{-1}[f] := f + \F_{x,y}^{-1}\!\left( \int_{\xi=\xi_1+\xi_2}
K_n(\xi_1,\xi_2) \mathbbm{1}_{A^c}(\xi, \xi_1,\xi_2)\,
\widehat f(\xi_1)\widehat f(\xi_2)\, d\xi_1 \right),
\end{equation}
i.e.\ \eqref{eq:gauge} with $K$ replaced by $K_n$ in \eqref{eq:Kmul_n}. The conclusions
of Lemma~\ref{LEM:D} hold for $\Gc_{\dl,n}$ with $\dl$, $\eps$, $R$
all uniform in $n$. The corresponding \emph{truncated gauged
equation} is
\begin{equation}\label{eq:GNV_n_dens}\tag{$\Gc_\dl$-NV$_n$}
\dt z_n + (\partial^3+\overline\partial^3)z_n
+ E(\overline\partial^{-1}\partial^2+\partial^{-1}\overline\partial^2)z_n
= \sum_{j=1}^\infty \sum_{\TT\in\Tr_j}\sum_{\l=1}^5
\Ncal_{\l,n}(\TT;z_n),
\end{equation}
where $\Ncal_{\l,n}(\TT)$ arises from $\Ncal_\l(\TT)$ in \eqref{z1}
by substituting each $K_\bul$ in \eqref{eq:multiplicadores} with
$K_{\bul,n}$.

\smallskip 
The multilinear estimates in \eqref{eq:bourg_bi_l2} in
Proposition~\ref{prop:bourgL2} apply verbatim to $\Ncal_{\l,n}$, with constants independent of $n$, by trivially estimating  $\ind_{[0,n]}$ by $1$.
 To establish the equivalence in Proposition~\ref{PRO:equiv}, we also require fixed-time estimates to replace those in Proposition~\ref{prop:bourgL2}, at the cost of additional regularity on one of the factors and a polynomial loss in the truncation parameter $n$, analogous to \cite[Lemma~6.4]{chcora}.

\begin{lemma}\label{lem:dens_estimativas}
Given $s > -1$, $s' \ge s + 4$, $|E|\ll 1$, $j,n\in\N$, $\TT\in \Tr_j$, $\l\in\{1, \ldots, 5\}$, $1\le k \le j+1$, and $f_1, \ldots, f_{j+1} \in H^{s'}(\R^2)$, the following estimate holds
\begin{equation}\label{eq:dens_estimativas}
\|\Ncal_{\l,n}(\TT)[f_1,\dots,f_{j+1}]\|_{H^s}
\le n \, C^j \dl^{(1+s)(j-3) \lor 0} \, \|f_k\|_{H^s}
\max_{\substack{1\le k'\le j+1\\ k'\ne k}}\!\bigg(
\|f_{k'}\|_{H^{s'}}\!\!\!\prod_{\substack{k''\ne k,k'}}\!\!\|f_{k''}\|_{H^s}
\bigg),
\end{equation}
for some $C=C(s, E)>0$ independent of $n$, $j$, $\l$ and $\TT$.
\end{lemma}

\begin{proof}
We present the proof for $s \le 0$, as the $s>0$ case follows from the earlier case together with the trivial bound $\jb{\xi}^s \le (j+1) \max_{\bul\in\TT^\infty} \jb{\xi_\bul}^s$, for $\TT \in \Tr_j$. 
We import the strategy of \cite[Lemma~6.4]{chcora}. Combining \eqref{wint0c},
\eqref{eq:multiplicadores}, Lemma~\ref{lem:pointwise_multipliers}, and \eqref{eq:Kmul_n} produces the
pointwise bound
\begin{equation}\label{eq:mul_n_pointwise}
|\mul_{\l}(\TT)(\vec\xi)| \prod_{\bul\in \TT^0} \ind_{[0,n]}(|\xi_\bul|)
\le C^j\, n\, \max_{\bul\in\TT^\infty}\jb{\xi_\bul}^2\,
\ind_{[0,n]}(|\xi|)
\ind_{\bigcap_{\bul\in\TT^0_{\l, \star}} A^c_\bul}
 \prod_{\bul\in\TT^0}\frac{1}{\jap{\xi_{\bul1}}\jap{\xi_{\bul2}}}
,
\end{equation}
with $ \subset \TT^0$, with $\star$ as in \eqref{wint0c}, is given by
\begin{align*}
\TT^0_{\l,\star} = 
\begin{cases}
\emptyset , & \l =1, 
\\
\TT^0 \setminus \{\star\}, & \l=2,3,5, 
\\
\TT^0 \setminus \{ \pb(\star)\}, & \l=4.
\end{cases}
\end{align*}

\noi The loss of two derivatives in \eqref{eq:mul_n_pointwise} is handled by the $H^{s'}$-factor on the r.h.s. of \eqref{eq:dens_estimativas}, with $s'\ge s+2$, exactly as in \cite{chcora}.

For the contributions to the l.h.s. of \eqref{eq:dens_estimativas} with $\l=1,2,4,5$, which are direct analogues of those in \cite[Lemma 6.4]{chcora} for $\l=1,\ldots,4$, respectively, we run the $p=2$ argument in the proof of \cite[Lemma 6.4]{chcora}.  More precisely, 
we first apply Cauchy-Schwarz inequality, 
followed by 
integration along the 
subtrees associated with $\bul\in\TT^{0,\ff}$ (which requires the additional loss $\max_{\bul\in\TT^\infty}\jb{\xi_{\bul}}^2$, imposing $s' \ge s+4$), 
followed by successively integrating $\jb{\xi_{\bul, \min}}^{-2(2+s)}$, with $|\xi_{\bul, \min}| = |\xi_{\bul1}| \land |\xi_{\bul2}|$
for $\bul\in\TT^0\setminus\TT^{0,\ff}$. The latter integrals for $\bul\in \TT^{0}_{\l, \star}$ contribute with the power $\dl^{(1+ s) (j-3) \lor 0}$ as $|\xi_{\bul, \min}| \ge \frac1\dl$.

For the NV-specific contributions $\Ncal_{3,n}$ associated with the Type $\III$, one 
can easily see from \eqref{eq:multiplicadores} that the convolution structure of this contribution matches that of $\Ncal_{2,n}$, and their multipliers satisfy the same upper bound \eqref{eq:mul_n_pointwise}, and the same strategy applies. This completes the proof. \qedhere

\end{proof}

Following \cite[Lemma~6.5]{chcora}, 
we extend Proposition~\ref{prop:lwp_rnv} to the truncated equation in \eqref{eq:GNV_n_dens} at high regularity $H^{s'}(\R^2)$, with a suitable notion of uniqueness on balls for the weaker norm $H^s(\R^2)$, with $s \le s'$.

\begin{lemma}\label{lem:zn_nv}
Let $|E| \ll 1$, $s > -1$, $s' \ge s+4$, $R>0$, and
$0<\dl\ll R^{ -1/(1+s)}$. For every $n\in\N$ and every
$z_{0,n}\in H^{s'}(\R^2)$ with $\|z_{0,n}\|_{H^s}\le R$, there exists
a unique solution $z_n\in X^{s',b}_{T_{\dl,R}}$ to
\eqref{eq:GNV_n_dens} with $z_n(0)=z_{0,n}$, where $T_{\dl,R}$ is
the local time of existence from Proposition~\ref{prop:lwp_rnv} (independent of $n$).
Moreover, for $t\in[0,T_{\dl,R}]$,
\begin{equation}\label{eq:bounds_zn_nv}
\|z_n(t)\|_{H^s}\le  \| z_n\|_{X^{s, b}_{T_{\dl, R}}} \le 2R,
\qquad
\|z_n(t)\|_{H^{s'}}\le  \| z_n\|_{X^{s', b}_{T_{\dl, R}}} \le 2\|z_{0,n}\|_{H^{s'}},
\end{equation}
and $z_n$ is unique in
$\{z\in C([0,T_{\dl,R}]; H^{s'}(\R^2)):\,
\|z\|_{L^\infty_{T_{\dl, R}} H^s_x}\le 10R\}$.
\end{lemma}

\begin{proof}
The well-posedness of \eqref{eq:gauge_n_dens} in $X^{s,b}_{T_{\dl,R}}$ for all $n\in\N$, for $T_{\dl,R} >0$ independent of $n$, follows from the proof of Proposition~\ref{prop:lwp_rnv}, noting that the estimates in Proposition~\ref{prop:bourgL2} also apply to $\Ncal_{\l,n}$ with constants uniform in $n$. The well-posedness in $X^{s',b}_{T_{\dl,R}}$, for the same time of existence, and the estimates in \eqref{eq:bounds_zn_nv} follow from a standard persistence of regularity argument as described in Remark~\ref{remark:persistence}.

For the sharper uniqueness statement, 
let
$z_n,z_n'$ as in the hypothesis. For $0<T \le T_{\dl, R}$, from the Duhamel formulation for \eqref{z1}, combined with \eqref{eq:dens_estimativas}, Remark~\ref{REM:tree-size}, and \eqref{eq:bounds_zn_nv}, we have 
\begin{align*}
& 
\|z_n-z_n'\|_{L^\infty_{T}H^s_x}
\\
&
\le   T n \big(\|z_n\|_{L^\infty_{T_{\dl,R}} H^{s'}_x}+\|z_n'\|_{L^\infty_{T_{\dl, R}} H^{s'}_x}\big)
\|z_n-z_n'\|_{L^\infty_{T}H^s_x}
\sum_{j=0}^\infty 
(2CR)^j \dl^{(s+1)(j-2) \lor 0} 
\\
&
\le C T R^2 n \big(\|z_n\|_{L^\infty_{T_{\dl,R}} H^{s'}_x}+\|z_n'\|_{L^\infty_{T_{\dl, R}} H^{s'}_x}\big) 
\|z_n-z_n'\|_{L^\infty_{T}H^s_x}, 
\end{align*}
for $0 < \dl \ll R^{-1/(s+1)}$. 
Choosing $T$ sufficiently small depending on $n$, $R$, and $\|z_{n}\|_{L^\infty_{T_{\dl, R}} H^{s'}_x}$, $\|z_{n}'\|_{L^\infty_{T_{\dl, R}} H^{s'}_x}$, we conclude from the above inequality that $z_n = z_n' $ on $[0,T]$. Uniqueness on the interval $[0,T_{\dl,R}]$ follows from iterating the above argument on successive subintervals of size $T$. 
\qedhere

\end{proof}

At the truncated level, we can prove the equivalence, for smooth data, of \eqref{eq:NV_n_dens} and
\eqref{eq:GNV_n_dens} through $\Gc_{\dl,n}$ on $[0,T_{\dl,R}]$.

\begin{lemma}\label{lem:equiv_n_nv}
Let $s'\gg 1$ be sufficiently large, $R>0$,
$0<\dl\ll R^{-1/(1 +s)}$, and $u_0\in H^{s'}(\R^2)$ with
$\|u_0\|_{H^s}<R$. Then,  there exists $T_{\dl,R}>0$ such that for every $n\in\N$:
\begin{enumerate}
\item[\textup{(i)}] 
There is a unique solution
$u_n\in X^{s',b}_{T_{\dl,R}}$ to \eqref{eq:NV_n_dens} on 
$[0,T_{\dl,R}]$, satisfying
\begin{equation}\label{eq:bounds_un}
\|u_n\|_{L^\infty_{T_{\dl,R}}H^s_x}<3R,
\qquad
\|u_n\|_{L^\infty_{T_{\dl,R}} H^{s'}_x}<3\|u_0\|_{H^{s'}};
\end{equation}
\item[\textup{(ii)}] If $z_n$ is the solution to
\eqref{eq:GNV_n_dens} from Lemma~\ref{lem:zn_nv} with
$z_{0,n}:=\Gc_{\dl,n}[u_0]$, then
\begin{equation}\label{eq:equiv_n_nv}
z_n=\Gc_{\dl,n}[u_n]
\quad\text{in}\quad C([0,T_{\dl,R}];H^{s'}(\R^2));
\end{equation}
\item[\textup{(iii)}] We have $\sup_{n\in\N }\|u_n\|_{X^{s',b}_{T_{\dl,R}}}\le
C(\|u_0\|_{H^{s'}})$. 
\end{enumerate}
\end{lemma}

\begin{proof}
The proof retains the three-step skeleton of
\cite[Lemma~6.6]{chcora}; we record only the points where the NV
features enter.

For $s' \gg 1$, from local well-posedness at high regularity (see \cite{angelopolous,kazeykinamunoz1,kazeykinamunoz2, adamsgrunrock} or Appendix~\ref{sec:s-34}), there exists $u_n \in X^{s',b}_{T_{n,R}}$ solution to \eqref{eq:NV_n_dens}, where $T_{n,R}>0$ is the maximal time of existence. 
For the truncated gauged data $z_{0,n} = \Gc_{\dl, n}[u_0] \in H^{s'}(\R^2)$ with $z_{0,n} \in B^{s',s}_{\frac{1}{1-\eps}R} $ for $0<\eps\ll 1$ sufficiently small (see Lemma~\ref{LEM:D}), Lemma~\ref{lem:zn_nv} guarantees the existence of a solution $z_n \in X^{s',b}_{T_{\dl,R}}$ to \eqref{z1} with data $z_{0,n}$, where $T_{\dl,R}>0$ is independent of $n$. 
Lastly,~let 
$$T_{0,n} : = \sup\{ 0< T \le \min(T_{\dl,R}, T_{n,R}) : \|u_n\|_{L^\infty_T H^s_x} \le 3R\}.$$

\smallskip
\noindent\textit{Step~1.}\
Setting $v_n(t):=U(-t)u_n(t)$, we can proceed with the reduction in Subsection~\ref{sec:nfe} for the truncated dynamics in \eqref{eq:NV_n_dens}, to obtain the finite normal form equation analogous to \eqref{NFE1} with $\BB_n,\NN_n,\MM_n$ defined with $K_{\bul,n} = \ind_{|\xi_\bul| \le n} K_\bul$ replacing $K_\bul$ (see \eqref{Kmul}). Note that all the relevant operations are justified  (namely the conditions (a)-(d) in \cite[p.61, proof of Lemma 6.6]{chcora}). Moreover, from the proof of  Lemma~\ref{lem:dens_estimativas}, we obtain summability in $j$ of the terms $\BB_n,\NN_n$ and convergence to 0 of $\MM_n$ as $J\to\infty$ in $L^\infty_{T_{0,n}} H^{s}_x$, and $v_n$ solves the frequency-truncated analogue of the infinite normal form equation in \eqref{NFE}. 
Then, setting $w_n$ as $v_n = w_n + \D_n[w_n,w_n]$, the cancellation in Lemmas~\ref{LEM:wbd}~and~\ref{LEM:wint} extends to the frequency-truncated setting, using
\eqref{w} (with $\D$ replaced by truncated version $\D_n$).
We now define $z'_n := U(t) w_n = \Gc_{\dl, n} [ u_n ]$
which solves \eqref{eq:GNV_n_dens} and from Lemma~\ref{LEM:D} (now read for
$\Gc_{\dl,n}$), we get 
\begin{equation}\label{eq:z'n_a_priori}
z_n'\in C([0,T_{0,n}];H^{s'}(\R^2)),
\qquad
\|z_n'\|_{L^\infty_{T_{0,n}}H^s_x}<\tfrac{3}{1-\eps} R.
\end{equation}

\smallskip
\noindent\textit{Step~2.}\
From \eqref{eq:bounds_zn_nv}, we have
\begin{align}\label{bdzn1}
z_n\in C([0,T_{0,n}];H^{s'}(\R^2)),
\qquad
\| z_n \|_{L^\infty_{T_{0,n}} H^s_x} \le \tfrac{2}{1-\eps} R, 
\end{align}
which together with \eqref{eq:z'n_a_priori} verify the assumptions for the uniqueness statement in Lemma~\ref{lem:zn_nv}, guaranteeing that $z_n = z'_n$ in $C([0,T_{0,n}]; H^{s'}(\R^2) )$.
Then, combined with Lemma~\ref{LEM:D}, \eqref{bdzn1}, and \eqref{eq:bounds_zn_nv}, we get
\begin{align} 
\|u_n\|_{L^\infty_{T_{0,n}}H^s_x}
& \le (1+\eps)\,\|z_n'\|_{L^\infty_{T_{0,n}}H^s_x}
= (1+\eps)\,\|z_n\|_{L^\infty_{T_{0,n}}H^s_x}
\le \tfrac{2(1+\eps)}{1-\eps}R<3R,
\\
\|u_n\|_{L^\infty_{T_{0,n}}H^{s'}_x}
& \le (1+\eps)\,\|z_n'\|_{L^\infty_{T_{0,n}}H^{s'}_x}
= (1+\eps)\,\|z_n\|_{L^\infty_{T_{0,n}}H^{s'}_x}
\le \tfrac{2(1+\eps)}{1-\eps} \|u_0\|_{H^{s'}},
\label{eq:un_bound_step2}
\end{align}
provided that $0<\eps < \frac15$.
Therefore, from a continuity argument and the definition of $T_{0,n}$, we conclude that $T_{0,n} \ge  T_{\dl, R}$, thus all earlier properties hold up to time $T_{\dl,R}$ and (i)-(ii) in Lemma~\ref{lem:equiv_n_nv} follow.

\smallskip
\noindent\textit{Step~3.}\
For
$M:=\frac{2(1+\eps)}{1-\eps}\|u_0\|_{H^{s'}}$, the high-regularity well-posedness of
\cite{angelopolous,kazeykinamunoz1,kazeykinamunoz2} (or Appendix~\ref{sec:s-34}) gives
$\|u_n\|_{X^{s',b}_{T}}\le C(M)$ for some $T=T(M)>0$. 
Since the bound in \eqref{eq:un_bound_step2} holds on the time interval $[0,T_{\dl,R}]$, we can iterate the well-posedness argument in $X^{s',b}$ on time intervals $[\l T, (\l+1) T]$, $\l=0, \ldots, \lceil \frac{T_{\dl, R}}{T} \rceil -1 $, and the estimate in (iii) follows from a gluing argument. \qedhere

\end{proof}

\begin{proof}[Proof of Proposition~\textup{\ref{PRO:equiv}}]
Take $u_0$ as in \eqref{equiv-0}, $z_0 = \Gc_\dl [u_0]$,
$z\in X^{s',b}_{T_{\dl,R}}$ the solution to \eqref{z1}  with data $z_0$ from Proposition~\ref{prop:lwp_rnv}, and let
$u\in X^{s',b}_T$ to be the maximal \eqref{eq:NV} solution evolving
from $u_0$
(see \cite{angelopolous,kazeykinamunoz1,kazeykinamunoz2} or Appendix~\ref{sec:s-34}). 
Moreover, take $u_n,z_n \in  X^{s',b}_{T_{\dl,R}} \subset C([0,T_{\dl,R}]; H^{s'}(\R^2)$ as in 
Lemma~\ref{lem:equiv_n_nv}. 
The proof reduces to showing
\begin{align}
\label{eq:conv_zn}
z_n &\to z \quad\text{in}\quad C([0,T_{\dl,R}]; H^{s'-1}(\R^2)),
\\
\label{eq:conv_Gun}
\Gc_{\dl,n}[u_n] &\to \Gc_\dl[u]
\quad\text{in}\quad C([0,T\wedge T_{\dl,R}); H^{s'-1}(\R^2)).
\end{align}
In fact, \eqref{eq:conv_zn}-\eqref{eq:conv_Gun} together with \eqref{eq:equiv_n_nv} yield $\Gc_\dl[u(t)]=z(t)$ in $H^{s'-1}(\R^2)$ on
$[0,T\wedge T_{\dl,R})$.
Moreover, the $L^\infty_{T_{\dl,R}} H^{s'}_x$ control on $z$
provided by Proposition~\ref{prop:lwp_rnv} and
Remark~\ref{remark:persistence}, transferred through
Lemma~\ref{LEM:D} as
$\|u(t)\|_{H^{s'}}\le \frac{1}{1-\eps} \|z\|_{L^\infty_{T_{\dl,R}} H^{s'}_x}$ for $0 \le t \le T\land T_{\dl,R}$,
prevents $u$ from blowing up before $T_{\dl,R}$. Thus
$T\ge T_{\dl,R}$ and the identity extends to
$[0,T_{\dl,R}]$, showing \eqref{equiv-1}.

We first show \eqref{eq:conv_zn} by estimating the l.h.s. in $X^{s'-1,b}_{T_{\dl, R}}$. 
Let $0< T_* \le T_{\dl, R}$ and $M=2\|u_0\|_{H^{s'}}$.
Subtracting the Duhamel formulations of \eqref{z1} and
\eqref{eq:GNV_n_dens} splits $z(t)-z_n(t)$ into three pieces:
\begin{align}
z(t)-z_n(t)
&= U(t)(\Gc_\dl-\Gc_{\dl,n})[u_0]
+ \sum_{j=1}^\infty \sum_{\TT\in \Tr_j} \sum_{\l=1}^5 \int_0^t\! U(t-t')
\big[\Ncal_\l(\TT;z)-\Ncal_{\l,n}(\TT;z)\big]\,dt'
\\
&\quad
+ \sum_{j=1}^\infty \sum_{\TT\in \Tr_j} \sum_{\l=1}^5 \!\int_0^t\! U(t-t')
\big[\Ncal_{\l,n}(\TT;z)-\Ncal_{\l,n}(\TT;z_n)\big]\,dt'.
\label{zzn-diff0}
\end{align}
The first contribution in \eqref{zzn-diff0} is handled by noting that 
\begin{align}
\Gc_\dl[u_0]+\D[\Gc_\dl[u_0], 
\Gc_\dl[u_0]]
=
u_0
=\Gc_{\dl,n}[u_0]+\mathbf{P}_n\D[\Gc_{\dl,n}[u_0],\Gc_{\dl,n}[u_0]], 
\end{align}
with $\D$ in \eqref{eq:gauge}, from Lemma~\ref{LEM:D} and its proof, we get
\begin{align}
& 
\| U(t) (\Gc_\dl-\Gc_{\dl,n})[u_0] \|_{X^{s'-1,b}_{T_*}} 
= 
\|(\Gc_\dl-\Gc_{\dl,n})[u_0]\|_{H^{s'-1}}
\\
&
\le\| \P_{>n}\D[\Gc_\dl[u_0],\Gc_\dl[u_0]]\|_{H^{s'-1}}
+ 
\| \D [ \Gc_\dl[u_0],\Gc_\dl[u_0]] - \D [\Gc_{\dl, n}[u_0],\Gc_{\dl, n}[u_0]] \|_{H^{s'-1}}
\\
&
\le \tfrac1n \|\D[\Gc_\dl[u_0],\Gc_\dl[u_0]]\|_{H^{s'}}
+ 2 C \dl^{s+1} \tfrac{1}{1-\eps} R   \| (\Gc_\dl-\Gc_{\dl,n})[u_0]\|_{H^{s'-1}}
\end{align}
where $\P_{>n} = \Id - \P_n$ denoting the projection to frequencies $|\xi| > n$, $0<\eps < \frac12$, $0<\dl \ll R^{- 1/(1+s)}$, and $s'-1 \ge s > -1 $. Then, we can hide the second term on the l.h.s. by taking  $\dl$ potentially smaller, and the convergence to 0 as $n\to\infty$ follows.

For the last terms in \eqref{zzn-diff0}, 
Proposition~\ref{prop:bourgL2} (applied to the truncated multipliers,
which holds without change), Remark~\ref{REM:tree-size}, and \eqref{eq:bounds_zn_nv} give
control of their $X^{s'-1,b}_{T_*}$-norm by
\begin{align}
& T_*^{0+} C \dl^{-3}
\| z - z_n\|_{X^{s'-1,b}_{T_{*}}}
( \| z \|_{X^{s'-1,b}_{T_{\dl,R}}} + \| z_n\|_{X^{s'-1,b}_{T_{\dl,R}}} )
\sum_{j=1}^\infty 
( C  \dl^{s+1} R )^{j-1}
\\
&
\le
T_*^{0+} 2 C \dl^{-3 } M 
\| z - z_n\|_{X^{s'-1,b}_{T_{*}}}
\sum_{j=1}^\infty 
( C  \dl^{s+1} R )^{j-1}
\\
&
\le T_*^{0+} 2  C_{\dl, R}
M 
\| z-z_n\|_{X^{s'-1,b}_{T_{*}}} 
,
\end{align}
for $s'-1 \ge s$ and some constant $C_{\dl, R}>0$ increasing in $R$ and decreasing in $\dl$, which we can further bound by $\frac12 \|z-z_n\|_{X^{s'-1,b}_{T_{*}}}$ by taking $T_*=T_*(M, \dl, R)>0$ sufficiently small, and subsequently hide on the l.h.s. of the estimate for \eqref{zzn-diff0}.

For the second set of contributions in \eqref{zzn-diff0}, we note that $\Ncal_\l(\TT)-\Ncal_{\l,n}(\TT)$ decomposes into finitely many multilinear operators of the form $\Ncal_{\l}$ with some indicators $\ind_{|\xi_\bul| \le n}$ on some nodes $\bul\in\TT^0$ and at least an indicator $\ind_{|\xi_\bul|>n}$ for one $\bul\in \TT^0$. The former indicators are irrelevant in the application of Proposition~\ref{prop:bourgL2}, while the latter
give that $n < |\xi_\bul| \les \max_{\bul\in\TT^\infty} |\xi_\bul|$ which allows for a gain of $\frac1n$ at the cost of placing one of the $z$ factors in $X^{s',b}_{T_{\dl,R}}$, thus giving
\begin{align*}
& T_*^{0+} \tfrac1n CR \dl^{-3} \|z \|_{X^{s',b}_{T_{\dl,R}}}  \sum_{j=1}^\infty ( C \dl^{s+1} R)^{j-1}
\le T_*^{0+} \tfrac1n  M R C_{\dl,R} .
\end{align*}
Combining the above estimates, for $0< T_* = T_*(M,\dl,R) \le T_{\dl,R}$, we have
\begin{align*}
\|z - z_n \|_{L^\infty_{T_*} H^{s'-1}_x} 
\le \|z-z_n\|_{X^{s'-1,b}_{T_*}} 
\le \tfrac1n C \big\{ \|\D[\Gc_\dl[u_0],\Gc_\dl[u_0]]\|_{H^{s'-1}} +  T_*^{0+}  2 M R C_{\dl,R} \big\} \to 0 , 
\end{align*}
as $n\to \infty$, thus proving \eqref{eq:conv_zn} on $[0,T_*]$. The result extends to $[0,T_{\dl,R}]$ by iterating the above process on sub-intervals of size $T_*$.

\smallskip 
Lastly, we show \eqref{eq:conv_Gun}. 
Writing $\Gc_{\dl, n} [u_n] - \Gc_{\dl}[u] = (\Gc_{\dl, n} [u_n] - \Gc_{\dl,n}[u]) + (\Gc_{\dl, n} [u] - \Gc_{\dl}[u])$, we handle the second contribution by considering the von Neumann expansion for $\Gc_\dl$ in \eqref{vonneumann1} and the analogue for $\Gc_{\dl,n}$ with an indicator $\ind_{|\xi_\bul| \le n} $ for each $\bul\in \TT^0$, which we write as
\begin{align}
\Gc_{\dl, n} [u] - \Gc_{\dl}[u]
&
= \sum_{j=1}^\infty \sum_{\TT\in\Tr_j} \big[ \Ical_n(\TT, - i \mulB(\TT); u) - \Ical(\TT, - i \mulB(\TT); u) \big], 
\end{align}
noting that the difference in the summand decomposes into finitely many terms of the form $\Ical(\TT; - i\mulB(\TT), u)$ with 
one indicator $\ind_{|\xi_\star|>n}$ for one node $\star\in \TT^0$
and indicators $\ind_{|\xi_\bul|\le n}$ on some nodes $\bul\in\TT^0$, which can be uniformly controlled in $n$. This mimics the discussion for the second contributions in \eqref{zzn-diff0}. By a straightforward modification of the proof of \eqref{vonneumann2a}, using that $n \le |\xi_\star| \les \max_{\bul\in\TT^\infty}|\xi_\bul| $, we obtain
\begin{align}
\| \Gc_\dl [u] - \Gc_{\dl,n}[u] \|_{L^\infty_{T_*} H^{s'-1}_x} 
\le \eps \| \P_{>n} u \|_{L^\infty_{T_*} H^{s'-1}_x}
\le \tfrac1n \eps \| u\|_{L^\infty_{T_*} H^{s'}_x} \to 0, 
\end{align}
as $n\to\infty$, for $0< T_* \le T \land T_{\dl, R}$. 
We estimate the first term using \eqref{vonneumann2b} (which extends to $\Gc_{\dl, n}$ with uniform-in-$n$ constants) and \eqref{eq:bounds_un} to get
\begin{align}
\| \Gc_{\dl, n} [u_n] - \Gc_{\dl, n}[u]\|_{L^\infty_{T_*} H^{s'-1}_x} \le (1 + 3 M ) \| u_n - u \|_{L^\infty_{T_*} H^{s'-1}_x}, 
\end{align}
and thus it remains to show the convergence $u_n \to u$ in $C_{T_*} H^{s'-1}_x$ for $0< T_* \le T \land T_{\dl, R}$ to conclude \eqref{eq:conv_Gun}. 
Taking $T_* = T_*(M)>0$ coming from the well-posedness in \cite{angelopolous,kazeykinamunoz1,kazeykinamunoz2} and Appendix~\ref{sec:s-34}, 
using the relevant bilinear estimate \eqref{eq:bilin_s-34}, which also extend to the truncated equation \eqref{eq:NV_n_dens}, and Lemma~\ref{lem:equiv_n_nv}(iii), we have
\begin{align*}
& 
\|u_n-u\|_{L^\infty_{T_*} H^{s'-1}_{x}}
\le 
\|u_n-u\|_{X^{s'-1,b}_{T_*}}
\\
&
\quad 
\le  
\|u_n(0)-u(0)\|_{ H^{s'-1}}
+
 C T_*^\ta \| \P_{>n} B[u,u]\|_{X^{s'-1,b}_{T_*}}
+ C T_*^\ta \| \P_n (B[u_n,u_n] - B[u,u] )\|_{X^{s'-1,b}_{T_*}}
\\
&
\quad 
\le 
\|u_n(0)-u(0)\|_{ H^{s'-1}}
+
CT^\ta_* \tfrac1n \| u\|^2_{X^{s',b}_{T}}
+ 
CT^\ta_* (\|u_n\|_{X^{s'-1,b}_{T_*}} + \|u\|_{X^{s'-1,b}_{T_*}}) \|u_n-u\|_{X^{s'-1,b}_{T_*}}
\\
&
\quad 
\le 
\|u_n(0)-u(0)\|_{ H^{s'-1}}
+
CT^\ta_* \tfrac1n \| u\|^2_{X^{s',b}_{T}}
+ 
CT^\ta_* ( C(M) + \|u\|_{X^{s',b}_{T}}) \|u_n-u\|_{X^{s',b}_{T_*}}, 
\end{align*}
for some $\ta>0$, where $B[u,u]$ denotes the nonlinearity in \eqref{eq:NV}; see also \eqref{Bdef}. 
Then, we can take $0< T_* \le T \land T_{\dl, R}$ sufficiently small, depending on $M$ and $\|u\|_{X^{s',b}_T}$, in order to hide the last contribution on the l.h.s., which together with $\| u_n(0) - u(0)\|_{H^{s'-1}}=0$ gives the convergence as $n\to\infty$ on $[0,T_*]$. The convergence on $[0,T\land T_{\dl,R}]$ then follows from iterating the above argument on sub-intervals of length $T_*$, completing the proof. \qedhere

\end{proof}

\appendix

\section{Improved local well-posedness for \texorpdfstring{$s>-\frac34$}{s>-3/4}}\label{sec:s-34}

In this short appendix, we extend the well-posedness result of Adams and Gr\"unrock~\cite{adamsgrunrock} for \eqref{eq:NV} with $E=0$ in $H^s(\R^2), s>-\frac34,$ to the nonzero energy case.
In particular, in this regularity range, we show that the solutions in Theorem \ref{thm:main_intro} also satisfy $u\in X^{s,b}_T$.

\begin{theorem}\label{thm:s-34}
Fix $E\in \R$ and $s>-\frac34$. Given $u_0\in H^s(\R^2)$, there exist $T=T(\|u_0\|_{H^s}, |E|)>0$ and a unique strong integral solution $u\in X^{s,b}_T$ to \eqref{eq:NV} with initial data $u_0$. Moreover, the data-to-solution map is analytic.
\end{theorem}

As before, using the scaling invariance \eqref{scaling}, it suffices to prove Theorem \ref{thm:s-34} for $|E|\ll1$. Moreover, 
Theorem~\ref{thm:s-34} follows from a standard contraction mapping argument in $X^{s,b}_T$ (see \cite{bourg1, bourg2} and the monographs \cite{linaresponce, tao_book}), once we establish a bilinear estimate for the operator $B$ defining the nonlinearity of \eqref{eq:NV}:
\begin{align}\label{Bdef}
B[u_1,u_2]^\wedge(\xi):=\int_{\xi=\xi_1+\xi_2} K(\xi,\xi_1,\xi_2)\Phi_0(\xi_1,\xi_2) \ft{u}_1(\xi_1) \ft{u}_2(\xi_2)d\xi_1, 
\end{align}
recalling $K, \Phi_0$ in \eqref{Kmul} and \eqref{Phi0}.

\begin{proposition}\label{prop:bilin_s-34}
Fix $|E|\ll 1$ and $s>-\frac34$. Then, for $B$ as in \eqref{Bdef}, we have
\begin{equation}\label{eq:bilin_s-34}
\left\| B[u_1,u_2]\right\|_{X^{s,-\frac12+}}\lesssim \|u_1\|_{X^{s,\frac12+}}\|u_2\|_{X^{s,\frac12+}}.
\end{equation}
\end{proposition}

\begin{proof}
First, we decompose $B$ as
\begin{align*}
B[u_1,u_2]^\wedge(\xi)&=\int_{\xi=\xi_1+\xi_2} K\Phi_0\mathbbm{1}_{A^c}\ft{u}_1(\xi_1)\ft{u}_2(\xi_2)d\xi_1 \\&\qquad+ \int_{\xi=\xi_1+\xi_2} K\Phi_0\mathbbm{1}_A \ft{u}_1(\xi_1) \ft{u}_2(\xi_2)d\xi_1, 
\end{align*}
where we used \eqref{Phi}, 
recalling $A,A^c$ in \eqref{A},  where we choose $\dl=1$. 
The bilinear estimate for the last term has already been shown in Proposition \ref{prop:bourgL2} for $s>-1$ (it corresponds to the Type $\I$ term, see Lemma~\ref{LEM:wint} and \eqref{eq:multiplicadores}). For the remaining term, we use Lemma \ref{lem:fre_bi}(i) (together with Lemma \ref{lem:eta1}) and prove the necessary frequency-restricted estimates.

To simplify the exposition, we take $-\frac34<s\le0$. From Lemma~\ref{lem:fre_bi}~(i),  it suffices to consider $1 < M \les |\al|$ and estimate the following for all distinct $k,\l \in \{\emptyset, 1,2\}$:
\begin{align}\label{app:fre}
\sup_{\xi_k} \int_{\xi = \xi_1+ \xi_2 } \frac{|K \Phi|^2 \jb{\xi}^{2s}}{\jb{\xi_1}^{2s} \jb{\xi_2}^{2s} } \ind_{A^c} \ind_{|\Phi - \al| < M} d\xi_\l . 
\end{align}
Using Lemma~\ref{lem:pointwise_multipliers} and assuming that $|\xi_1| \ge |\xi_2|$, by symmetry, we bound the multiplier above as follows
\begin{align}
\frac{|K \Phi|^2 \jb{\xi}^{2s}}{\jb{\xi_1}^{2s-} \jb{\xi_2}^{2s} } \ind_{|\Phi - \al| < M}
\les 
|\al| \frac{\jb{\xi}^{2s} |\xi| }{\jb{\xi_1}^{2s+1-} \jb{\xi_2}^{2s+1} } \ind_{|\Phi - \al| < M}. 
\end{align}
Moreover, we consider only the worst-case scenario $|\xi_1| \ge |\xi_2| \ge |\xi|$. 

We first fix $\xi_2$, i.e., $k=2$ in \eqref{app:fre}. 
If $\xi_1 \not\simeq \pm \xi$, by Lemma~\ref{lem:prop_resonance}~(i), we have
\begin{align}
\eqref{app:fre} \vert_{k=2} \les |\al| \int \frac{\jb{\xi}^{2s+1}}{ \jb{\xi_1}^{2s+3-} \jb{\xi_2}^{2s+1}} \ind_{|\phi-\al|<M} d\phi d\zeta \les |\al | M \int_\R  \jb{\zeta}^{-2s-3+} d\zeta \les |\al | M , 
\end{align}
for $s>-1$. 
Otherwise, $\xi_1 \simeq \ld \xi$ for $\ld \in\{\pm 1\}$. From the assumption that $|\xi_1| \ge |\xi_2|$, we have $\xi_1 \simeq -\xi$ and Lemma~\ref{lem:prop_resonance}~(ii) applies to give
\begin{align}
\eqref{app:fre} \vert_{k=2} \les |\al| \sup_{\xi_2} \int_{|q| < 1} \jb{\xi_2}^{-2s+1+} \ind_{| |\xi_2|^3 \Re q^2 -\wt \al|<M} dq \les 
|\al| M \sup_{\xi_2} \jb{\xi_2}^{-2s-2+} \les |\al| M.
\end{align}

For $\xi_1$ fixed, we consider the cases $\xi_2 \not\simeq \pm \xi$ or $\xi_2 \not\simeq \ld \xi$ and proceed as above, swapping the roles of $\xi_1$ and $\xi_2$. 

Lastly, we consider $\xi$ fixed. If $\xi_1 \not\simeq \pm \xi_2 $, we use Lemma~\ref{lem:prop_resonance}~(i) and perform the change of variables $\xi_2 \mapsto (\phi, \zeta_2) = (\Phi, \Im \xi_2)$, to find
\begin{align}
\eqref{app:fre} \vert_{k=\emptyset} \les |\al| \int \jb{\zeta_2}^{-4s-4+} \ind_{|\phi-\al| < M } d\phi d\zeta_2 \les |\al| M , 
\end{align}
for $s>-\frac34$. 
If $\xi_1 \simeq \xi_2$, then $|\xi_1|\sim|\xi_2|\sim|\xi|$ and we apply Lemma~\ref{lem:prop_resonance}:
\begin{align}
\eqref{app:fre} \vert_{k=\emptyset} \les |\al| \sup_\xi \int_{|q_2|<1} \jb{\xi}^{-2s+1+} \ind_{||\xi|^3 Re q^2-\al| < M } dq \les |\al| M \sup_\xi \jb{\xi}^{-2s-2+} \les |\al | M .
\end{align}
Lastly, if $\xi_ 1\simeq - \xi_2$, we have $|\xi | = |\xi_1+ \xi_2| \ll |\xi_1| \sim |\xi_2|$, and we instead verify the condition in \eqref{eq:frequad_2d} instead of \eqref{eq:frequad}; see Remark~\ref{remark:frebi}. From the frequency relations, we have
$$
|\nabla_{\xi_1}\Phi|\gtrsim  |\xi|\big(|\xi_1|-\tfrac{|E|}{|\xi_1|} \big) \gtrsim |\xi||\xi_1|.
$$
Writing $\xi_j = (\eta_j, \zeta_j)$, $j\in\{\emptyset, 1,2\}$, we can suppose that $|\partial_{\eta_1} \Phi| \ges |\xi| |\xi_1|$. Then, fixing $\zeta, \zeta_1,\zeta_2$, we apply the change of variables $\eta_1 \mapsto \phi = \Phi(\eta_1)$:

\begin{align*}
|\al|  \sup_{\eta} \int \frac{\jap{\xi}^{2s}|\xi|}{\jap{\xi_1}^{4s+2-}}\ind_{|\Phi-\alpha|<M} d\eta_1 \lesssim   |\al| \jap{\zeta}^{2s} \sup_{\eta} \int \frac{1}{\jap{\xi_1}^{4s+3-}}\ind_{|\phi-\alpha|<M} d\phi\lesssim \jap{\zeta}^{-2s-3+} |\al| M, 
\end{align*}
for $s>-\frac34$. Then, the r.h.s. is integrable in $\zeta$ and the estimate follows from Remark \ref{remark:frebi}, concluding the proof of \eqref{eq:bilin_s-34}. \qedhere

\end{proof}

We now turn to the failure of the bilinear estimate for $s\le -\frac78$. The proof follows from a straightforward adaptation of the classical counterexample of Kenig, Ponce and Vega for the Korteweg-de Vries equation \cite[Theorem 1.3]{kpv_bilin}.

\begin{proposition}\label{PRO:ill}
Take $E=0$, $b\in\R$, and $-\frac12 < b'  \le 0$. Then, for $B$ as in \eqref{Bdef}, the estimate
\begin{equation}\label{eq:bilin_s-78}
\left\| B[u_1,u_2]\right\|_{X^{s,b'}}\lesssim \|u_1\|_{X^{s,b}}\|u_2\|_{X^{s,b}}.
\end{equation}
fails for $s \le - \frac78$.
\end{proposition}
\begin{proof}

Given $N\gg 1$, consider the sets
\begin{align*}
R_1& =\big\{(\tau,\xi)\in\R^3: |\xi-N|<N^{-\frac12},\ |\tau-\phi_0(N)-\nabla\phi_0(N)\cdot (\xi-N)|<1\big\}
, 
\\
R_2& =\big\{(\tau,\xi)\in\R^3: |\xi+N|<2N^{-\frac12},\ |\tau-\phi_0(-N)-\nabla\phi_0(-N)\cdot (\xi+N)|<2\big\}
, 
\\
R & =\big\{(\tau,\xi)=(\tau, \eta, \zeta)\in\R^3: |\xi|<N^{-\frac12}, \ \eta<-\tfrac{1}{2}N^{-\frac12},\ |\tau-\nabla\phi_0(N)\cdot \xi|<1\big\},
\end{align*}
whose measure is comparable to $N^{-1}$. If $(\tau_1,\xi_1)\in R_1$ and $(\tau,\xi)\in R$, then $(\tau_2,\xi_2)=(\tau-\tau_1, \xi-\xi_1)\in R_2$: indeed,
\begin{align*}
|\xi_2+N|& \le |\xi| + |\xi_1-N| < 2N^{-\frac12}
\\
|\tau_2-\phi_0(-N)-\nabla\phi_0(-N)\cdot (\xi+N)| &= |\tau-\tau_1+\phi_0(N)-\nabla\phi_0(N)(\xi-(\xi_1-N))| \\&\le |\tau -\nabla\phi_0(N)\xi| + |\tau_1-\phi_0(N)-\nabla\phi_0(N)(\xi_1-N)|<2.
\end{align*}
Similarly, if $(\tau_1,\xi_1)\in -R_1$ and $(\tau, \xi)\in R$, then $(\tau_2,\xi_2)\in -R_2$.

Define $\sigma_j=\tau_j-\phi_0(\xi_j)$, $j=\emptyset, 1, 2$. Then, for $(\tau_1,\xi_1)\in R_1$, 
\begin{align*}
|\sigma_1|&\le 1 + |\phi_0(\xi_1)-\phi_0(N)-\nabla\phi_0(N)\cdot (\xi_1-N)|\\&\le 1 + |D^2\phi_0(N)||\xi_1-N|^2 + C|\xi_1-N|^3 \lesssim 1, 
\end{align*}
since $|D^2\phi_0(N)|\sim N$ and $|\xi_1-N|<N^{-\frac12}$. By symmetry, $|\sigma_2|\lesssim 1$ for $(\tau_2, \xi_2)\in R_2$.
The same computation yields $|\s_1|, |\s_2|\les 1$ for $(\tau_1,\xi_1) \in -R_1$ and $(\tau_2,\xi_2) \in -R_2$. 

\medskip

Set
$$
\widehat{u_1}=\ind_{R_1}+\ind_{-R_1},\quad \widehat{u_2}=\ind_{R_2}+\ind_{-R_2}, \quad \widehat{u}=\ind_R.
$$
If \eqref{eq:bilin_s-78} holds, then
\begin{align*}
\I: =\int_{\substack{\xi_2=\xi-\xi_1\\\tau_2=\tau-\tau_1}} \frac{\jap{\sigma}^{b'}\jap{\xi}^sK\Phi_0}{\jap{\xi_1}^s\jap{\xi_2}^s\jap{\sigma_1}^b\jap{\sigma_2}^b}\widehat{u_1}(\tau_1,\xi_1)\widehat{u_2}(\tau_2,\xi_2)\widehat{u}(\tau,\xi)& d\tau_1 d\tau d\xi_1 d\xi \les \prod_{j=\emptyset, 1,2} \| \ft{u}_j\|_{L^2_{\tau, \xi}}, 
\end{align*}
The r.h.s. is comparable to $N^{-\frac32}$. For the l.h.s.,  first observe that, due to the convolution relation and $(\tau,\xi)\in R$, the only nonzero contributions are $(\tau_1, \xi_1)\in R_1$, $(\tau_2, \xi_2)\in R_2$ or $(\tau_1, \xi_1)\in -R_1$, $(\tau_2, \xi_2)\in -R_2$.  In any case, a simple computation shows that
$$
K\Phi_0\gtrsim N^{-\frac12}, \quad |\Phi_0|\sim N^{\frac{3}{2}},\quad |\sigma|\le |\Phi_0|+|\sigma_1|+|\sigma_2| \lesssim N^{\frac{3}{2}}.
$$
Then
\begin{align*}
\I \gtrsim &\int_{R\times R_1} \frac{N^{\frac{3b'}{2}}\cdot N^{-\frac12}}{N^{2s}} d\tau d\xi d\tau_1 d\xi_1 \sim N^{\frac{3b'}{2}-2s-\frac{5}{2}}. 
\end{align*}
We conclude that, for all $N\gg 1$, if $b' > -\frac12$, then
$$
N^{\frac{3b'}{2}-2s-\frac{5}{2}}  \ge   N^{-\frac{13}{4} - 2s } \gg  N^{-\frac32}
$$
 where the last inequality implies $s \le - \frac78$, completing the proof. 
\end{proof}

\bibliographystyle{plain}
\bibliography{biblio}

\end{document}